\documentclass[11pt,reqno]{amsart}
\usepackage{amsmath}
\usepackage{amsfonts}
\usepackage{amssymb}
\usepackage{amsgen}
\usepackage{mathrsfs}
\usepackage{cite}
\usepackage[all]{xy}

\newcommand{\supp}{\mathop{\mathrm{supp}}}
\newcommand{\rad}{\mathop{\mathrm{rad}}}
\newcommand{\modu}{\mathrm{mod}\text{-}}
\newcommand{\bil}{\mathop{\mathrm{bil}}\nolimits}
\newcommand{\Res}{\mathop{\mathrm{Res}}\nolimits}
\newcommand{\Ind}{\mathop{\mathrm{Ind}}\nolimits}
\newcommand{\Coind}{\mathop{\mathrm{Coind}}\nolimits}
\newcommand{\rk}{\mathop{\mathrm{rk}}\nolimits}
\newcommand{\J}{\mathrel{\mathscr J}} 
\newcommand{\R}{\mathrel{\mathscr R}} 
\newcommand{\eL}{\mathrel{\mathscr L}} 

\newcommand{\OO}{\mathcal O}
\newcommand{\res}{\mathop{\mathrm{res}}\nolimits}
\newcommand{\ind}{\mathop{\mathrm{ind}}\nolimits}
\newcommand{\coind}{\mathop{\mathrm{coind}}\nolimits}
\newcommand{\pv}[1]{\mathbf {#1}}
\newcommand{\inv}{^{-1}}
\newcommand{\p}{\varphi}

\newcommand{\ov}[1]{\ensuremath{\overline {#1}}}
\newcommand{\til}[1]{\ensuremath{\widetilde {#1}}}

\newcommand{\Irr}{\mathop{\mathrm{Irr}}\nolimits}
\newcommand{\Sub}{\mathop{\mathrm{Sub}}\nolimits}
\newcommand{\Pos}{\mathop{\mathrm{Pos}}\nolimits}
\newcommand{\Cong}{\mathop{\mathrm{Cong}}\nolimits}
\newcommand{\RLM}{\mathop{\mathsf{RLM}}\nolimits}
\newcommand{\soc}[1]{\mathrm{Soc}(#1)}

\newcommand{\Thmname}{Theorem}
\newcommand{\Propname}{Proposition}
\newcommand{\Lemmaname}{Lemma}
\newcommand{\Definitionname}{Definition}
\newtheorem{Thm}{\Thmname}[section]
\newtheorem{Prop}[Thm]{\Propname}
\newtheorem{Lemma}[Thm]{\Lemmaname}
{\theoremstyle{definition}
}
{\theoremstyle{remark}
\newtheorem{Rmk}[Thm]{Remark}}
\newtheorem{Cor}[Thm]{Corollary}
{\theoremstyle{remark}
}

\newtheorem*{Lemma*}{Lemma}

\numberwithin{equation}{section}

\title[A Theory of Transformation Monoids]{A Theory of Transformation Monoids: Combinatorics and Representation Theory}

\author{Benjamin Steinberg}
\address{School of Mathematics and Statistics\\
Carleton University \\
1125 Colonel By Drive\\
Ottawa, Ontario  K1S 5B6 \\
Canada}
\thanks{The author was supported in part by NSERC}
\email{bsteinbg@math.carleton.ca}
\date{April 17, 2010}

\subjclass[2010]{20M20, 20M30, 20M35}

\keywords{transformation monoids, orbitals, semigroup representation theory, Markov chains, synchronizing automata}

\begin{document}
\begin{abstract}
The aim of this paper is to develop a theory of finite transformation monoids and in particular to study primitive transformation monoids.   We introduce the notion of orbitals and orbital digraphs for transformation monoids and prove a monoid version of D.~Higman's celebrated theorem characterizing primitivity in terms of connectedness of orbital digraphs.

A thorough study of the module (or representation) associated to a transformation monoid is initiated.  In particular, we compute the projective cover of the transformation module over a field of characteristic zero in the case of a transitive transformation or partial transformation monoid.  Applications of probability theory and Markov chains to transformation monoids are also considered and an ergodic theorem is proved in this context.  In particular, we obtain a generalization of a lemma of P.~Neumann, from the theory of synchronizing groups, concerning the partition associated to a transformation of minimal rank.
\end{abstract}
\maketitle
\tableofcontents

\section{Introduction}
The principal task here is to initiate a theory of finite transformation monoids that is similar in spirit to the theory of finite permutation groups that can be found, for example, in~\cite{dixonbook,cameron}.  I say similar in spirit because attempting to study transformation monoids by analogy with permutation groups is like trying to study finite dimensional algebras by analogy with semisimple algebras.  In fact, the analogy between finite transformation monoids and finite dimensional algebras is quite apt, as the theory will show.  In particular, an analogue of Green's theory~\cite[Chapter 6]{Greenpoly} of induction and restriction functors relating an algebra $A$ with algebras of the form $eAe$ with $e$ idempotent plays a key role in this paper, whereas there is no such theory in permutation groups as there is but one idempotent.

There are many worthy books that touch upon --- or even focus on --- transformation monoids~\cite{CP,Higginsbook,howiebook,GM,Lipscomb}, as well as a vast number of research articles on the subject.  But most papers in the literature focus on specific transformation monoids (such as the full transformation monoid, the symmetric inverse monoid, the monoid of order preserving transformations, the monoid of all partial transformations, etc.) and on combinatorial issues, e.g., generalizations of cycle notation, computation of the submonoid generated by the idempotents~\cite{Howie}, computation of generators and relations, computation of Green's relations, construction of maximal submonoids satisfying certain properties, etc.

The only existing theory of finite transformation and partial transformation monoids as a general object is the Krohn-Rhodes wreath product decomposition theory~\cite{PDT,KRannals,Arbib}, whose foundations were laid out in the book of Eilenberg~\cite{Eilenberg}. See also~\cite{qtheor} for a modern presentation of the Krohn-Rhodes theory, but with a focus on abstract rather than transformation semigroups.

The Krohn-Rhodes approach is very powerful, and in particular has been very successful in dealing with problems in automata theory, especially those involving classes of languages.  However, the philosophy of Krohn-Rhodes is that the task of classifying monoids (or transformation monoids) up to isomorphism is hopeless and not worthwhile.  Instead, one uses a varietal approach~\cite{Eilenberg} similar in spirit to the theory of varieties of groups~\cite{Neumann}.  But there are some natural problems in automata theory where one really has to stick with a given transformation monoid and cannot perform the kind of decompositions underlying the Krohn-Rhodes theory.  One such problem is the \v{C}ern\'y conjecture, which has a vast literature~\cite{Pincerny,pincernyconjecture,synchgroups,dubuc,cerny,volkovc1,rystsov1,rystsov2,AMSV,Trahtman,traht2,volkovc2,Kari,volkovc3,rystcom,rystrank,mycerny,Karicounter,VolkovLata,PerrinBeal,strongtrans,strongtrans2,mortality,beal,Salomcerny,averaging}.
In the language of transformation monoids, it says that if $X$ is a set of maps on $n$ letters such that some product of elements of $X$ is a constant map, then there is a product of length at most $(n-1)^2$ that is a constant map.  The best known upper bound is cubic~\cite{twocomb}, whereas it is known that one cannot do better than $(n-1)^2$~\cite{cerny}.

Markov chains can often be fruitfully studied via random mappings:  one has a transformation monoid $M$ on the state set $\Omega$ and a probability $P$ on $M$.  One randomly chooses an element of $M$ according to $P$ and has it act on $\Omega$.  A theory of transformation monoids, in particular of the associated matrix representation, can then be used to analyze the Markov chain.  This approach has been adopted with great success by Bidigare, Hanlon and Rockmore~\cite{BHR}, Diaconis and Brown~\cite{DiaconisBrown1,Brown1,Brown2} and Bj\"orner~\cite{bjorner1,bjorner2}; see also my papers~\cite{mobius1,mobius2}.  This is another situation to which the Krohn-Rhodes theory does not seem to apply.

This paper began as an attempt to systematize and develop some of the ideas that have been used by various authors while working on the \v{C}ern\'y conjecture.  The end result is the beginnings of a theory of transformation monoids.  My hope is that the theory initiated here will lead toward some progress on the \v{C}ern\'y conjecture.   However, it is also my intent to interest combinatorialists, group theorists and representation theorists in transformation monoids and convince them that there is quite a bit of structure there.  For this reason I have done my best not to assume any background knowledge in semigroup theory and to avoid usage of certain semigroup theoretic notions and results, such as Green's relations~\cite{Green} and Rees's theorem~\cite{CP}, that are not known to the general public.  In particular, many standard results in semigroup theory are proved here in a novel way, often using transformation monoid ideas and in particular an analogue of Schur's lemma.

The first part of the paper is intended to systemize the foundations of the theory of transformation monoids.  A certain amount of what is here should be considered folklore, although probably some bits are new.  I have tried to indicate what I believe to be folklore or at least known to the \textit{cognoscenti}.  In particular, some of Sections~\ref{sthree} and~\ref{sfour} can be viewed as a specialization of Sch\"utzenberger's theory of unambiguous matrix monoids~\cite{berstelperrinreutenauer}. The main new part here is the generalization of Green's theory~\cite{Greenpoly} from the context of modules to transformation monoids.  A generalization of Green's results to semirings, with applications to the representation theory of finite semigroups over semirings, can be found in~\cite{ZurJohnBen}.

The second part of the paper is a first step in the program of understanding primitive transformation monoids.  In part, they can be understood in terms of primitive groups in much the same way that irreducible representations of monoids can be understood in terms of irreducible representations of groups via Green's theory~\cite{Greenpoly,myirreps} and the theory of Munn and Ponizovsky~\cite[Chapter 5]{CP}.  The tools of orbitals and orbital digraphs are introduced, generalizing the classical theory from permutation groups~\cite{dixonbook,cameron}.

The third part of the paper commences a detailed study of the modules associated to a transformation monoid.  In particular, the projective cover of the transformation module is computed for the case of a transitive action by partial or total transformations.  The paper ends with applications of Markov chains to the study of transformation semigroups.

\section{Actions of monoids on sets}
Before turning to transformation monoids, i.e., monoids acting faithfully on sets, we must deal with some ``abstract nonsense'' type preliminaries concerning monoid actions on sets and formalize notation and terminology.

\subsection{$M$-sets}
Fix a monoid $M$.  A (right) \emph{action} of $M$ on a set $\Omega$ is, as usual, a map $\Omega\times M\to \Omega$, written $(\alpha,m)\mapsto \alpha m$, satisfying, for all $\alpha\in \Omega$, $m,n\in M$,
\begin{enumerate}
\item $\alpha 1=\alpha$;
\item $(\alpha m)n = \alpha (mn)$.
\end{enumerate}
Equivalently, an action is a homomorphism $M\to T_{\Omega}$, where $T_{\Omega}$ is the monoid of all self-maps of $\Omega$ acting on the right.  In this case, we say that $\Omega$ is an \emph{$M$-set}. The action is \emph{faithful} if the corresponding morphism is injective.  Strictly speaking, there is a unique action of $M$ on the empty set, but in this paper we tacitly assume that we are dealing only with actions on non-empty sets.

A \emph{morphism} $f\colon \Omega\to \Lambda$ of $M$-sets is a map such that $f(\alpha m)=f(\alpha)m$ for all $\alpha\in \Omega$ and $m\in M$.  The set of morphisms from $\Omega$ to $\Lambda$ is denoted $\hom_M(\Omega,\Lambda)$. The category of right $M$-sets will be denoted $\pv{Set}^{M^{\mathrm {op}}}$ following category theoretic notation for presheaf categories~\cite{Mac-CWM}.

The $M$-set obtained by considering the right action of $M$ on itself by right multiplication is called the \emph{regular} $M$-set. It is a special case of a free $M$-set.  An $M$-set $\Omega$ is \emph{free} on a set $X$ if there is a map $\iota\colon X\to M$ so that given a function $g\colon X\to \Lambda$ with $\Lambda$ an $M$-set, there is a unique morphism of $M$-sets $f\colon \Omega\to \Lambda$ such that
\[\xymatrix{X\ar[r]^{\iota}\ar[rd]_g&\Omega\ar@{..>}[d]^f\\ & \Lambda}\]
commutes.  The free $M$-set on $X$ exists and can explicitly be realized as $X\times M$ where the action is given by $(x,m')m = (x,m'm)$ and the morphism $\iota$ is $x\mapsto (x,1)$.  The functor $X\mapsto X\times M$ from $\pv {Set}$ to $\pv{Set}^{M^{\mathrm {op}}}$ is left adjoint to the forgetful functor.  In concrete terms, an $M$-set $\Omega$ is free on a subset $X\subseteq \Omega$ if and only if, for all $\alpha\in \Omega$, there exists a unique $x\in X$ and $m\in M$ such that $\alpha=xm$.  We call $X$ a \emph{basis} for the $M$-set $\Omega$.  Note that if $M$ is a group, then $\Omega$ is free if and only if $M$ acts \emph{freely} on $\Omega$, i.e., $\alpha m=\alpha$, for some $\alpha\in \Omega$, implies $m=1$.  In this case, any transversal to the $M$-orbits is a basis.

Group actions are to undirected graphs as monoid actions are to directed graphs (digraphs).  Just as a digraph has both weak components and strong components, the same is true for monoid actions.  Let $\Omega$ be an $M$-set.  A non-empty subset $\Delta$ is \emph{$M$-invariant} if $\Delta M\subseteq M$; we do not consider the empty set as an $M$-invariant subset.  An $M$-invariant subset of the form $\alpha M$ is called \emph{cyclic}.  The cyclic sub-$M$-sets form a poset $\Pos(\Omega)$ with respect to inclusion.  The assignment $\Omega\to \Pos(\Omega)$ is a functor $\pv{Set}^{M^{\mathrm {op}}}\to \pv{Poset}$.  A cyclic subset will be called \emph{minimal} if it is minimal with respect to inclusion.

Associated to $\Pos(\Omega)$ is a preorder on $\Omega$ given by $\alpha\leq_{\Omega} \beta$ if and only if $\alpha M\subseteq \beta M$.  If $\Omega$ is clear from the context, we drop the subscript and simply write $\leq$.  From this preorder arise two naturally defined equivalence relations: the symmetric-transitive closure $\simeq$ of $\leq$ and the intersection $\sim$ of $\leq$ and $\geq$.  More precisely, $\alpha\simeq \beta$ if and only if there is a sequence $\alpha=\omega_0,\omega_1,\ldots, \omega_n=\beta$ of elements of $\Omega$ such that, for each $0\leq i\leq n-1$, either $\omega_i\leq \omega_{i+1}$ or $\omega_{i+1}\leq \omega_i$.  On the other hand, $\alpha\sim \beta$ if and only if $\alpha\leq \beta$ and $\beta\leq \alpha$, that is, $\alpha M=\beta M$.  The equivalence classes of $\simeq$ shall be called \emph{weak orbits}, whereas the equivalence classes of $\sim$ shall be called \emph{strong orbits}.  These correspond to the weak and strong components of a digraph.  If $M$ is a group, then both notions coincide with the usual notion of an orbit.

Notice that weak orbits are $M$-invariant, whereas a strong orbit is $M$-invariant if and only if it is a minimal cyclic subset $\alpha M$.   The action of $M$ will be called \emph{weakly transitive} if it has a unique weak orbit and shall be called \emph{transitive}, or \emph{strongly transitive} for emphasis, if it has a unique strong orbit. Observe that $M$ is transitive on $\Omega$ if and only if there are no proper $M$-invariant subsets of $\Omega$. Thus transitive $M$-sets can be thought of as analogues of irreducible representations; on the other hand weakly transitive $M$-sets are the analogues of indecomposable representations since it is easy to see that the action of $M$ on $\Omega$ is weakly transitive if and only if $\Omega$ is not the coproduct (disjoint union) of two proper $M$-invariant subsets.  The regular $M$-set is weakly transitive, but if $M$ is finite then it is transitive if and only if $M$ is a group.
The weak orbit of an element $\alpha\in \Omega$ will be denoted $\OO_w(\alpha)$ and the strong orbit $\OO_s(\alpha)$.  The set of weak orbits will be denoted $\pi_0(\Omega)$ (in analogy with connected components of graphs; and in any event this designation can be made precise in the topos theoretic sense) and the set of strong orbits shall be denoted $\Omega/M$. Note that $\Omega/M$ is naturally a poset isomorphic to $\Pos(\Omega)$ via the bijection $\OO_s(\alpha)\mapsto \alpha M$.  Also note that $\pi_0(\Omega)$ is in bijection with $\pi_0(\Pos(\Omega))$ where we recall that if $P$ is a poset, then the set $\pi_0(P)$ of connected components of $P$ is the set of equivalence classes of the symmetric-transitive closure of the partial order (i.e., the set of connected components of the Hasse diagram of $P$).

We shall also have need to consider $M$-sets with zero.  An element $\alpha\in \Omega$ is called a \emph{sink} if $\alpha M=\{\alpha\}$.  An \emph{$M$-set with zero}, or \emph{pointed $M$-set}, is a pair $(\Omega,0)$ where $\Omega$ is an $M$-set and $0\in M$ is a distinguished sink\footnote{This usage of the term ``pointed transformation monoid'' differs from that of~\cite{qtheor}.}.  An $M$-set with zero $(\Omega,0)$ is called \emph{$0$-transitive} if $\alpha M=\Omega$ for all $\alpha\neq 0$.  Notice that an $M$-set with zero is the same thing as an action of $M$ by partial transformations (just remove or adjoin the zero) and that $0$-transitive actions correspond to transitive actions by partial functions. Morphisms of $M$-sets with zero must preserve the zero and, in particular, in this context $M$-invariant subsets are assumed to contain the zero.  The category of $M$-sets with zero will be denoted $\pv{Set}_*^{M^{\mathrm{op}}}$ as it is the category of all contravariant functors from $M$ to the category of pointed sets.

\begin{Prop}\label{uniquesink}
Suppose that $\Omega$ is a $0$-transitive $M$-set.  Then $0$ is the unique sink of $\Omega$.
\end{Prop}
\begin{proof}
Suppose that $\alpha\neq 0$.  Then $0\in \Omega=\alpha M$ shows that $\alpha$ is not a sink.
\end{proof}

A strong orbit $\OO$ of $M$ on $\Omega$ is called \emph{minimal} if it is minimal in the poset $\Omega/M$, or equivalently the cyclic poset $\omega M$ is minimal for $\omega\in \OO$.   The union of all minimal strong orbits of $M$ on $\Omega$ is $M$-invariant and is called the \emph{socle} of $\Omega$, denoted $\soc \Omega$. If $M$ is a group, then $\soc \Omega=\Omega$. The case that $\Omega=\soc \Omega$ is analogous to that of a completely reducible representation: one has that $\Omega$ is a coproduct of transitive $M$-sets.  If $\Omega$ is an $M$-set with zero, then a minimal non-zero strong orbit is called \emph{$0$-minimal}.  In this setting we define the socle to be the union of all the $0$-minimal strong orbits together with zero; again it is an $M$-invariant subset.

A \emph{congruence} or \emph{system of imprimitivity} on an $M$-set $\Omega$ is an equivalence relation $\equiv$ such that $\alpha \equiv \beta$ implies $\alpha m\equiv \beta m$ for all $\alpha,\beta\in \Omega$ and $m\in M$.  In this case, the quotient $\Omega/{\equiv}$ becomes an $M$-set in the natural way and the quotient map $\Omega\to \Omega/{\equiv}$ is a morphism.  The standard isomorphism theorem holds in this context.  If $\Delta\subseteq \Omega$ is $M$-invariant, then one can define a congruence $\equiv_{\Delta}$ by putting $\alpha\equiv_{\Delta} \beta$ if $\alpha=\beta$ or $\alpha,\beta\in \Delta$. In other words, the congruence $\equiv_{\Delta}$ crushes $\Delta$ to a point.  The quotient $M$-set is denoted $\Omega/\Delta$.  The class of $\Delta$, often denoted by $0$, is a sink and it is more natural to view $\Omega/\Delta$ as an $M$-set with zero.  The reader should verify that if
\begin{equation}\label{series}
\Omega = \Omega_0\supset\Omega_1\supset \Omega_2\supset\cdots \supset \Omega_k
\end{equation}
is an unrefinable chain of $M$-invariant subsets, then the successive quotients $\Omega_i/\Omega_{i+1}$ are in bijection with the strong orbits of $M$ on $\Omega$.  If we view $\Omega_i/\Omega_{i+1}$ as an $M$-set with zero, then it is a $0$-transitive $M$-set corresponding to the natural action of $M$ on the associated strong orbit by partial maps.  Of course, $\Omega_k$ will be a minimal strong orbit and hence a minimal cyclic sub-$M$-set.

For example, if $N$ is a submonoid of $M$, there are two natural congruences on the regular $M$-set associated to $N$:  namely, the partition of $M$ into weak orbits of the left action of $N$ and the partition of $M$ into the strong orbits of the left action of $N$.  To the best of the author's knowledge, only the latter has every been used in the literature and most often when $M=N$.

More generally, if $\Omega$ is an $M$-set, a relation $\rho$ on $\Omega$ is said to be \emph{stable} if $\alpha\mathrel{\rho} \beta$ implies $\alpha m\mathrel{\rho}\beta m$ for all $m\in M$.

If $\Upsilon$ is any set, then we can make it into an $M$-set via the trivial action $\alpha m=\alpha$ for all $\alpha\in \Upsilon$ and $m\in M$; such $M$-sets are called \emph{trivial}.  This gives rise to a functor $\Delta\colon \pv{Set}\to \pv{Set}^{M^{\mathrm{op}}}$.  The functor $\pi_0\colon \pv{Set}^{M^{\mathrm {op}}}\to \pv{Set}$ provides the left adjoint.  More precisely, we have the following important proposition that will be used later when applying module theory.

\begin{Prop}\label{connectedcomp}
Let $\Omega$ be an $M$-set and $\Upsilon$ a trivial $M$-set.  Then a function $f\colon \Omega\to \Upsilon$ belongs to $\hom_M(\Omega,\Upsilon)$ if and only if $f$ is constant on weak orbits.  Hence $\hom_M(\Omega,\Upsilon)\cong \pv{Set}(\pi_0(\Omega),\Upsilon)$.
\end{Prop}
\begin{proof}
As the weak orbits are $M$-invariant, if we view $\pi_0(\Omega)$ as a trivial $M$-set, then the projection map $\Omega\to \pi_0(\Omega)$ is an $M$-set morphism.  Thus any map $f\colon \Omega\to \Upsilon$ that is constant on weak orbits is an $M$-set morphism. Conversely, suppose that $f\in \hom_M(\Omega,\Upsilon)$ and assume $\alpha \leq \beta\in \Omega$.  Then $\alpha =\beta m$ for some $m\in M$ and so $f(\alpha)=f(\beta m)=f(\beta)m=f(\beta)$.  Thus the relation $\leq$ is contained in $\ker f$.  But $\simeq$ is the equivalence relation generated by $\leq$, whence $f$ is constant on weak orbits.  This completes the proof.
\end{proof}

\begin{Rmk}
The right adjoint of the functor $\Delta$ is the so-called ``global sections'' functor $\Gamma\colon \pv{Set}^{M^{\mathrm{op}}}\to \pv {Set}$ taking an $M$-set $\Omega$ to the set of $M$-invariants of $\Omega$, that is, the set of global fixed points of $M$ on $\Omega$.
\end{Rmk}

We shall also need some structure theory about automorphisms of $M$-sets.
\begin{Prop}\label{schurforsets}
Let $\Omega$ be a transitive $M$-set.  Then every endomorphism of $\Omega$ is surjective.  Moreover, the fixed point set of any non-trivial endomorphism of $\Omega$ is empty.  In particular, the automorphism group of $\Omega$ acts freely on $\Omega$.
\end{Prop}
\begin{proof}
If $f\colon \Omega\to \Omega$ is an endomorphism, then $f(\Omega)$ is $M$-invariant and hence coincides with $\Omega$.  Suppose that $f$ has a fixed point.  Then the fixed point set of $f$ is an $M$-invariant subset of $\Omega$ and thus coincides with $\Omega$.  Therefore, $f$ is the identity.
\end{proof}

In particular, the endomorphism monoid of a finite transitive $M$-set is its automorphism group.

\subsection{Green-Morita theory}
An important role in the theory to be developed is the interplay between $M$ and its subsemigroups of the form $eMe$ with $e$ an idempotent of $M$. Notice that $eMe$ is a monoid with identity $e$.  The group of units of $eMe$ is denoted $G_e$ and is called the \emph{maximal subgroup} of $M$ at $e$.  The set of idempotents of $M$ shall be denoted $E(M)$; more generally, if $X\subseteq M$, then $E(X)=E(M)\cap X$.
First we need to define the tensor product in the context of $M$-sets (cf.~\cite{actsbook,Mac-CWM}).

Let $\Omega$ be a right $M$-set and $\Lambda$ a left $M$-set.  A map $f\colon \Omega\times \Lambda\to \Phi$ of sets is \emph{$M$-bilinear} if $f(\omega m,\lambda) = f(\omega, m\lambda)$ for all $\omega\in \Omega$, $\lambda \in \Lambda$ and $m\in M$.  The universal bilinear map is $\Omega\times \Lambda \to \Omega\otimes_M \Lambda$ given by $(\omega,\lambda)\mapsto \omega\otimes \lambda$.  Concretely, $\Omega\otimes_M \Lambda$ is the quotient of $\Omega\times \Lambda$ by the equivalence relation generated by the relation $(\omega m,\lambda)\approx (\omega, m\lambda)$ for $\omega\in \Omega$, $\lambda\in \Lambda$ and $m\in M$.  The class of $(\omega,\lambda)$ is denoted $\omega\otimes \lambda$.  Suppose that $N$ is a monoid and that $\Lambda$ is also right $N$-set.  Moreover, assume that the left action of $M$ commutes with the right action of $N$; in this case we call $\Lambda$ a \emph{bi-$M$-$N$-set}. Then $\Omega\otimes_M \Lambda$ is a right $N$-set via the action $(\omega\otimes \lambda)n = \omega\otimes (\lambda n)$.  That this is well defined follows easily from the fact that the relation $\approx$ is stable for the right $N$-set structure because the actions of $M$ and $N$ commute.

For example, if $N$ is a submonoid of $M$ and $\{\ast\}$ is the trivial $N$-set, then $\{\ast\}\otimes_N M$ is easily verified to be isomorphic as an $M$-set to the quotient of the regular $M$-set by the weak orbits of the left action of $N$ on $M$.

If $\Upsilon$ is a right $N$-set and $\Lambda$ a bi-$M$-$N$ set, then $\hom_N(\Lambda,\Upsilon)$ is a right $M$-set via the action $(fm)(\lambda) = f(m\lambda)$.
The usual adjunction between tensor product and hom holds in this setting.  We just sketch the proof idea.

\begin{Prop}\label{adjunction}
Let $\Omega$ be a right $M$-set, $\Lambda$ a bi-$M$-$N$-set and $\Upsilon$ a right $N$-set.  Then there is a natural bijection \[\hom_N(\Omega\otimes_M \Lambda,\Upsilon)\cong \hom_M(\Omega,\hom_N(\Lambda,\Upsilon))\] of sets.
\end{Prop}
\begin{proof}
Both sides are in bijection with $M$-bilinear maps $f\colon \Omega\times \Lambda\to \Upsilon$ satisfying $f(\omega,\lambda n) = f(\omega,\lambda)n$ for $\omega\in \Omega$, $\lambda\in \Lambda$ and $n\in N$.
\end{proof}

Something we shall need later is the description of $\Omega\otimes_M \Lambda$ when $\Lambda$ is a free left $M$-set.

\begin{Prop}\label{freesetbasisfortensor}
Let $\Omega$ be a right $M$-set and let $\Lambda$ be a free left $M$-set with basis $B$.  Then $\Omega\otimes_M\Lambda$ is in bijection with $\Omega\times B$.  More precisely, if $\lambda\in \Lambda$, then one can uniquely write $\lambda = m_{\lambda}b_{\lambda}$ with $m_{\lambda}\in M$ and $b_{\lambda}\in B$.  The isomorphism takes $\omega\otimes \lambda$ to $(\omega m_{\lambda},b_{\lambda})$.
\end{Prop}
\begin{proof}
It suffices to show that the map $f\colon \Omega\times \Lambda\to \Omega\times B$ given by $(\omega,\lambda)\mapsto (\omega m_{\lambda},b_{\lambda})$ is the universal $M$-bilinear map.  It is bilinear because freeness implies that if $n\in M$, then since $n\lambda= nm_{\lambda}b_{\lambda}$, one has $m_{n\lambda} = nm_{\lambda}$ and $b_{n\lambda}=b_{\lambda}$.  Thus \[f(\omega,n\lambda)=(\omega nm_{\lambda},b_{\lambda}) = f(\omega n,\lambda)\] and so $f$ is $M$-bilinear.

Suppose now that $g\colon \Omega\times \Lambda\to \Upsilon$ is $M$-bilinear.  Then define $h\colon \Omega\times B\to \Upsilon$ by $h(\omega,b) = g(\omega,b)$.  Then \[h(f(\omega,\lambda))=h(\omega m_{\lambda},b_{\lambda}) = g(\omega m_{\lambda},b_{\lambda}) = g(\omega,\lambda)\] where the last equality uses $M$-bilinearity of $g$ and that $m_{\lambda}b_{\lambda}=\lambda$.  This completes the proof.
\end{proof}

We are now in a position to present the analogue of the Morita-Green theory~\cite[Chapter 6]{Greenpoly} in the context of $M$-sets.  This will be crucial for analyzing transformation monoids, in particular, primitive ones.
The following result is proved in an identical manner to its ring theoretic counterpart.

\begin{Prop}\label{restrictionfunctor}
Let $e\in E(M)$ and let $\Omega$ be an $M$-set.  Then there is a natural isomorphism $\hom_M(eM,\Omega)\cong \Omega e$.
\end{Prop}
\begin{proof}
Define $\p\colon \hom_M(eM,\Omega)\to \Omega e$ by $\p(f)=f(e)$.  This is well defined because $f(e)=f(ee)=f(e)e\in \Omega e$.  Conversely, if $\alpha\in \Omega e$, then one can define a morphism $F_{\alpha}\colon eM\to \Omega$ by $F_{\alpha}(m) = \alpha m$. Observe that $F_\alpha(e)=\alpha e=\alpha$ and so $\p(F_\alpha)=\alpha$.  Thus to prove these constructions are inverses it suffices to observe that if $f\in \hom_M(eM,\Omega)$ and $m\in eM$, then  $f(m)=f(em)=f(e)m=F_{\p(f)}(m)$ for all $m\in eM$.
\end{proof}

We shall need a stronger form of this proposition for the case of principal right ideals generated by idempotents.  Associate to $M$ the category $M_E$ (known as the \emph{idempotent splitting} of $M$) whose object set is $E(M)$ and whose hom sets are given by $M_E(e,f) = fMe$.  Composition \[M_E(f,g)\times M_E(e,f)\to M_E(e,g),\] for $e,f,g\in E(M)$, is given by $(m,n)\mapsto mn$.  This is well defined since $gMf\cdot fMe\subseteq gMe$.  One easily verifies that $e\in M_E(e,e)$ is the identity at $e$.  The endomorphism monoid $M_E(e,e)$ of $e$ is $eMe$.  The idempotent splitting plays a crucial role in semigroup theory~\cite{Tilson,qtheor}.  The following result is well known to category theorists.

\begin{Prop}\label{categoryofprojectiveMsets}
The full subcategory $\pv C$ of $\pv{Set}^{M^{\mathrm {op}}}$ with objects the right $M$-sets $eM$ with $e\in E(M)$ is equivalent to the idempotent splitting $M_E$.  Consequently, the endomorphism monoid of the $M$-set $eM$ is $eMe$ (with its natural left action on $eM$).
\end{Prop}
\begin{proof}
Define $\psi\colon M_E\to \pv C$ on objects by $\psi(e)=eM$; this map is evidentally surjective.
We already know (by Proposition~\ref{restrictionfunctor}) that, for each pair of idempotents $e,f$ of $M$, there is a bijection $\psi_{e,f}\colon fMe\to \hom_M(eM,fM)$ given by $\psi_{e,f}(n) = F_n$ where $F_n(m)=nm$. So to verify that the family $\{\psi_{e,f}\}$, together with the object map $\psi$, provides an equivalence of categories, we just need to verify functoriality, that is, if $n_1\in fMe$ and $n_2\in gMf$, then $F_{n_2}\circ F_{n_1}=F_{n_2n_1}$ and $F_e=1_{eM}$.  For the latter, clearly $F_e(m)=em=m$ for any $m\in eM$.  As to the former, $F_{n_2}(F_{n_1}(m)) = F_{n_2}(n_1m) = n_2(n_1m)=F_{n_2n_1}(m)$.

For the final statement, because $M_E(e,e)=eMe$ it suffices just to check that the actions coincide.  But if $m\in eM$ and $n\in eMe$, then the corresponding endomorphism $F_n\colon eM\to eM$ takes $m$ to $nm$.
\end{proof}

As a consequence, we see that if $e,f\in E(M)$, then  $eM\cong fM$ if and only if there exists $m\in eMf$ and $m'\in fMe$ such that $mm'=e$ and $m'm=f$.  In semigroup theoretic lingo, this is the same thing as saying that $e$ and $f$ are $\mathscr D$-equivalent~\cite{CP,qtheor,Higginsbook,Green}.    If $e,f\in E(M)$ are $\mathscr D$-equivalent, then because $eMe$ is the endomorphism monoid of $eM$ and $fMf$ is the endomorphism monoid of $fM$, it follows that $eMe\cong fMf$ (and hence $G_e\cong G_f$) as $eM\cong fM$.  The reader familiar with Green's relations~\cite{Green,CP} should verify that the elements of $fMe$ representing isomorphisms $eM\to fM$ are exactly those $m\in M$ with $f\R m\eL e$.

It is a special case of more general results from category theory that if $M$ and $N$ are monoids, then $\pv{Set}^{M^{\mathrm {op}}}$ is equivalent to $\pv{Set}^{N^{\mathrm{op}}}$ if and only if $M_E$ is equivalent to $N_E$, if and only if there exists $f\in E(N)$ such that $N= NfN$ and $M\cong fNf$;  see also~\cite{Talwar3}.  In particular, for finite monoids $M$ and $N$ it follows that $\pv{Set}^{M^{\mathrm {op}}}$ and $\pv{Set}^{N^{\mathrm{op}}}$ are equivalent if and only if $M\cong N$ since the ideal generated by a non-identity idempotent of a finite monoid is proper.  The proof goes something like this.  The category $M_E$ is equivalent to the full subcategory on the projective indecomposable objects of $\pv{Set}^{M^{\mathrm {op}}}$ and hence is taken to $N_E$ under any equivalence $\pv{Set}^{M^{\mathrm {op}}}\to \pv{Set}^{N^{\mathrm{op}}}$.  If the object $1$ of $M_E$ is sent to $f\in E(N)$, then $M\cong fNf$ and $N=NfN$.  Conversely, if $f\in E(N)$ with $fNf\cong M$ and $NfN=N$, then $fN$ is naturally a bi-$M$-$N$-set using that $M\cong fNf$.  The equivalence $\pv{Set}^{M^{\mathrm {op}}}\to \pv{Set}^{N^{\mathrm{op}}}$ then sends an $M$-set $\Omega$ to $\Omega\otimes _M fN$.

Fix now an idempotent $e\in E(M)$.  Then $eM$ is a left $eMe$-set and so $\hom_M(eM,\Omega)\cong \Omega e$ is a right $eMe$-set.  The action on $\Omega e$ is given simply by restricting the action of $M$ to $eMe$.  Thus there results a restriction functor $\res_e\colon \pv{Set}^{M^{\mathrm {op}}}\to \pv{Set}^{eMe^{\mathrm {op}}}$ given by \[\res_e(\Omega)=\Omega e.\]  It is easy to check that this functor is exact in the sense that it preserves injectivity and surjectivity.  It follows immediately from the isomorphism $\res_e(-)\cong \hom_M(eM,(-))$ that $\res_e$ has a left adjoint, called \emph{induction}, $\ind_e\colon \pv{Set}^{eMe^{\mathrm {op}}}\to \pv{Set}^{M^{\mathrm {op}}}$ given by \[\ind_e(\Omega) = \Omega\otimes_{eMe} eM.\]  Observe that $\Omega\cong \ind_e(\Omega)e$ as $eMe$-sets via the map $\alpha\mapsto \alpha\otimes e$ (which is the unit of the adjunction).  As this map is natural, the functor $\res_e\ind_e$ is naturally isomorphic to the identity functor on $\pv{Set}^{eMe^{\mathrm {op}}}$.

Let us note that if $\Omega$ is a right $M$-set, then each element of $\Omega\otimes_M Me$ can be uniquely written in the form $\alpha\otimes e$ with $\alpha\in \Omega$.  Thus the natural map $\Omega\otimes_M Me\to \Omega e$ sending $\alpha\otimes e$ to $\alpha e$ is an isomorphism.  Hence Proposition~\ref{restrictionfunctor} shows that $\res_e$ also has a right adjoint $\coind_e\colon \pv{Set}^{eMe^{\mathrm {op}}}\to \pv{Set}^{M^{\mathrm {op}}}$, termed \emph{coinduction},  defined by putting \[\coind_e(\Omega) = \hom_{eMe}(Me,\Omega).\]  Note that $\coind_e(\Omega) e\cong \Omega$ as $eMe$-sets via the map sending $f$ to $f(e)$ (which is the counit of the adjunction) and so $\res_e\coind_e$ is also naturally isomorphic to the identity functor on $\pv{Set}^{eMe^{\mathrm {op}}}$.

The module theoretic analogues of these constructions are essential to much of representation theory, especially monoid representation theory~\cite{Greenpoly,myirreps,rrbg}.

\begin{Prop}\label{inductionprop}
Let $\Omega$ be an $eMe$-set.  Then $\ind_e(\Omega)eM=\ind_e(\Omega)$.
\end{Prop}
\begin{proof}
Indeed, $\alpha\otimes m = (\alpha\otimes e)m\in \ind_e(\Omega)eM$ for $m\in eM$.
\end{proof}

Let us now investigate these constructions in more detail.  First we consider how the strong and weak orbits of $M$ and $Me$ interact.

\begin{Prop}\label{relatedorbits}
Let $\alpha,\beta \in \Omega e$.  Then $\alpha\leq_\Omega \beta$ if and only if $\alpha\leq_{\Omega e}\beta$.  In other words, there is an order embedding $f\colon \Pos(\Omega e)\to \Pos (\Omega)$ taking $\alpha eMe$ to $\alpha M$.
\end{Prop}
\begin{proof}
Trivially, $\alpha\in \beta eMe$ implies $\alpha M\subseteq \beta M$.  Conversely, suppose that $\alpha M\subseteq \beta M$.  Then $\alpha eMe = \alpha Me\subseteq \beta Me=\beta eMe$.
\end{proof}

As an immediate consequence, we have:

\begin{Cor}\label{restrictorbit}
The strong orbits of $\Omega e$ are the sets of the form $\OO_s(\alpha)\cap \Omega e$ with $\alpha\in \Omega e$.  Consequently, if $\Omega$ is a transitive $M$-set, then $\Omega e$ is a transitive $eMe$-set.
\end{Cor}

The relationship between weak orbits of $\Omega$ and $\Omega e$ is a bit more tenuous.

\begin{Prop}\label{weakorbitrestriction}
There is a surjective map $\p\colon \pi_0(\Omega e)\to \pi_0(\Omega)$.  Hence if $\Omega e$ is weakly transitive, then $\Omega$ is weakly transitive.
\end{Prop}
\begin{proof}
The order embedding $\Pos(\Omega e)\to \Pos(\Omega)$ from Proposition~\ref{relatedorbits} induces a map $\p\colon \pi_0(\Omega e)\to \pi_0(\Omega)$ that sends the weak orbit of $\alpha\in \Omega e$ under $eMe$ to its weak orbit $\OO_w(\alpha)$ under $M$.  This map is onto, because $\OO_w(\omega)=\OO_w(\omega e)$ for any $\omega\in \Omega$.
\end{proof}

In general, the map $\p$ in Proposition~\ref{weakorbitrestriction}
is not injective.  For example, let $\Omega = \{1,2,3\}$ and let $M$ consist of the identity map on $\Omega$ together with the maps \[e=\begin{pmatrix} 1 & 2& 3\\ 2 &2 &3\end{pmatrix}, \quad f=\begin{pmatrix} 1 & 2& 3\\ 3 & 2& 3\end{pmatrix}.\]  Then $M$ is weakly transitive on $\Omega$, but $eMe = \{e\}$, $\Omega e=\{2,3\}$ and $eMe$ is not weakly transitive on $\Omega e$.

%

Next we relate the substructures and the quotient structures of $\Omega$ and $\Omega e$ via Galois connections.  The former is the easier one to deal with.  If $\Omega$ is an $M$-set, then $\Sub_M(\Omega)$ will denote the poset of $M$-invariant subsets.

\begin{Prop}\label{substructureGalois}
There is a surjective map of posets \[\psi\colon \Sub_M(\Omega)\to \Sub_{eMe}(\Omega e)\] given by $\Lambda\mapsto \Lambda e$.  Moreover, $\psi$ admits an injective left adjoint given by $\Delta\mapsto \Delta M$.  More concretely, this means that $\Delta M$ is the least $M$-invariant subset $\Lambda$ such that $\Lambda e=\Delta$.
\end{Prop}
\begin{proof}
If $\Lambda$ is $M$-invariant, then $\Lambda eeMe\subseteq \Lambda e$ and hence $\Lambda e\in \Sub_{eMe}(\Omega e)$.
Clearly, $\psi$ is an order preserving map.  If $\Delta\subseteq \Omega e$ is $eMe$-invariant, then $\Delta M$ is $M$-invariant and $\Delta = \Delta e\subseteq \Delta Me=\Delta eMe\subseteq \Delta$.  Thus $\psi$ is surjective.  Moreover, if $\Lambda\in \Sub_M(\Omega)$ satisfies $\Lambda e=\Delta$, then $\Delta M\subseteq \Lambda eM\subseteq \Lambda$.  This completes the proof.
\end{proof}

We now show that induction preserves transitivity.

\begin{Prop}\label{preservetransitive}
Let $\Omega$ be a transitive $eMe$-set.  Then $\ind_e(\Omega)$ is a transitive $M$-set.
\end{Prop}
\begin{proof}
Since $\ind_e(\Omega)e\cong \Omega$ is transitive, if $\Lambda\subseteq \ind_e(\Omega)$ is $M$-invariant, then we have $\Lambda e=\ind_e(\Omega)e$.  Thus Propositions~\ref{inductionprop} and~\ref{substructureGalois} yield $\ind_e(\Omega)=\ind_e(\Omega)eM\subseteq \Lambda$ establishing the desired transitivity.
\end{proof}

It is perhaps more surprising that similar results also hold for the congruence lattice.  If $\Omega$ is an $M$-set, denote by $\Cong_M(\Omega)$ the lattice of congruences on $\Omega$.  If $\equiv$ is a congruence on $\Omega e$, then we define a congruence $\equiv'$ on $\Omega$ by $\alpha\equiv' \beta$ if and only if $\alpha me\equiv \beta me$ for all $m\in M$.

\begin{Prop}\label{enlargecongruence}
Let $\equiv$ be a congruence on $\Omega e$.  Then:
\begin{enumerate}
\item $\equiv'$ is a congruence on $\Omega$;
\item $\equiv'$ restricts to $\equiv$ on $\Omega e$;
\item $\equiv'$ is the largest congruence on $\Omega$ satisfying (2).
\end{enumerate}
\end{Prop}
\begin{proof}
Trivially, $\equiv'$ is an equivalence relation.  To see that it is a congruence, suppose $\alpha\equiv' \beta$ and $n\in M$.  Then, for any $m\in M$, we have $\alpha nme\equiv \beta nme$ by definition of $\equiv'$.  Thus $\alpha n\equiv' \beta n$ and so $\equiv'$ is a congruence.

To prove (2), suppose that $\alpha,\beta\in \Omega e$.  If $\alpha\equiv' \beta$, then $\alpha=\alpha e\equiv \beta e=\beta$ by definition of $\equiv'$.  Conversely, if $\alpha\equiv \beta$ and $m\in M$, then $\alpha me=\alpha eme\equiv \beta eme=\beta me$.  Thus $\alpha\equiv' \beta$.

Finally, suppose that $\approx$ is a congruence on $\Omega$ that restricts to $\equiv$ on $\Omega e$ and assume $\alpha\approx \beta$.  Then for any $m\in M$, we have $\alpha me,\beta me\in \Omega e$ and $\alpha me\approx \beta me$.  Thus $\alpha me\equiv \beta me$ by hypothesis and so $\alpha\equiv' \beta$.  This completes the proof.
\end{proof}

Let us reformulate this result from a categorical viewpoint.

\begin{Prop}\label{quotientGalois}
The map $\varrho\colon \Cong_M(\Omega)\to \Cong_{eMe}(\Omega e)$ induced by restriction is a surjective morphism of posets.  Moreover, it admits an injective right adjoint given by ${\equiv}\mapsto {\equiv'}$.
\end{Prop}

\section{Transformation monoids}\label{sthree}
A \emph{transformation monoid} is a pair $(\Omega, M)$ where $\Omega$ is a set and $M$ is a submonoid of $T_{\Omega}$.  Notice that if $e\in E(M)$, then $(\Omega e,eMe)$ is also a transformation monoid.  Indeed, if $m,m'\in eMe$ and restrict to the same function on $\Omega e$, then for any $\alpha\in \Omega$, we have $\alpha m=\alpha em=\alpha em'=\alpha m'$ and hence $m=m'$.

A transformation monoid $(\Omega, M)$ is said to be \emph{finite} if $\Omega$ is finite.  Of course, in this case $M$ is finite, too. In this paper, we are primarily interested in the theory of finite transformation monoids.  If $|\Omega|=n$, then we say that $(\Omega,M)$ has \emph{degree} $n$.

\subsection{The minimal ideal}
For the moment assume that $(\Omega, M)$ is a finite transformation monoid.  Following standard semigroup theory notation going back to Sch\"utzenberger, if $m\in M$, then $m^{\omega}$ denotes the unique idempotent that is a positive power of $m$.  Such a power exists because finiteness implies $m^k=m^{k+n}$ for some $k>0$ and $n>k$. Then $m^{a+n}=m^a$ for any $a\geq k$ and so if $r$ is the unique natural number $k\leq r\leq k+n-1$ that is divisible by $n$, then $(m^{r})^2=m^{2r}=m^r$.   Uniqueness follows because $\{m^a\mid a\geq k\}$ is easily verified to be a cyclic group with identity $m^r$.
For the basic structure theory of finite semigroups, the reader is referred to~\cite{Arbib} or~\cite[Appendix A]{qtheor}.

If $M$ is a monoid, then a \emph{right ideal} $R$ of $M$ is a non-empty subset $R$ so that $RM\subseteq R$; in other words, right ideals are $M$-invariant subsets of the (right) regular $M$-set. Left ideals are defined dually.   The strong orbits of the regular $M$-set are called \emph{$\R$-classes} in the semigroup theory literature.  An \emph{ideal} is a subset of $M$ that is both a left and right ideal.  If $M$ is a monoid, then $M^{\mathrm {op}}$ denotes the monoid obtained by reversing the multiplication.  Notice that $M^{\mathrm {op}}\times M$ acts on $M$ by putting $x(m,m') = mxm'$.  The ideals are then the $M^{\mathrm {op}}\times M$-invariant subsets; note that this action is weakly transitive.  The strong orbits of this action are called \emph{$\J$-classes} in the semigroup literature.

If $\Lambda$ is an $M$-set and $R$ is a right ideal of $M$, then observe that $\Lambda R$ is an $M$-invariant subset of $\Lambda$.

A key property of finite monoids that we shall use repeatedly is stability.  A monoid $M$ is \emph{stable} if, for any $m,n\in M$, one has that:
\begin{align*}
MmnM=MmM &\iff mnM=mM;\\  MnmM=MmM&\iff Mnm=Mm.
\end{align*}
A proof can be found, for instance, in~\cite[Appendix A]{qtheor}.  We offer a different (and easier) proof here for completeness.

\begin{Prop}
Finite monoids are stable.
\end{Prop}
\begin{proof}
We handle only the first of the two conditions. Trivially, $mnM=mM$ implies $MmnM=MmM$.  For the converse, assume $MmnM=MmM$.  Clearly, $mnM\subseteq mM$.  Suppose that $u,v\in M$ with $umnv=m$.  Then $mM\subseteq umnM$ and hence $|mM|\leq |umnM|\leq |mnM|\leq |mM|$.  It follows that $mM=mnM$.
\end{proof}

An important consequence is the following.  Let $G$ be the group of units of a finite monoid $M$.  By stability, it follows that every right/left unit of $M$ is a unit and consequently $M\setminus G$ is an ideal. Indeed, suppose $m$ has a right inverse $n$, i.e., $mn=1$.  Then $MmM=M=M1M$ and so by stability $Mm=M$.  Thus $m$ has a left inverse and hence an inverse.  The following result is usually proved via stability, but we use instead the techniques of this paper.

\begin{Prop}\label{D=J}
Let $M$ be a finite monoid and suppose that $e,f\in E(M)$.  Then $eM\cong fM$ if and only if $MeM=MfM$.  Consequently, if $e,f\in E(M)$ with $MeM=MfM$, then $eMe\cong fMf$ and hence $G_e\cong G_f$.
\end{Prop}
\begin{proof}
If $eM\cong fM$, then by Proposition~\ref{categoryofprojectiveMsets} that there exist $m\in fMe$ and $m'\in eMf$ with $m'm=e$ and $mm'=f$.  Thus $MeM=MfM$.

Conversely, if $MeM=MfM$, choose $u,v\in M$ with $uev=f$ and put $m=fue$, $m'=evf$.  Then $m\in fMe$, $m'\in eMf$ and $mm'=fueevf=f$.  Thus the morphism $F_m\colon eM\to fM$ corresponding to $m$ (as per Proposition~\ref{categoryofprojectiveMsets}) is surjective and in particular $|fM|\leq |eM|$.  By symmetry, $|eM|\leq |fM|$ and so $F_m$ is an isomorphism by finiteness.

The last statement follows since $eM\cong fM$ implies that $eMe\cong fMf$ by Proposition~\ref{categoryofprojectiveMsets} and hence $G_e\cong G_f$.
\end{proof}

A finite monoid $M$ has a unique minimal ideal $I(M)$.  Indeed, if $I_1,I_2$ are ideals, then $I_1I_2\subseteq I_1\cap I_2$ and hence the set of ideals of $M$ is downward directed and so has a unique minimum by finiteness.  Trivially, $I(M)= MmM=I(M)mI(M)$ for any $m\in I(M)$ and hence $I(M)$ is a simple semigroup (meaning it has no proper ideals).  Such semigroups are determined up to isomorphism by Rees's theorem~\cite{CP,qtheor,Rees} as Rees matrix semigroups over groups.  However, we shall not need the details of this construction in this paper.

If $m\in I(M)$, then $m^{\omega}\in I(M)$ and so $I(M)$ contains idempotents.  Let $e\in E(I(M))$.  The following proposition is a straightforward consequence of the structure theory of theory of finite semigroups.  We include a somewhat non-standard proof using transformation monoids.

\begin{Prop}\label{schutzrep}
Let $M$ be a finite monoid and $e\in E(I(M))$. Then
\begin{enumerate}
\item $eM$ is a transitive $M$-set;
\item $eMe=G_e$;
\item $G_e$ is the automorphism group of $eM$.  In particular, $eM$ is a free left $G_e$-set;
\item If $f\in E(I(M))$, then $fM\cong eM$ and hence $G_e\cong G_f$.
\end{enumerate}
\end{Prop}
\begin{proof}
If $m\in eM$, then $m=em$ and hence, as $MemM=I(M)=MeM$, stability yields $eM=emM=mM$.  Thus $eM$ is a transitive $M$-set. Since $eM$ is finite, Proposition~\ref{schurforsets} shows that the endomorphism monoid of $eM$ coincides with its automorphism group, which moreover acts freely on $eM$.  But the endomorphism monoid is $eMe$ by Proposition~\ref{categoryofprojectiveMsets}.  Thus $eMe=G_e$ and $eM$ is a free left $G_e$-set.  For the final statement, observe that $MeM=I(M)=MfM$ and apply Proposition~\ref{D=J}.
\end{proof}

It is useful to know the following classical characterization of the orbits of $G_e$ on $eM$.
\begin{Prop}\label{Lclasses}
Let $e\in E(I(M))$ and $m,m'\in eM$. Then $G_em=G_em'$ if and only if $Mm=Mm'$.
\end{Prop}
\begin{proof}
This is immediate from the dual of Proposition~\ref{relatedorbits} and the fact that $eMe=G_e$.
\end{proof}

An element $s$ of a semigroup $S$ is called (von Neumann) \emph{regular} if $s=sts$ for some $t\in S$.  For example, every element of $T_{\Omega}$ is regular~\cite{CP}.  It is well known that, for a finite monoid $M$, every element of $I(M)$ is regular in the semigroup $I(M)$.   In fact, we have the following classical result.

\begin{Prop}\label{unionofgroups}
Let $M$ be a finite monoid. Then the disjoint union \[I(M)=\biguplus_{e\in E(I(M))} G_e\] is valid.  Consequently, each element of $I(M)$ is regular in $I(M)$.
\end{Prop}
\begin{proof}
Clearly maximal subgroups are disjoint.  Suppose $m\in I(M)$ and choose $k>0$ so that $e=m^k$ is idempotent.  Then because \[MeM=Mmm^{k-1}M=I(M)=MmM,\] we have by stability that $eM=mM$.  Thus $em=m$ and similarly $me=m$.  Hence $m\in eMe=G_e$.  This establishes the disjoint union.  Clearly, if $g$ is in the group $G_e$, then $gg\inv g=g$ and so $g$ is regular.
\end{proof}

The next result is standard.  Again we include a proof for completeness.

\begin{Prop}\label{LRdontchange}
Let $N$ be a submonoid of $M$ and suppose that $n,n'\in N$ are regular in $N$.  Then $nN=n'N$ if and only if $nM=n'M$ and dually $Nn=Nn'$ if and only if $Mn=Mn'$.
\end{Prop}
\begin{proof}
We handle only the case of right ideals. Trivially, $nN=n'N$ implies $nM=n'M$.  For the converse, suppose $nM=n'M$. Write $n'=n'bn'$ with $b\in N$. Assume that $n=n'm$ with $m\in M$.  Then $n'bn=n'bn'm=n'm=n$ and so $nN\subseteq n'N$.  A symmetric argument establishes $n'N\subseteq nN$.
\end{proof}

In the case $M\leq T_{\Omega}$, the minimal ideal has a (well-known) natural description.
Let $\Omega$ be a finite set and let $f\in T_{\Omega}$.  Define the \emph{rank} of $f$  \[\rk(f)=|f(\Omega)|\] by analogy with linear algebra.  It is well known and easy to prove that $T_{\Omega}fT_{\Omega}=T_{\Omega}gT_{\Omega}$ if and only if $\rk(f)=\rk(g)$~\cite{CP,Higginsbook}.  By stability it follows that $f\in G_{f^{\omega}}$ if and only if $\rk(f)=\rk(f^2)$. The next theorem should be considered folklore.

\begin{Thm}\label{minimalideal}
Let $(\Omega,M)$ be a transformation monoid with $\Omega$ finite.  Let $r$ be the minimum rank of an element of $M$.  Then \[I(M)=\{m\in M\mid \rk(m)=r\}.\]
\end{Thm}
\begin{proof}
Let $J=\{m\in M\mid \rk(m)=r\}$; it is clearly an ideal and so $I(M)\subseteq J$.  Suppose $m\in J$. Then $m^2\in J$ and so $\rk(m^2)=r=\rk (m)$.  Thus $m$ belongs to the maximal subgroup of $T_{\Omega}$ at $m^{\omega}$ and so $m^k=m$ for some $k>0$.  It follows that $m$ is regular in $M$.  Suppose now that $e\in E(I(M))$.  Then we can find $u,v\in M$ with $umv=e$.  Then $eume=e$ and so $eumM=eM$.  Because $\rk(eum)=r=\rk(m)$, it follows that $T_{\Omega}eum=T_{\Omega}m$ by stability.  But $eum$ and $m$ are regular in $M$ (the former by Proposition~\ref{unionofgroups}) and thus $Meum=Mm$ by Proposition~\ref{LRdontchange}.  Thus $m\in I(M)$ completing the proof that $J=I(M)$.
\end{proof}

We call the number $r$ from the theorem the \emph{min-rank} of the transformation monoid $(\Omega,M)$.  Some authors call this the rank of $M$, but this conflicts with the well-established usage of the term ``rank'' in permutation group theory.

In $T_{\Omega}$ one has $fT_{\Omega}=gT_{\Omega}$ if and only if $\ker f=\ker g$ and $T_{\Omega}f=T_{\Omega}g$ if and only if $\Omega f=\Omega g$~\cite{CP,Higginsbook}. Therefore, Proposition~\ref{LRdontchange} immediately yields:

\begin{Prop}
Let $(\Omega,M)$ be a finite transformation monoid and suppose $m,m'\in I(M)$. Then $mM=m'M$ if and only if $\ker m=\ker m'$ and $Mm=Mm'$ if and only if $\Omega m=\Omega m'$.
\end{Prop}

The action of $M$ on $\Omega$ induces an action of $M$ on the power set $P(\Omega)$.   Define \[\min_M\nolimits(\Omega) = \{\Omega m\mid m\in I(M)\}\]
to be the set of images of elements of $M$ of minimal rank.

\begin{Prop}\label{minsetinvariant}
The set $\min_M(\Omega)$ is an $M$-invariant subset of $P(\Omega)$.
\end{Prop}
\begin{proof}
Observe that $\min_M(\Omega) = \{\Omega\}I(M)$ and the latter set is trivially $M$-invariant.
\end{proof}

Let $s\in I(M)$ and suppose that $\ker s = \{P_1,\ldots,P_r\}$.  Then if $X\in \min_M(\Omega)$, the fact that $r=|Xs|=|X|$ implies that $|X\cap P_i|\leq 1$ for $i=1,\ldots, r$.  But since $\ker s$ is a partition into $r=|X|$ blocks, we conclude that $|X\cap P_i|=1$ for all $i=1,\ldots, r$.  We state this as a proposition.

\begin{Prop}\label{kernelpartition}
Let $X\in \min_M(\Omega)$ and $s\in I(M)$.  Suppose that $P$ is a block of $\ker s$.  Then $|X\cap P|=1$.  In particular, right multiplication by $s$ induces a bijection $X\to Xs$.
\end{Prop}

We now restate some of our previous results specialized to the case of minimal idempotents.  See also~\cite{berstelperrinreutenauer}.
\begin{Prop}\label{minimalfacts}
Let $(\Omega,M)$ be a finite transformation monoid and let $e\in E(I(M))$.  Then:
\begin{enumerate}
\item $(\Omega e,G_e)$ is a permutation group of degree the min-rank of $M$;
\item $|\Omega e/G_e|\geq |\pi_0(\Omega)|$;
\item If $M$ is transitive on $\Omega$, then $(\Omega e,G_e)$ is a transitive permutation group.
\end{enumerate}
\end{Prop}

Another useful and well-known fact is that if $(\Omega,M)$ is a finite transitive transformation monoid, then $I(M)$ is transitive on $\Omega$.

\begin{Prop}\label{minidealistrans}
Let $(\Omega,M)$ be a finite transitive transformation monoid.  Then the semigroup $I(M)$ is transitive on $\Omega$ (i.e., there are no proper $I(M)$-invariant subsets).
\end{Prop}
\begin{proof}
If $\alpha\in \Omega$, then $\alpha I(M)$ is $M$-invariant and so $\alpha I(M)=\Omega$.
\end{proof}

In the case that the maximal subgroup $G_e$ of the minimal ideal is trivial and the action of $M$ on $\Omega$ is transitive, one has that each element of $I(M)$ acts as a constant map and $\Omega\cong eM$.  This fact should be considered folklore.

\begin{Prop}\label{constantmapcase}
Let $(\Omega,M)$ be a finite transitive transformation monoid and let $e\in E(I(M))$.  Suppose that $G_e$ is trivial.  Then $I(M)=eM$, $\Omega\cong eM$ and $I(M)$ is the set of constant maps on $\Omega$.
\end{Prop}
\begin{proof}
If $f\in E(I(M))$, then $G_f\cong G_e$ implies $G_f$ is trivial.  Proposition~\ref{unionofgroups} then implies that $I(M)$ consists only of idempotents.  By Proposition~\ref{minimalfacts}, the action of $G_f$ on $\Omega f$ is transitive and hence $|\Omega f|=1$; say $\Omega f=\{\omega_f\}$.   Thus each element of $I(M)$ is a constant map.  In particular, $ef=f$ for all $f\in I(M)$ and hence $eM=I(M)$.  By transitivity of $I(M)$ on $\Omega$ (Proposition~\ref{minidealistrans}), we have that each element of $\Omega$ is the image of a constant map from $I(M)$.  Consequently, we have a bijection $eM\to \Omega$ given by $f\mapsto \omega_f$ (injectivity follows from faithfulness of the action on $\Omega$).  The map is a morphism of $M$-sets because if $m\in M$, then $fm\in I(M)$ and $\Omega fm= \{\omega_fm\}$ and so $\omega_{fm}=\omega_fm$ by definition.  This shows that $\Omega\cong eM$.
\end{proof}

Let us relate $I(M)$ to the socle of $\Omega$.

\begin{Prop}\label{socle}
Let $(\Omega,M)$ be a finite transformation monoid.  Then $\Omega I(M)=\soc \Omega$.  Hence the min-ranks of $\Omega$ and $\soc \Omega$ coincide.
\end{Prop}
\begin{proof}
Let $\alpha\in \soc \Omega$. Then $\alpha M$ is a minimal cyclic sub-$M$-set and hence a transitive $M$-set.  Therefore, $\alpha M=\alpha I(M)$ by transitivity of $M$ on $\alpha M$ and so $\alpha \in \Omega I(M)$.  Conversely, suppose that $\alpha\in \Omega I(M)$, say $\alpha =\omega m$ with $\omega\in \Omega$ and $m\in I(M)$.  Let $\beta\in \alpha M$. We show that $\beta M=\alpha M$, which will establish the minimality of $\alpha M$.  Suppose that $\beta = \alpha n$ with $n\in M$.  Then $\beta = \omega mn$ and $mn\in I(M)$.  Stability now yields $mM=mnM$ and so we can find $n'\in M$ with $mnn'=m$.  Thus $\beta n'=\omega mnn'=\omega m=\alpha$.  It now follows that $\alpha M$ is minimal and hence $\alpha\in \soc \Omega$.
\end{proof}

\subsection{Wreath products}
We shall mostly be interested in transitive (and later $0$-transitive) transformation semigroups.  In this section we relate transitive transformation monoids to induced transformation monoids and give an alternative description of certain tensor products in terms of wreath products.  This latter approach underlies the Sch\"utzenberger representation of a monoid~\cite{Schutzrep,CP,qtheor}.  Throughout this section, $M$ is a finite monoid.

Not all finite monoids have a faithful transitive representation.  A monoid $M$ is called \emph{right mapping} with respect to its minimal ideal if it acts faithfully on the right of $I(M)$~\cite{Arbib,qtheor}.  Regularity implies that if $e_1,\ldots,e_k$ are idempotents forming a transversal to the $\R$-classes of $I(M)$, then $I(M) = \biguplus_{i=1}^m e_kM$. (Indeed, if $mnm=m$, then $mn$ is idempotent and $mM=mnM$.) But all these right $M$-sets are isomorphic (Proposition~\ref{schutzrep}).  Thus $M$ is right mapping with respect to $I(M)$ if and only if $M$ acts faithfully on $eM$ for some (equals any) idempotent of $I(M)$ and so in particular $M$ has a faithful transitive representation.  The converse is true as well.

\begin{Prop}\label{inducedquotient}
Let $(\Omega, M)$ be a transformation monoid and let $e\in E(M)$.  Suppose that $\Omega = \Omega eM$, e.g., if $M$ is transitive.  Then $M$ acts faithfully on $eM$ and there is a surjective morphism $f\colon \ind_e(\Omega e)\to \Omega$ of $M$-sets.
\end{Prop}
\begin{proof}
The counit of the adjunction yields a morphism $f\colon \ind_e(\Omega e)\to \Omega$, which is surjective because \[f(\ind_e(\Omega e)) = f(\ind_e(\Omega e)eM) = \Omega eM=\Omega\] where we have used Proposition~\ref{inductionprop} and that $f$ takes $\ind_e(\Omega e)e$ bijectively to $\Omega e$.  Trivially, if $m,m'\in M$ act the same on $eM$, then they act the same on $\ind_e(\Omega e)=\Omega_e\otimes_{eMe} eM$.  It follows from the surjectivity of $f$ that $m,m'$ also act the same on $\Omega$ and so $m=m'$.
\end{proof}

As a consequence we see that a finite monoid $M$ has a faithful transitive representation if and only if it is right mapping with respect to its minimal ideal.

Suppose that $(\Omega,M)$ and $(\Lambda,N)$ are transformation monoids.  Then $N$ acts on the left of the monoid $M^{\Lambda}$ by endomorphisms by putting $nf(\lambda) = f(\lambda n)$.  The corresponding semidirect product $M^{\Lambda}\rtimes N$ acts faithfully on $\Omega\times \Lambda$ via the action \[(\omega,\lambda)(f,n) = (\omega f(\lambda),\lambda n).\]  The resulting transformation monoid $(\Omega\times \Lambda, M^{\Lambda}\rtimes N)$ is called the \emph{transformation wreath product} and is denoted $(\Omega, M)\wr (\Lambda, N)$.  The semidirect product $M^{\Lambda}\rtimes N$ is denoted $M\wr (\Lambda,N)$.  The wreath product is well known to be associative on the level of transformation monoids~\cite{Eilenberg}.

Suppose now that $M$ is finite and $e\in E(I(M))$.  Notice that since $G_e$ acts on the left of $eM$ by automorphisms, the quotient set $G_e\backslash eM$ has the structure of a right $M$-set given by $G_en\cdot m = G_enm$.  The resulting transformation monoid is denoted $(G_e\backslash eM,\RLM(M))$ in the literature~\cite{qtheor,Arbib}.  The monoid $\RLM(M)$ is called \emph{right letter mapping} of $M$.

Let's consider the following slightly more general situation.  Suppose that $G$ is a group and $M$ is a monoid.  Let $\Lambda$ be a right $M$-set and suppose that $G$ acts freely on the left of $\Lambda$ by automorphisms of the $M$-action. Then $M$ acts naturally on the right of $G\backslash \Lambda$. Let $B$ be a transversal to $G\backslash \Lambda$; then $\Lambda$ is a free $G$-set on $B$.    Suppose that $\Omega$ is a right $G$-set.  Then Proposition~\ref{freesetbasisfortensor} shows that $\Omega\otimes_G \Lambda$ is in bijection with $\Omega\times B$ and hence in bijection with $\Omega\times G\backslash \Lambda$.  If we write $\ov{G\lambda}$ for the representative from $B$ of the orbit $G\lambda$ and define $g_{\lambda}\in G$ by $\lambda =g_{\lambda}\ov{G\lambda}$, then the bijection is $\omega\otimes \lambda\to (\omega g_{\lambda},\ov{G\lambda})\mapsto (\omega g_{\lambda},G\lambda)$.  The action of $M$ is then given by $(\omega,G\lambda)m = (\omega g_{\ov{G\lambda} m},G\lambda m)$.  This can be rephrased in terms of the wreath product, an idea going back to Frobenius for groups and Sch\"utzenberger for monoids~\cite{CP,CP2}; see also~\cite{selfsimilar} for a recent exposition in the group theoretic context.

\begin{Prop}\label{generalizedschutz}
Let $(\Lambda,M)$ be a transformation monoid and suppose that $G$ is a group of automorphisms of the $M$-set $\Lambda$ acting freely on the left.  Let $\Omega$ be a right $G$-set.  Then:
\begin{enumerate}
\item If $\Omega$ is a transitive $G$-set and $\Lambda$ is a transitive $M$-set, then $\Omega\otimes_G \Lambda\cong \Omega\times G\backslash \Lambda$ is a transitive $M$-set.
\item If $\Omega$ is a faithful $G$-set, then the action of $M$ on $\Omega\otimes _G \Lambda\cong \Omega\times G\backslash \Lambda$ is faithful and is contained in the wreath product \[(\Omega,G)\wr (G\backslash \Lambda,\ov M)\] where $\ov M$ is the quotient of $M$ by the kernel of its action on $G\backslash \Lambda$.
\end{enumerate}
\end{Prop}
\begin{proof}
We retain the notation from just before the proof.  We begin with (1). Let $(\alpha_0,G\lambda_0)$ and $(\alpha_1,G\lambda_1)$ be elements of $\Omega\times  G\backslash \Lambda$.  Without loss of generality, we may assume $\lambda_0,\lambda_1\in B$.  By transitivity we can choose $m\in M$ with $\lambda_0m=\lambda_1$.  Then $(\alpha_0,G\lambda_0)m = (\alpha_0,G\lambda_1)$.   Then by transitivity of $G$, we can find $g\in G$ with $\alpha'g=\alpha_1$.  By transitivity of $M$, there exists $m'\in M$ such that $g\lambda_1=\lambda_1m'$.  Then $\ov{G\lambda_1m'}=\lambda_1$ and $g_{\lambda_1 m'}=g$.  Therefore, \[(\alpha_0,G\lambda_1)m' = (\alpha_0 g_{\lambda_1m'},G\lambda_1) = (\alpha_0g,G\lambda_1)=(\alpha_1,G\lambda_1).\]  This establishes the transitivity of $M$ on $\Omega\otimes_G \Lambda$.

To prove (2), first suppose that $m\neq m'$ are elements of $M$.  Then we can find $\lambda\in \Lambda$ such that $\lambda m\neq \lambda m'$.  Then $g\lambda m\neq g\lambda m'$ for all $g\in G$ and so we may assume that $\lambda\in B$.  If $G\lambda m\neq G\lambda m'$, we are done.  Otherwise, $\lambda m = g_{\lambda m}\ov{G\lambda m}$ and $\lambda m' = g_{\lambda m'}\ov{G\lambda m}$ and hence $g_{\lambda m}\neq g_{\lambda m'}$.  Thus by faithfulness of the action of $G$, we have $\alpha\in \Omega$ such that $\alpha g_{\lambda m}\neq \alpha g_{\lambda m'}$.  Therefore, we obtain \[(\alpha,G\lambda)m =(\alpha g_{\lambda m},G\lambda m)\neq (\alpha g_{\lambda m'},G\lambda m) = (\alpha,G\lambda)m'\] establishing the faithfulness of $M$ on $\Omega\otimes_G \Lambda$.

Finally, we turn to the wreath product embedding.  Write $\ov m$ for the class of $m\in M$ in the monoid $\ov M$.  For $m\in M$, we define $f_m\colon G\backslash\Lambda\to \Omega$ by $f_m(G\lambda) = g_{\ov{G\lambda}m}$.  Then $(f_m,\ov m)$ is an element of the semidirect product $G^{G\backslash\Lambda}\rtimes \ov M$ and if $\alpha\in \Omega$ and $\lambda\in \Lambda$, then \[(\alpha,G\lambda)(f_m,\ov m) = (\alpha f_m(G\lambda),G\lambda m)=(\alpha g_{\ov{G\lambda}m},G\lambda m) = (\alpha,G\lambda)m\] as required.  Since the action of $M$ on $\Omega\times G\backslash \Lambda$ is faithful, this embeds $M$ into the wreath product.
\end{proof}

A particularly important case of this result is when $(\Omega,M)$ is a transitive transformation monoid and $G$ is a group of $M$-set automorphisms of $\Omega$; the action of $G$ is free by Proposition~\ref{schurforsets}.  Observing that $\Omega=G\otimes_G \Omega$, we have the following corollary.

\begin{Cor}
Let $(\Omega,M)$ be a transitive transformation monoid and $G$ a group of automorphisms of $(\Omega,M)$.  Then $\Omega$ is in bijection with $G\times G\backslash \Omega$ and the action of $M$ on $\Omega$ is contained in the wreath product $(G,G)\wr (G\backslash \Omega,\ov M)$ where $\ov M$ is the quotient of $M$ by the kernel of its action on $G\backslash \Omega$.
\end{Cor}

Another special case is the following slight generalization of the classical Sch\"utzenberger representation~\cite{CP,Arbib,qtheor}, which pertains to the case $\Omega=G_e$ (as $\ind_e(G_e)\cong eM$); cf.~\cite{CP2}.

\begin{Cor}
Suppose that $M$ is a finite right mapping monoid (with respect to $I(M))$ and let $e\in E(I(M))$. If $\Omega$ is a transitive $G_e$-set, then $\ind_e(\Omega)$ is a transitive $M$-set.  Moreover, if $\Omega$ is faithful, then $\ind_e(\Omega)$ is a faithful $M$-set and $(\ind_e(\Omega),M)$ is contained inside of the wreath product $(\Omega,G_e)\wr (G_e\backslash eM,\RLM(M))$.
\end{Cor}

Thus faithful transitive representations of a right mapping monoid $M$ are, up to division~\cite{Arbib,Eilenberg,qtheor}, the same things as wreath products of the right letter mapping representation with transitive faithful permutation representations of the maximal subgroup of $I(M)$.

\section{Finite $0$-transitive transformation monoids}\label{sfour}
In this section we begin to develop the corresponding theory for finite $0$-transitive transformation monoids.  Much of the theory works as in the transitive case once the correct adjustments are made.  For this reason, we will not tire the reader by repeating analogues of all the previous results in this context.  What we call a $0$-transitive transformation monoid is called by many authors a \emph{transitive partial transformation monoid}.

Assume now that $(\Omega, M)$ is a finite $0$-transitive transformation monoid.  The zero map, which sends all elements of $\Omega$ to $0$, is denoted $0$.

\begin{Prop}\label{haszero}
Let $(\Omega, M)$ be a finite $0$-transitive transformation monoid.  Then the zero map belongs to $M$ and $I(M)=\{0\}$.
\end{Prop}
\begin{proof}
Let $e\in E(I(M))$.  First note that $0\in \Omega e$.  Next observe that if $0\neq \alpha\in \Omega e$, then $\alpha eMe= \alpha Me=\Omega e$ and hence $G_e=eMe$ is transitive on $\Omega e$.  But $0$ is a fixed point of $G_e$ and so we conclude that $\Omega e=\{0\}$ and hence $e=0$.  Then trivially $I(M)=MeM=\{0\}$.
\end{proof}

An ideal $I$ of a monoid $M$ with zero is called \emph{$0$-minimal} if $I\neq 0$ and the only ideal of $M$ properly contained in $I$ is $\{0\}$.  It is easy to see that $I$ is $0$-minimal if and only if $MaM=I$ for all $a\in I\setminus \{0\}$, or equivalently, the action of $M^{\mathrm {op}}\times M$ on $I$ is $0$-transitive.  In a finite monoid $M$ with zero, a $0$-minimal ideal is regular (meaning all its elements are regular in $M$) if and only if $I^2=I$~\cite{CP,qtheor}.  We include a proof for completeness.

\begin{Prop}\label{regularideal}
Suppose that $I$ is a $0$-minimal ideal of a finite monoid $M$.  Then $I$ is regular if and only if $I^2=I$.  Moreover, if $I\neq I^2$, then $I^2=0$.
\end{Prop}
\begin{proof}
If $I$ is regular and $0\neq m\in I$, then we can write $m=mnm$ with $n\in M$ and so $m=m(nm)\in I^2$.  It follows $I^2=I$.  Conversely, if $I^2=I$ and $m\in I\setminus \{0\}$, then we can write $m=ab$ with $a,b\in I\setminus \{0\}$.  Then $MmM=MabM=MaM=MbM$ and so stability yields $mM=aM$ and $Mm=Mb$.  Therefore, we can write $a=mx$ and $b=ym$ and hence $m=mxym$ is regular.

For the final statement, suppose $I\neq I^2$.  Then $I^2$ is an ideal strictly contained in $I$ and so $I^2=0$.
\end{proof}

Of course if $I$ is regular, then it contains non-zero idempotents.  Using this one can easily show~\cite{CP,qtheor} that each element of $I$ is regular in the semigroup $I$. In fact, $I$ is a $0$-simple semigroup and hence its structure is determined up to isomorphism by Rees's theorem~\cite{CP,qtheor,Rees}.

If $\Omega$ is an $M$-set and $\Lambda$ is an $M$-set with $0$, then the map sending each element of $\Omega$ to $0$ is an $M$-set map, which we again call the zero map and denote by $0$.

\begin{Prop}\label{Schuragain}
Let $\Omega$ be an $M$-set and $\Lambda$ a $0$-transitive $M$-set.  Then every non-zero morphism $f\colon \Omega\to \Lambda$ of $M$-sets is surjective.
\end{Prop}
\begin{proof}
If $f\colon \Omega\to \Lambda$ is a non-zero morphism, then $0\neq f(\Omega)$ is $M$-invariant and hence equals $\Lambda$ by $0$-transitivity.
\end{proof}

As a corollary we obtain an analogue of Schur's lemma.

\begin{Cor}\label{freeactionagain}
Let $\Omega$ be a finite $0$-transitive $M$-set.  Then every non-zero endomorphism of $\Omega$ is an automorphism.  Moreover, $\mathrm{Aut}_M(\Omega)$ acts freely on $\Omega\setminus \{0\}$.
\end{Cor}
\begin{proof}
By Proposition~\ref{Schuragain}, any non-zero endomorphism of $\Omega$ is surjective and hence is an automorphism.  Since any automorphism of $\Omega$ fixes $0$ (as it is the unique sink by Proposition~\ref{uniquesink}), it follows that $\Omega\setminus \{0\}$ is invariant under $\mathrm{Aut}_M(\Omega)$.  If $f\in \mathrm{Aut}_M(\Omega)$, then its fixed point set is $M$-invariant and hence is either $0$ or all of $\Omega$.  This shows that the action of $\mathrm{Aut}_M(\Omega)$ on $\Omega \setminus \{0\}$ is free.
\end{proof}

We can now prove an analogue of Proposition~\ref{schutzrep} for $0$-minimal ideals.  Again this proposition is a well-known consequence of the classical theory of finite semigroups.  See~\cite{berstelperrinreutenauer} for the corresponding result in the more general situation of unambiguous representations of monoids.

\begin{Prop}\label{schutzrep2}
Let $M$ be a finite monoid with zero, let $I$ be a regular $0$-minimal ideal and let $e\in E(I)\setminus \{0\}$. Then:
\begin{enumerate}
\item $eM$ is a $0$-transitive $M$-set;
\item $eMe=G_e\cup \{0\}$;
\item $G_e$ is the automorphism group of the $M$-set $eM$ and so in particular, $eM\setminus \{0\}$ is a free left $G_e$-set;
\item If $f\in E(I)\setminus \{0\}$, then $fM\cong eM$ and hence $G_e\cong G_f$; moreover, one has $fMe\setminus \{0\}$ and $eMf\setminus \{0\}$ are in bijection with $G_e$.
\end{enumerate}
\end{Prop}
\begin{proof}
Trivially $0\in eM$.
Suppose that $0\neq m\in eM$.  Then $m=em$ and hence, as $MmM=MemM=MeM$, stability yields $mM=eM$.  Thus $eM$ is a $0$-transitive $M$-set. Since $eM$ is finite, Corollary~\ref{freeactionagain} shows that the endomorphism monoid of $eM$ consists of the zero morphism and its group of units, which acts freely on $eM\setminus \{0\}$.  But the endomorphism monoid is $eMe$ by Proposition~\ref{categoryofprojectiveMsets}.  Thus $eMe=G_e\cup \{0\}$ and $eM\setminus \{0\}$ is a free left $G_e$-set.

Now we turn to the last item.  Since $MeM=I=MfM$, we have that $eM\cong fM$ by Proposition~\ref{D=J}.  Clearly the automorphism group $G_e$ of $eM$ is in bijection with the set of isomorphisms $eM\to fM$; but this latter set is none other than $fMe\setminus \{0\}$.  The argument for $eMf\setminus \{0\}$ is symmetric.
\end{proof}

Of course the reason for developing all this structure is the folklore fact that a finite $0$-transitive transformation monoid has a unique $0$-minimal ideal, which moreover is regular.  Any element of this ideal will have minimal non-zero rank.

\begin{Thm}\label{unique0min}
Let $(\Omega,M)$ be a finite $0$-transitive transformation monoid.  Then $M$ has a unique $0$-minimal ideal $I$; moreover, $I$ is regular and acts $0$-transitively (as a semigroup) on $\Omega$.
\end{Thm}
\begin{proof}
We already know that $0\in M$ by Proposition~\ref{haszero}.  Let $I$ be a $0$-minimal ideal of $M$ (it has one by finiteness).  Then $\Omega I$ is $M$-invariant.  It is also non-zero since $I$ contains a non-zero element of $M$.  Thus $\Omega I=\Omega$.  Therefore, $\Omega I^2=\Omega I=\Omega$ and so $I^2\neq 0$.  We conclude by Proposition~\ref{regularideal} that $I$ is regular.  This also implies the $0$-transitivity of $I$ because if $0\neq \alpha\in \Omega$, then $\alpha I\supseteq \alpha MI=\Omega I=\Omega$.
Finally, suppose that $I'$ is any non-zero ideal of $M$.  Then $\Omega I'\neq 0$ and is $M$-invariant.  Thus $\Omega=\Omega I' = \Omega II'$ and so $0\neq II'\subseteq I\cap I'$.  By $0$-minimality, we conclude $I=I\cap I'\subseteq I'$ and hence $I$ is the unique $0$-minimal ideal of $M$.
\end{proof}

We also have the following analogue of Proposition~\ref{minimalfacts}(3).

\begin{Prop}\label{transitivityofgroup}
Let $(\Omega, M)$ be a finite $0$-transitive transformation monoid with $0$-minimal ideal $I$ and let $0\neq e\in E(I)$.  Then $(\Omega e\setminus \{0\},G_e)$ is a transitive permutation group.
\end{Prop}
\begin{proof}
If $0\neq \alpha\in \Omega e$, then $\alpha eMe=\alpha Me=\Omega e$.  But $eMe=G_e\cup \{0\}$ and hence $\alpha G_e=\Omega e\setminus \{0\}$ (as $0$ is a fixed point for $G_e$).
\end{proof}

Again, in the case that $G_e$ is trivial, one can say more, although not as much as in the transitive case.

\begin{Prop}\label{aperiodicbottom}
Let $(\Omega, M)$ be a finite $0$-transitive transformation monoid with $0$-minimal ideal $I$ and let $0\neq e\in E(I)$.  Suppose that $G_e$ is trivial.  Then each element of $I\setminus\{0\}$ has rank $2$ and $\Omega\cong eM$.
\end{Prop}
\begin{proof}
First observe that since $G_e$ is trivial, Proposition~\ref{transitivityofgroup}  implies that $\Omega e$ contains exactly one non-zero element.  Thus, for each $m\in I\setminus \{0\}$, there is a unique non-zero element $\omega_m\in \Omega$ so that $\Omega m= \{0,\omega_m\}$, as all non-zero elements of $I$ have the same rank and have $0$ in their image.  We claim that $0\mapsto 0$ and $m\mapsto \omega_m$ gives an isomorphism between $eM$ and $\Omega$.  First we verify injectivity.  Since $m\in eM\setminus \{0\}$ implies $eM=mM$, all elements of $eM\setminus \{0\}$ have the same kernel.  This kernel is a partition $\{P_1,P_2\}$ of $\Omega$ with $0\in P_1$.  Then all elements of $eM$ send $P_1$ to $0$ and hence each element of $eM$ is determined by where it sends $P_2$.  Thus $m\mapsto \omega_m$ is injective on $eM$.  Clearly it is a morphism of $M$-sets because if $m\in eM\setminus \{0\}$ and $n\in M$, then either $mn=0$ and hence $\omega_mn\in \Omega mn=\{0\}$ or $\{0,\omega_{mn}\}=\Omega mn=\{0,\omega_mn\}$.  Finally, to see that the map is surjective observe that $\omega_ee=\omega_e$ and so $\{0\}\neq \omega_eeM$.  The $0$-transitivity of $M$ then yields $\omega_eeM=\Omega$.  But then if $0\neq \alpha\in \Omega$, we can find $m\in eM\setminus \{0\}$ so that $\alpha=\omega_em=\omega_{em}=\omega_m$.  This completes the proof.
\end{proof}

One can develop a theory of induced and coinduced $M$-sets with zero and wreath products in this context and prove analogous results, but we avoid doing so for the sake of brevity.  We do need one result on congruences.

\begin{Prop}\label{enlargecongruence2}
Let $(\Omega,M)$ be a finite $0$-transitive transformation monoid with $0$-minimal ideal $I$ and let $0\neq e\in E(I)$.  Suppose that $\equiv$ is a congruence on $(\Omega e\setminus \{0\},G_e)$.  Then there is a unique largest congruence $\equiv'$ on $\Omega$ whose restriction to $\Omega e\setminus \{0\}$ is $\equiv$.
\end{Prop}
\begin{proof}
First extend $\equiv$ to $\Omega e$ by setting $0\equiv 0$.  Then $\equiv$ is a congruence for $eMe=G_e\cup \{0\}$ and any congruence $\sim$ whose restriction to $\Omega e\setminus \{0\}$ equals $\equiv$ satisfies $0\sim 0$.  The result now follows from Proposition~\ref{enlargecongruence}.
\end{proof}

A monoid $M$ that acts faithfully on the right of a $0$-minimal ideal $I$ is said to be \emph{right mapping} with respect to $I$~\cite{Arbib,qtheor}.  In this case $I$ is the unique $0$-minimal ideal of $M$, it is regular and $M$ acts faithfully and $0$-transitively on $eM$ for any non-zero idempotent $e\in E(I)$.  Conversely, if $(\Omega, M)$ is finite $0$-transitive, then one can verify (similarly to the transitive case) that if $0\neq e\in E(I)$, where $I$ is the unique $0$-minimal ideal of $M$, then $M$ acts faithfully and $0$-transitively on $eM$ and hence is right mapping with respect to $I$.  Indeed, if $0\neq \omega\in \Omega e$, then $\omega eM$ is non-zero and $M$-invariant, whence $\Omega =\omega eM$.  Thus if $m,m'\in M$ act the same on $eM$, then they also act the same on $\Omega$.  Alternatively, one can use induced modules in the category of $M$-sets with zero to prove this.

\section{Primitive transformation monoids}
A transformation monoid $(\Omega,M)$ is \emph{primitive} if it admits no non-trivial proper congruences.  In this section, we assume throughout that $|\Omega|$ is finite.  Trivially, if $|\Omega|\leq 2$ then $(\Omega,M)$ is primitive, so we shall also tacitly assume that $|\Omega|\geq 3$,

\begin{Prop}\label{primitive}
Suppose that $(\Omega, M)$ is a primitive transformation monoid with $2<|\Omega|$.  Then
$M$ is either transitive or $0$-transitive.  In particular, $M$ is weakly transitive.
\end{Prop}
\begin{proof}
If $\Delta$ is an $M$-invariant subset, then consideration of $\Omega/\Delta$ shows that either $\Delta=\Omega$ or $\Delta$ consists of a single point.  Singleton invariant subsets are exactly sinks.  However, if $\alpha,\beta$ are sinks, then $\{\alpha,\beta\}$ is an $M$-invariant subset.  Because $|\Omega|>2$, we conclude that $\Omega$ has at most one sink.

First suppose that $\Omega$ has no sinks.  Then if $\alpha\in \Omega$, one has that $\alpha M\neq \{\alpha\}$ and hence by primitivity $\alpha M=\Omega$.  As $\alpha$ was arbitrary, we conclude that $M$ is transitive.

Next suppose that $\Omega$ has a sink $0$.  We already know it is unique. Hence if $0\neq \alpha\in M$, then $\alpha M\neq \{\alpha\}$ and so $\alpha M=\Omega$.  Thus $M$ is $0$-transitive.

The final statement follows because any transitive or $0$-transitive action is trivially weakly transitive.
\end{proof}

The following results constitute a transformation monoid analogue of Green's results relating simple modules over an algebra $A$ with simple modules over $eAe$ for an idempotent $e$, cf.~\cite[Chapter 6]{Greenpoly}.

\begin{Prop}\label{restricttoidempotents}
Let $(\Omega, M)$ be a primitive transformation monoid and $e\in E(M)$.  Then $(\Omega e, eMe)$ is a primitive transformation monoid.  Moreover, if $|\Omega e|>1$, then $\Omega\cong \ind_e(\Omega e)/{='}$ where $='$ is the congruence on $\ind_e(\Omega e)$ associated to the trivial congruence $=$ on $\ind_e(\Omega e)e\cong \Omega e$ as per Proposition~\ref{enlargecongruence}.
\end{Prop}
\begin{proof}
Suppose first that $(\Omega e,eMe)$ admits a non-trivial proper congruence $\equiv$.  Then Proposition~\ref{enlargecongruence} shows that $\equiv'$ is a non-trivial proper congruence on $\Omega$.  This contradiction shows that $(\Omega e,eMe)$ is primitive.

Next assume $|\Omega e|>1$.  The counit of the adjunction provides a morphism \[f\colon \ind_e(\Omega e)\to \Omega.\]  As the image is $M$-invariant and contains $\Omega e$, which is not a singleton, it follows that $f$ is surjective.  Now $\ker f$ must be a maximal congruence by primitivity of $\Omega$.  However, the restriction of $f$ to $\ind_e(\Omega e)e\cong \Omega e$ is injective.  Proposition~\ref{enlargecongruence} shows that $='$ is the largest such congruence on $\ind_e(\Omega e)$.  Thus $\ker f$ is $='$, as required.
\end{proof}

Of course, the case of interest is when $e$ belongs to the minimal ideal.

\begin{Cor}\label{primitivegroup}
Suppose that $(\Omega,M)$ is a primitive transitive transformation monoid and that $e\in E(I(M))$.  Then $(\Omega e,G_e)$ is a primitive permutation group.  If $G_e$ is non-trivial, then $\Omega = \ind_e(\Omega e)/{='}$.
\end{Cor}

This result is analogous to the construction of the irreducible representations of $M$~\cite{myirreps}.

In the transitive case if $G_e$ is trivial, then we already know that $\Omega\cong eM=\ind_e(\Omega e)$  (since $|\Omega e|=1$) and that $I(M)$ consists of the constant maps on $\Omega$ (Proposition~\ref{constantmapcase}).  In this case, things can be quite difficult to analyze.  For instance, let $(\Omega,G)$ be a permutation group and let $(\Omega,\ov G)$ consist of $G$ along with the constant maps on $\Omega$.  Then it is easy to see that $(\Omega,G)$ is primitive if and only if $(\Omega,\ov G)$ is primitive.  The point here is that any equivalence relation is stable for the ideal of constant maps and so things reduce to $G$.

Sometimes it is more convenient to work with the coinduced action.  The following is dual to Proposition~\ref{restricttoidempotents}.

\begin{Prop}\label{primitivecoinduced}
Let $(\Omega, M)$ be a primitive transformation monoid and let $e\in E(M)$ with $|\Omega e|>1$.  Then there is an embedding $g\colon \Omega\rightarrow \coind_e(\Omega e)$ of $M$-sets.  The image of $g$ is $\coind_e(\Omega e)eM$, which is the least $M$-invariant subset containing $\coind_e(\Omega e)e\cong \Omega e$.
\end{Prop}
\begin{proof}
The unit of the adjunction provides the map $g$ and moreover, $g$ is injective on $\Omega e$.  Because $|\Omega e|>1$, it follows that $g$ is injective by primitivity.  For the last statement, observe that $\Omega eM=\Omega$ by primitivity because $|\Omega e|>1$.  Thus $g(\Omega) = g(\Omega e)eM = \coind_e(\Omega e)eM$.
\end{proof}

We hope that the theory of primitive permutation groups can be used to understand transitive primitive transformation monoids in the case the maximal subgroups of $I(M)$ are non-trivial.

Next we focus on the case of a $0$-transitive transformation monoid.

\begin{Prop}\label{primitivegroup2}
Let $(\Omega, M)$ be a $0$-transitive primitive transformation monoid with $0$-minimal ideal $I$ and suppose $0\neq e\in E(I)$.  Then one has that $(\Omega e\setminus \{0\},G_e)$ is a primitive permutation group.
\end{Prop}
\begin{proof}
If $(\Omega e\setminus \{0\},G_e)$ admits a non-trivial proper congruence, then so does $\Omega$ by Proposition~\ref{enlargecongruence2}.
\end{proof}

Again one can prove that $(\Omega,M)$ is a quotient of an induced $M$-set with zero and embeds in a coinduced $M$-set with zero when $|\Omega e\setminus \{0\}|>1$.  In the case that $G_e$ is trivial, we know from Proposition~\ref{aperiodicbottom} that $\Omega\cong eM$ and each element of the $0$-minimal ideal $I$ acts on $\Omega$ by rank $2$ transformations (or equivalently by rank $1$ partial transformations on $\Omega\setminus \{0\}$).

Recall that a monoid $M$ is an \emph{inverse monoid} if, for each $m\in M$, there exists a unique $m^*\in M$ with $mm^*m=m$ and $m^*mm^*$.  Inverse monoids abstract monoids of partial injective maps, e.g., Lie pseudogroups\cite{Lawson}.  It is a fact that the idempotents of an inverse monoid commute~\cite{Lawson,CP}.  We shall use freely that in an inverse monoid one has $eM=mM$ with $e\in E(M)$ if and only if $mm^*=e$ and dually $Me=Mm$ if and only if $m^*m=e$.  We also use that $(mn)^*=n^*m^*$~\cite{Lawson}.

The next result describes all finite $0$-transitive transformation inverse monoids (transitive inverse monoids are necessarily groups).  This should be considered folklore, although the language of tensor products is new in this context; more usual is the language of wreath products.  The corresponding results for the matrix representation associated to a transformation inverse monoid can be found in~\cite{mobius2}.

\begin{Thm}\label{inversecase}
Let $(\Omega,M)$ be a finite transformation monoid with $M$ an inverse monoid.
\begin{enumerate}
\item If $M$ is transitive on $\Omega$, then $M$ is a group.
\item  If $\Omega$ is a $0$-transitive $M$-set, then $M$ acts on $\Omega\setminus \{0\}$ by partial injective maps and $\Omega\cong (\Omega e\setminus \{0\})\otimes_{G_e}eM$ where $e$ is a non-zero idempotent of the unique $0$-minimal ideal $I$ of $M$.
\end{enumerate}
\end{Thm}
\begin{proof}
Suppose first that $M$ is transitive on $\Omega$.  It is well known that the minimal ideal $I(M)$ of a finite inverse monoid is a group~\cite{CP,Arbib,qtheor}.  Let $e$ be the identity of this group.  Then since $I(M)$ is transitive on $\Omega$, we have $\Omega=\Omega e$.  Thus $e$ is the identity of $M$ and so $M=I(M)$ is a group.

Next suppose that $M$ is $0$-transitive on $\Omega$.  Let $I$ be the $0$-minimal ideal of $M$ and let $e\in E(I)\setminus \{0\}$.  We claim that $\alpha e\neq 0$ implies $\alpha\in \Omega e$.  Indeed, if $\alpha e\neq 0$, then $\alpha eI=\Omega$ and so we can write $\alpha=\alpha em$ with $m\in I$.  Then $\alpha eme=\alpha e\neq 0$.  Thus $eme$ is a non-zero element of $eMe=G_e\cup \{0\}$.  Therefore, $e=(eme)^*eme=em^*eme$ and hence $m^*me=m^*mem^*eme=em^*eme=e$.  But $em^*m=m^*me=e$ and thus $e\in m^*mMm^*m=G_{m^*m}\cup \{0\}$.  We conclude $e=m^*m$.  Thus $\alpha =\alpha em=\alpha emm^*m=\alpha eme$ and so $\alpha \in \Omega e$.  Of course, this is true for any idempotent of $E(I)\setminus \{0\}$, not just for $e$.

Now let $f\in E(M)\setminus \{0\}$ and suppose that $\omega f\neq 0$.  We claim $\omega f=\omega$.  Indeed, choose $\alpha\in \Omega f\setminus \{0\}$.  Then $\alpha I=\Omega$ by $0$-transitivity and so we can write $\omega=\alpha m$ with $m\in I$. Then $\omega f=\alpha mf$.  Because $\alpha= \alpha mf(mf)^*(mf)$ it follows that $\alpha mf(mf)^*\neq 0$.  The previous paragraph applied to $mf(mf)^*\in E(I)\setminus \{0\}$ yields $\alpha=\alpha mf(mf)^*=\alpha mfm^*$.  Therefore, $\omega =\alpha m=\alpha mfm^*m=\alpha mf=\omega f$.

Suppose next that $\omega_1,\omega_2\in \Omega\setminus \{0\}$ and $m\in M$ with $\omega_1m=\omega_2m\neq 0$.  Then $\omega_1mm^*=\omega_2mm^*\neq 0$ and so by the previous paragraph $\omega_1=\omega_1mm^*=\omega_2mm^*=\omega_2$.  We conclude that the action of $M$ on $\Omega\setminus \{0\}$ by partial maps is by partial injective maps.

Let $e\in E(I)\setminus \{0\}$ and put $\Lambda =\Omega e\setminus \{0\}$.  Then $(\Lambda,G_e)$ is a transitive permutation group by Proposition~\ref{transitivityofgroup}.  Consider $\Lambda\otimes_{G_e}eM$.  Observe that if $\alpha,\beta\in \Lambda$ and $\alpha g=\beta$ with $g\in G_e$, then $\beta\otimes 0=\alpha g\otimes 0=\alpha \otimes g0=\alpha\otimes 0$.  Thus $\Lambda\times \{0\}$ forms an equivalence class of $\Lambda\otimes_{G_e}eM$ that we denote by $0$.  It is a sink for the right action of $M$ on $\Lambda\otimes_{G_e}eM$ and hence we can view the latter set as a right $M$-set with zero.

Define $F\colon \Lambda\otimes_{G_e}eM\to \Omega$ by $\alpha\otimes m\mapsto \alpha m$.  This is well defined because the map $\Lambda\times eM\to \Omega$ given by $(\alpha,m)\mapsto \alpha m$ is $G_e$-bilinear.  The map $F$ is a morphism of $M$-sets with zero because $F(\alpha\otimes m)m' = \alpha mm' = F(\alpha\otimes mm')$ and $0$ is sent to $0$.  Observe that $F$ is onto.  Indeed, fix $\alpha\in\Lambda$. Then since $\alpha eM=\alpha M=\Omega$ by $0$-transitivity, given $\omega\in \Omega\setminus \{0\}$, we can find $m\in eM$ with $\omega =\alpha m$.  Thus $\omega=F(\alpha\otimes m)$.  We conclude that $F$ is surjective.

To show injectivity, first observe that if $F(\alpha\otimes m)=0$, then $m=0$.  Indeed, assume $m\neq 0$.  Then $m\in eM\setminus \{0\}$ implies that $eM=mM$ and hence $mm^*=e$.  Thus $0=\alpha mm^*=\alpha e=\alpha$.  This contradiction shows that $m=0$ and hence only $0$ maps to $0$.  Next suppose that $F(\alpha\otimes m)=F(\beta\otimes n)$ with $m,n\in eM\setminus \{0\}$.  Then $\alpha m=\beta n$. From $mm^*=e$, we obtain $0\neq \alpha =\alpha e=\alpha mm^* =\beta nm^*$ and $nm^*\in eMe\setminus \{0\}=G_e$.  Then $e=nm^*mn^*$ and so $nm^*m=nm^*mn^*n=en=n$.  Therefore,
  $\alpha\otimes m=\beta nm^*\otimes m=\beta \otimes nm^*m=\beta\otimes n$ completing the proof that $F$ is injective.
\end{proof}

This theorem shows that the study of ($0$-)transitive representations of finite inverse monoids reduces to the case of groups.  It also reduces the classification of primitive inverse transformation monoids to the case of permutation groups.

\begin{Cor}
Let $(\Omega,M)$ be a primitive finite transformation monoid with $M$ an inverse monoid.  Then either $(\Omega,M)$ is a primitive permutation group, or it is $0$-transitive and $G_e=\{e\}$ for any non-zero idempotent $e$ of the unique $0$-minimal ideal of $M$.  In the latter case, $(\Omega,M)\cong (eM,M)$.
\end{Cor}
\begin{proof}
A primitive transformation monoid is either transitive or $0$-transitive (Proposition~\ref{primitive}).  By Theorem~\ref{inversecase}, if $(\Omega,M)$ is transitive, then it is a primitive permutation group.  Otherwise, the theorem provides an isomorphism $(\Omega,M)\cong (\Omega e\setminus \{0\}\otimes_{G_e}eM,M)$ where $e$ is a non-zero idempotent in the $0$-minimal ideal of $M$.  Suppose that $|G_e|>1$.  Since $\Omega e\setminus \{0\}$ is a faithful $G_e$-set, we conclude $|\Omega e\setminus \{0\}|>1$.  Functoriality of the tensor product yields a non-injective, surjective $M$-set morphism \[(\Omega,M)\to (\{\ast\}\otimes_{G_e}eM,M)\cong (G_e\backslash eM,M).\]  As $0$ and $e$ are in different orbits of $G_e$, this morphism is non-trivial.  This contradiction establishes that $G_e$ is trivial.  We conclude that $(\Omega,M)\cong (eM,M)$ by Proposition~\ref{aperiodicbottom}.
\end{proof}

\begin{Rmk}
A finite primitive transformation monoid $(\Omega,M)$ can only have a non-trivial automorphism group $G$ if $M$ is a group.  Indeed, consideration of $G\backslash \Omega$ shows that either $G$ is trivial or transitive.  But if $G$ is transitive, then $M$ is a monoid of endomorphisms of a finite transitive $G$-set and hence is a permutation group. 
\end{Rmk}

\section{Orbitals}

Let us recall that if $(\Omega,G)$ is a transitive permutation group, then the orbits of $G$ on $\Omega^2=\Omega\times \Omega$ are called \emph{orbitals}.  The diagonal orbital $\Delta$ is called the \emph{trivial orbital}.  The \emph{rank} of $G$ is the number of orbitals.  For instance, $G$ has rank $2$ if and only if $G$ is $2$-transitive.  Associated to each non-trivial orbital $\OO$ is an orbital digraph $\Gamma(\OO)$ with vertex set $\Omega$ and edge set  $\OO$.  Moreover, there is a vertex transitive action of $G$ on $\Gamma(\OO)$.  A classical result of D.~Higman is that the weak and strong components of an orbital digraph coincide and that $G$ is primitive if and only if each orbital digraph is connected~\cite{dixonbook,cameron}.  The goal of this section is to obtain the analogous results for transformation monoids.  The inspiration for how to do this comes out of Trahtman's paper~\cite{Trahtman} on the \v{C}ern\'y conjecture for aperiodic automata.  He considers there certain strong orbits of $M$ on $\Omega^2$ and it turns out that these have the right properties to play the role of orbitals.

After coming up with the definition of orbital presented below, I did an extensive search of the literature with Google and found the paper of Scozzafava~\cite{orbitoids}.  In this paper, if $(\Omega,M)$ is a finite transformation monoid, then a minimal strong orbit is termed an \emph{orbitoid}.  Scozzafava then views the orbitoids of $M$ on $\Omega^2$ as the analogue of orbitals.  He provides two pieces of evidence to indicate that his notion of orbital is ``correct''.  The first is that the number of orbitoids of $M$ on $\Omega^2$ is to equal the number of orbitoids of a point stabilizer on $\Omega$, generalizing the case of permutation groups.  The second is that from an orbitoid of $\Omega^2$, one obtains an action of $M$ on a digraph by graph endomorphisms.  However, this approach does not lead to a generalization of Higman's theorem characterizing primitivity of permutation groups in terms of connectedness of non-trivial orbital digraphs.  Suppose for instance that $G$ is a transitive permutation group on $\Omega$ and $M$ consists of $G$ together with the constant maps on $\Omega$.  Then the unique orbitoid of $M$ on $\Omega^2$ is the diagonal $\Delta$ and so one has no non-trivial orbitals in the sense of~\cite{orbitoids}.  On the other hand, it is easy to see that $M$ is primitive if and only if $G$ is primitive.  In fact, it is clear that if $M$ contains constant maps, then there is no non-trivial digraph on $\Omega$ preserved by $M$ if we use the standard notion of digraph morphism.  Our first step is to define the appropriate category of digraphs in which to work.

\subsection{Digraphs and cellular morphisms}
A (simple) \emph{digraph} $\Gamma$ consists of a set of vertices $V$ and an anti-reflexive relation $E$ on $V\times V$.  If $v,w\in V$, then there is an edge from $v$ to $w$, denoted $(v,w)$, if $(v,w)\in E$.  A \emph{walk} $p$ of length $m$ in a digraph is a sequence of vertices $v_0,v_1,\ldots, v_m$ such that, for each $0\leq i\leq m-1$, one has $(v_i,v_{i+1})$ is an edge, or $v_i=v_{i+1}$.  In particular, for each vertex $v$, there is an \emph{empty walk} of length $0$ consisting of only the vertex $v$.   A walk is called \emph{simple} if it never visits a vertex twice.  The walk $p$ is closed if $v_0=v_m$.  A closed non-empty walk is called a \emph{cycle} if the only repetition occurs at the final vertex.  If $v_0,v_1,\ldots, v_m$ is a walk, then a \emph{deletion} is a removal of a subwalk $v_i,v_{i+1}$ with $v_i=v_{i+1}$.  A walk that admits no deletions is called \emph{non-degenerate}; we consider empty walks as non-degenerate. Deletion is confluent and so from any walk $v_0,\ldots,v_m$, we can obtain a unique \emph{non-degenerate} walk $(v_0,\ldots,v_m)^{\wedge}$ by successive deletions (the resulting path may be empty).

Define a preorder on the vertices of $\Gamma$ by putting $v\leq w$ if there is a walk from $w$ to $v$.  Then the symmetric-transitive closure $\simeq$ of $\leq$ is an equivalence relation on the vertices.  If this relation has a single equivalence class, then the digraph $\Gamma$ is said to be \emph{weakly connected} or just \emph{connected} for short.  In general, the \emph{weak components} of $\Gamma$ are the maximal weakly connected subgraphs of $\Gamma$.  They are disjoint from each other and have vertex sets the $\simeq$-equivalence classes (with the induced edge sets).  The digraph $\Gamma$ is \emph{strongly connected} if $v\leq w$ and $w\leq v$ hold for all vertices $v,w$.  In general, the \emph{strong components} are the maximal strongly connected subgraphs.  A strong component is said to be \emph{trivial} if it contains no edges; otherwise it is \emph{non-trivial}.  A digraph is said to be \emph{acyclic} if all its strong components are trivial.  In this case, the preorder $\leq$ is in fact a partial order on the vertex set.  It is easy to see that if a strong component is non-trivial, then each of its edges belongs to a cycle.  Conversely, a digraph in which each edge belongs to a cycle is strongly connected.  In particular, a digraph is acyclic if and only if it contains no cycles, whence the name.

Usually morphisms of digraphs are required to send edges to edges, but we need to consider here a less stringent notion of morphism.  Namely, we allow maps with degeneracies, i.e., that map edges to vertices.  More precisely, if $\Gamma=(V,E)$ and $\Gamma'=(V',E')$ are digraphs, then a \emph{cellular morphism} is a map $f\colon V\to V'$ such that if $(v,w)\in E$, then either $f(v)=f(w)$ or $(f(v),f(w))\in E'$.  The reason for the term ``cellular'' is that if we view graphs as $1$-dimensional CW-complexes, then it is perfectly legal to map a cell to a lower dimensional cell. If $p=v_0,\ldots, v_m$ is a walk in $\Gamma$, then $f(p)=f(v_0),\ldots, f(v_m)$ is a walk in $\Gamma'$; however, non-degenerate walks can be mapped to degenerate walks.   It is trivial to see that if $f\colon \Gamma\to \Gamma'$ is a morphism, then $f$ takes weak components of $\Gamma$ into weak components of $\Gamma'$ and strong components of $\Gamma$ into strong components of $\Gamma'$.

A cycle $C$ in a digraph $\Gamma$ is \emph{minimal} if it has minimal length amongst all cycles of $\Gamma$.  Minimal cycles exist in non-acyclic digraphs and have length at least $2$ because we do not allow loop edges.

\begin{Prop}\label{minimaltominimal}
Let $f\colon \Gamma\to \Gamma$ be a cellular endomorphism of $\Gamma$ and let $C$ be a minimal cycle of $\Gamma$.  Then either $f(C)$ is a minimal cycle or $f(C)^{\wedge}$ is empty.
\end{Prop}
\begin{proof}
Let $m$ be the length of $C$ and suppose $f(C)^{\wedge}$ is non-empty.  Then $f(C)^{\wedge}$ is a closed path of length at most $m$.  If it is not a cycle, then it contains a proper subwalk that is a cycle of length smaller than the length of $C$, a contradiction.  Thus $f(C)^{\wedge}$ is a cycle.  But then minimality of $C$ implies that $f(C)^{\wedge}$ has length $m$.  Thus $f(C)=f(C)^{\wedge}$ is a minimal cycle.
\end{proof}

By an action of a monoid $M$ on a digraph $\Gamma=(V,E)$, we mean an action by cellular morphisms.  In other words, $M$ acts on $V$ in such a way that the reflexive closure of $E$ is stable for the action of $M$.   We say the action is \emph{vertex transitive} if $M$ is transitive on $V$; we say that it is \emph{edge transitive} if either $M$ acts transitively on $E$ or $M$ acts $0$-transitively on $(E\cup \Delta)/\Delta$ where $\Delta = \{(v,v)\mid v\in V\}$ is the diagonal.  Equivalently, for each pair of edges $e,f\in E$, there is an element $m\in M$ with $em=f$ where in the setting of monoid actions on digraphs, we shall use the notation of right actions.

\begin{Lemma}\label{edgetransitive}
Suppose that $\Gamma$ is a non-acyclic digraph admitting an edge transitive monoid $M$ of cellular endomorphisms.  Then every edge of $\Gamma$ belongs to a minimal cycle.
\end{Lemma}
\begin{proof}
Let $C$ be a minimal cycle of $\Gamma$ and fix an edge $e$ of $C$. Suppose now that $f$ is an arbitrary edge of $\Gamma$.  By edge transitivity, there exists $m\in M$ with $em=f$.  Since $(Cm)^{\wedge}$ is non-empty (it contains the edge $f$), it follows that $Cm$ is minimal by Proposition~\ref{minimaltominimal}.  This completes the proof.
\end{proof}

An immediate corollary of the lemma is the following result.

\begin{Cor}\label{acyclicorstrong}
Suppose that $\Gamma$ is a digraph admitting an edge transitive monoid of cellular endomorphisms.  Then either $\Gamma$ is acyclic or each weak component of $\Gamma$ is strongly connected.
\end{Cor}
\begin{proof}
If $\Gamma$ is not acyclic, then Lemma~\ref{edgetransitive} shows that each edge of $\Gamma$ belongs to a cycle.  It is then immediate that each weak component is strongly connected (since the relation $\leq$ is symmetric in this case).
\end{proof}

\subsection{Orbital digraphs}
Suppose now that $(\Omega,M)$ is a transitive transformation monoid.  Then $M$ acts on $\Omega^2=\Omega\times \Omega$ by $(\alpha,\beta)m = (\alpha m,\beta m)$.  Notice that $\Delta=\{(\alpha,\alpha)\mid \alpha\in \Omega\}$ is a (minimal) strong orbit.  We call $\Delta$ the \emph{trivial orbital} of $M$.  A strong orbit $\OO\neq \Delta$ is an \emph{orbital} if it is minimal in the poset $(\Omega^2/M)\setminus \{\Delta\}$, or equivalently if it is $0$-minimal in $\Omega^2/\Delta$.  Such orbitals are called \emph{non-trivial}.  This coincides with the usual group theoretic notion when $M$ is a group~\cite{dixonbook,cameron}.  Non-trivial orbitals were first studied by Trahtman~\cite{Trahtman} under a different name in the context of the \v{C}ern\'y conjecture.  The number of orbitals of $M$ is called the \emph{rank} of $M$ because this is the well-established terminology in group theory. From now on we assume that $\Omega$ is finite in this section.

For permutation groups, it is well known~\cite{dixonbook,cameron} that the number of orbitals is equal to the number of suborbits (recall that a suborbit is an orbit of the point stabilizer).  This is not the case for transformation monoids.  For example, if $\Omega$ is a finite set of size $n$ and $M$ consists of the identity map and the constant maps, then there are $n^2-n+1$ orbitals, which is larger than the number of points of $\Omega$.

If $\OO$ is a non-trivial orbital, then the corresponding \emph{orbital digraph} $\Gamma(\OO)$ has vertex set $\Omega$ and edge set $\OO$.  Since $\OO$ is a strong orbit, it follows that $M$ acts edge transitively on $\Gamma(\OO)$ by cellular morphisms.  Hence we have the following immediate consequence of Corollary~\ref{acyclicorstrong}.

\begin{Thm}\label{posetorstrongconnected}
Let $(\Omega,M)$ be a transformation monoid and let $\OO$ be a non-trivial orbital.  Then the orbital digaph $\Gamma(\OO)$ is either acyclic or each weak component of $\Gamma(\OO)$ is strongly connected.
\end{Thm}

It was shown by Trahtman~\cite{Trahtman} that if $M$ is aperiodic, then $\Gamma(\OO)$ is always acyclic (using  different terminology: he speaks neither of digraphs nor orbitals).  Here we recall that a finite monoid $M$ is \emph{aperiodic} if each of its maximal subgroups $G_e$ with $e\in E(M)$ is trivial, or equivalently, if $M$ satisfies an identity of the form $x^n=x^{n+1}$.  On the other hand, if $M$ is a non-trivial group, then each weak component of $\Gamma(\OO)$ is strongly connected~\cite{dixonbook,cameron}.

Let $\OO$ be a non-trivial orbital.  Then $M$ either acts transitively on $\OO$ (if it is a minimal strong orbit) or $0$-transitively on the set $\til{\OO}=(\OO\cup \Delta)/\Delta$. In either case, the action need not be faithful.
For example if $I(M)$ consists of the constant maps on $\Omega$, then all of $I(M)$ acts as the zero map on $\til \OO$.  Let $M(\OO)$ be the faithful quotient.    If $M$ is aperiodic, then so is $M(\OO)$.

\begin{Thm}\label{superTraht}
Let $(\Omega,M)$ be a finite transformation monoid and suppose that $\OO$ is a non-trivial orbital.
\begin{enumerate}
 \item If $M$ acts transitively on $\OO$, then $\Gamma(\OO)$ is acyclic if $G_e$ is trivial for $e\in E(I(M(\OO)))$.
 \item If $M$ acts $0$-transitively on $\til \OO$ and $e\in E(I)\setminus \{0\}$ where $I$ is the $0$-minimal ideal of $M(\OO)$, then $\Gamma(\OO)$ is acyclic if $G_e$ is trivial.
\end{enumerate}
\end{Thm}
\begin{proof}
We handle (2) only as (1) is similar, but simpler.  Suppose that $G_e$ is trivial, but that $\Gamma(\OO)$ is not acyclic.  Since $G_e$ is transitive on $\OO e\setminus \Delta$, we have  $|\OO e\setminus \Delta|=1$.  Let $(\alpha,\beta)\in \OO e\setminus \Delta$.  By Lemma~\ref{edgetransitive}, $(\alpha,\beta)$ belongs to some minimal cycle $C$.  Let $m\in M$ with $m$ mapping to $e$ in $M(\OO)$.  Then $(\alpha,\beta)m=(\alpha,\beta)$ and so $Cm$ is a minimal cycle by Proposition~\ref{minimaltominimal}.  If $X$ is the set of edges of $C$, this yields that $Xe$ is a subset of $\OO e\setminus \Delta$ of size greater than $1$.  This contradiction shows that $\Gamma(\OO)$ is acyclic.
\end{proof}

Theorem~\ref{superTraht} admits the following corollary, due to Trahtman with a different formulation.

\begin{Cor}[Trahtman~\cite{Trahtman}]
Let $(\Omega, M)$ be a transitive finite transformation monoid with $M$ aperiodic.  The each non-trivial orbital digraph $\Gamma(\OO)$ is acyclic and hence defines a non-trivial partial order on $\Omega$ that is stable for the action of $M$.
\end{Cor}

If $(\Omega, G)$ is a finite transitive permutation group, then a classical result of D.~Higman says that $G$ is primitive if and only if each non-trivial orbital digraph is strongly connected (equals weakly connected in this context)~\cite{dixonbook,cameron}.  We now prove the transformation monoid analogue.  It is this result that justifies our choice of the notion of an orbital.

\begin{Thm}\label{orbitalthm}
A finite transitive transformation monoid  $(\Omega,M)$ is primitive if and only if each of its non-trivial orbital digraphs is weakly connected.
\end{Thm}
\begin{proof}
Suppose first that $(\Omega,M)$ is primitive and let $\OO$ be a non-trivial orbital.  Then the partition of $\Omega$ into the weak components of $\Gamma(\OO)$ is a non-trivial congruence. Indeed, as $M$ acts by cellular morphisms, it preserves the weak components; moreover, $\Gamma(\OO)$ has at least one edge so not all weak components are trivial.  It follows by primitivity that there is just one weak component, i.e., $\Gamma(\OO)$ is weakly connected.

Conversely, assume that each non-trivial orbital digraph is weakly connected and let $\equiv$ be a non-trivial congruence on $\Omega$.  Then $\equiv$ is an $M$-invariant subset of $\Omega^2$ strictly containing the diagonal $\Delta$.  By finiteness, we conclude that $\equiv$ contains a minimal strong orbit of $\Omega^2\setminus \{\Delta\}$, that is, there is a non-trivial orbital $\OO$ with $\OO\subseteq {\equiv}$.  The weak components of $\Gamma(\OO)$ are the equivalence classes of the equivalence relation generated by $\OO$ and hence each weak component of $\Gamma(\OO)$ is contained in a single $\equiv$-class.  But $\Gamma(\OO)$ is weakly connected, so $\Omega$ is contained in a single $\equiv$-class, that is, $\equiv$ is not a proper congruence.  This completes the proof that $(\Omega,M)$ is primitive.
\end{proof}

As a corollary, we obtain the following.

\begin{Cor}
Let $(\Omega,M)$ be a primitive finite transitive transformation monoid with $M$ aperiodic.  Then $\Omega$ admits a stable connected partial order.
\end{Cor}

Later on, it will be convenient to have a name for the set of weak orbits of $M$ on $\Omega^2$.  We shall call them \emph{weak orbitals}.

\section{Transformation modules}
Our goal now is to study the representations associated to a transformation monoid.  The theory developed here has a different flavor from the group case because there is an interesting duality that arises.

Fix for this section a finite transformation monoid $(\Omega,M)$ and a field $K$ of characterstic $0$.  Let $KM$ be the corresponding monoid algebra.  Associated to the $M$-set $\Omega$ are a right $KM$-module and a left $KM$-module together with a dual pairing.  This pairing has already been implicitly exploited in a number of papers in the \v{C}ern\'y conjecture literature, e.g.,~\cite{dubuc,Kari,mycerny,averaging}.

The \emph{transformation module} associated to $(\Omega,M)$ is the right $KM$-module $K\Omega$.  That is we take a $K$-vector space with basis $\Omega$ and extend the action of $M$ on $\Omega$ linearly: formally, for $m\in M$, define \[\left(\sum_{\omega\in \Omega} c_{\omega}\omega\right) m = \sum_{\omega\in\Omega}c_{\omega}\omega m.\]  The \emph{dual transformation module} is the space $K^{\Omega}$ of $K$-valued functions on $\Omega$ with the left $KM$-module structure given by $mf(\omega) = f(\omega m)$ for $m\in M$ and $f\colon \Omega\to M$.  When $M$ is a group, these two representations are the same under the natural correspondence between left modules and right modules, but for monoids these modules are simply dual to each other.

There is a non-degenerate pairing $\langle\ ,\ \rangle\colon K\Omega\times K^{\Omega}\to K$ given by
\begin{equation}\label{pairing}
\langle \alpha,f\rangle = f(\alpha)
\end{equation}
for $\alpha\in \Omega$.   The pairing on general linear combinations is given by \[\left\langle \sum_{\alpha\in \Omega}c_{\alpha}\alpha,f\right\rangle = \sum_{\alpha\in \Omega}c_{\alpha}f(\alpha).\]

Observe that $K^{\Omega}$ has basis the Dirac functions $\delta_{\omega}$ with $\omega\in \Omega$.  If $m\in M$, then one verifies that \[m\delta_{\omega} = \sum_{\alpha\in \omega m\inv} \delta_{\alpha}\] and more generally if $S\subseteq \Omega$ and $I_S$ denotes the indicator (or characteristic) function of $S$, then \[mI_S=I_{Sm\inv}.\]  Intuitively, the action of $M$ on  $K^{\Omega}$ is by inverse images and this is why $K\Omega$ and $K^\Omega$ contain the same information in the case of groups.

The following adjointness holds.

\begin{Prop}\label{adjoint}
The left and right actions of $m\in M$ on $K\Omega$ and $K^{\Omega}$ are adjoint.  That is, for $v\in K\Omega$, $f\in K^{\Omega}$ and $m\in M$, one has \[\langle vm,f\rangle = \langle v,mf\rangle\]
\end{Prop}
\begin{proof}
It suffices by linearity to handle the case $v=\alpha\in \Omega$.  Then \[\langle \alpha m,f\rangle = f(\alpha m) =mf(\alpha) = \langle \alpha,mf\rangle,\] as required.
\end{proof}

As a consequence, we see that $K^{\Omega}$ is dual to $K\Omega$, that is, \[K^{\Omega}\cong \hom_K(K\Omega,K)\] as left $KM$-modules.  We remark that the bases $\Omega$ and $\{\delta_{\omega}\mid \omega\in \Omega\}$ are dual with respect to the pairing \eqref{pairing}.

If $|\Omega|=n$ and we fix an ordering $\Omega =\{\omega_1,\ldots, \omega_n\}$, then it is convenient to identify elements of $K\Omega$ with row vectors in $K^n$ and elements of $K^{\Omega}$ with column vectors (by associating $f$ with the column vector $(f(\omega_1),\ldots,f(\omega_n))^T$).  The dual pairing then turns into the usual product of a row vector with a column vector.  If $\rho\colon M\to M_n(K)$ is the matrix representation afforded by the right $KM$-module $K\Omega$, then the action on column vectors is the matrix representation afforded by the left $KM$-module $K^{\Omega}$.

We mention the following trivial observation.

\begin{Prop}
Let $(\Omega,M)$ be a finite transformation monoid and suppose that $\OO_1,\ldots,\OO_s$ are the weak orbits of $M$.  Then $K\Omega\cong \bigoplus_{i=1}^sK\OO_i$ and $K^\Omega\cong \bigoplus_{i=1}^sK^{\OO_i}$ where we identify $K^{\OO_i}$ with those functions $\Omega\to K$ supported on $\OO_i$, for $1\leq i\leq s$.
\end{Prop}

Thus for most purposes, it suffices to restrict our attention to the weakly transitive case.

\subsection{The subspace of $M$-invariants}
Let $V$ be a left/right $KM$-module.  Then $V^M$ denotes the subspace of \emph{$M$-invariants}, that is, of all vectors fixed by $M$.  If $K$ is the trivial left/right $KM$-module, then $V^M\cong \hom_{KM}(K,V)$.  Unlike the case of groups, it is not in general true that $\hom_{KM}(K,V)\cong \hom_{KM}(V,K)$.  In fact, we shall see in a moment that in most cases $K\Omega^M=\{0\}$, whereas the $K$-dimension of $\hom_{KM}(K\Omega,K)$ is the number of weak orbits of $M$.  It is also the case that the module $K\Omega$ is almost never semisimple and quite often the multiplicity of the trivial module as a composition factor of $K\Omega$ is strictly greater than the number of weak orbits of $M$.

The following result generalizes a standard result from permutation group theory.

\begin{Prop}\label{fixedset}
Consider $\hom_{KM}(K\Omega,K)$ where $K$ is given the structure of a trivial $KM$-module and $\hom_M(\Omega,K)$ where $K$ is given the structure of a trivial $M$-set.  Then there are $K$-vector space isomorphisms \begin{equation}\label{eq:fixedset}
\hom_{KM}(K\Omega,K)\cong \hom_M(\Omega,K)=(K^{\Omega})^M\cong K^{\pi_0(\Omega)}.
\end{equation}
More precisely, $f\in (K^{\Omega})^M$ if and only if it is constant on weak orbits of $\Omega$.  Consequently, $\dim_K \hom_{KM}(K\Omega,K)=\dim_K  (K^{\Omega})^M$ is the number of weak orbits of $M$ on $\Omega$.
\end{Prop}
\begin{proof}
A $K$-linear map $T\colon K\Omega\to K$ is the same thing as a map $\Omega\to K$ because $\Omega$ is a basis for $K\Omega$.  Clearly, $T$ is a $KM$-module morphism if and only if the associated mapping $\Omega\to K$ is an $M$-set morphism.  This provides the first isomorphism of \eqref{eq:fixedset}.  Proposition~\ref{connectedcomp} shows that $f\colon \Omega\to K$ is an $M$-set morphism if and only if it is constant on weak orbits yielding the isomorphism of the second and fourth terms of \eqref{eq:fixedset}.  Finally, observe that $f\colon \Omega\to K$ is an $M$-set map if and only if $f(\omega m)=f(\omega)$ for all $\omega\in \Omega$, $m\in M$.  But this is equivalent to asking $mf=f$ for all $m\in M$, i.e., $f\in (K^{\Omega})^M$.  This completes the proof.
\end{proof}

The situation for $K\Omega^M$ is quite different. It is well known that a finite monoid $M$ admits a surjective maximal group image homomorphism $\sigma\colon M\to G(M)$ where $G(M)$ is a finite group.  This map is characterized by the universal property that if $\p\colon M\to H$ is a homomorphism from $M$ into a group $H$, then there is a unique homomorphism $\psi\colon G(M)\to H$ so that
\[\xymatrix{M\ar[r]^\sigma\ar[rd]_\p & G(M)\ar@{..>}[d]^\psi \\ & H}\]
commutes.  Using the fact that a finite monoid is a group if and only if it has a unique idempotent, one can describe $G(M)$ as the quotient of $M$ by the least congruence for which all idempotents are equivalent and $\sigma$ as the quotient map.  Alternatively, it is the quotient by the intersection of all congruences on $M$ whose corresponding quotient is a group.

\begin{Prop}\label{invariantsonright0}
If $(\Omega,M)$ is a transformation monoid, then $K\Omega^M\neq 0$ if and only if there is an $M$-invariant subset $\Lambda$ fixed by all idempotents of $M$.
\end{Prop}
\begin{proof}
Suppose first that $\Lambda$ is an $M$-invariant subset fixed by all idempotents of $M$.  Then $\Lambda$ is naturally a $G(M)$-set and $K\Lambda^M=K\Lambda^{G(M)}$.  Group representation theory then yields that $\dim_K K\Lambda^{G(M)}$ is the number of orbits of $G(M)$ on $\Lambda$, and so is non-zero.  Thus $K\Omega^M\neq 0$.

Next suppose that $v\in K\Omega^M$.  Then $ve=v$ for all idempotents $e\in E(M)$, so $v\in \bigcap_{e\in E(M)}K\Omega e$.  Suppose that $v=\sum_{\lambda\in \Lambda}c_{\lambda}\lambda$ with $c_{\lambda}\neq 0$ for all $\lambda\in \Lambda$.  Then $\Lambda\subseteq \Omega e$ for all $e\in E(M)$.  Also, if $m\in M$, then $vm=v$ implies that $\Lambda m=\Lambda$ and so $\Lambda$ is $M$-invariant.
\end{proof}

As corollaries, we obtain the following results.

\begin{Cor}\label{invariantsonright}
Suppose that all elements of $I(M)$ have the same image $\Lambda$.  Then $K\Omega^M\neq 0$.
\end{Cor}
\begin{proof}
Let $e$ be an idempotent of $M$ and choose $m\in I(M)$.  Then $me\in I(M)$ and so $\Lambda = \Omega me\subseteq \Omega e$.  Thus all idempotents of $M$ fix $\Lambda$.  Since $\min_M(\Omega) = \{\Lambda\}$, it follows that $\Lambda$ is $M$-invariant.  The result now follows from Proposition~\ref{invariantsonright0}.
\end{proof}

The next corollary shows that in the transitive setting only groups admit non-trivial invariants.

\begin{Cor}
Suppose that $(\Omega, M)$ is transitive.  Then $K\Omega^M\neq 0$ if and only if $M$ is a group.  In particular, if $(\Omega,M)$ is transitive, then the module $K\Omega$ is semisimple if and only if $M$ is a group.
\end{Cor}
\begin{proof}
If $M$ is a group, then $\dim_K K\Omega^M$ is the number of orbits of $M$ on $G$ and hence is non-zero.  For the converse, suppose $K\Omega^M\neq 0$.  Then Proposition~\ref{invariantsonright0} implies that there is an $M$-invariant subset $\Lambda\subseteq \Omega$ such that every idempotent of $M$ fixes $\Lambda$.  But transitivity implies $\Lambda=\Omega$.  Thus the unique idempotent of $M$ is the identity.  We conclude that $M$ is a group.

The final statement follows because if $K\Omega$ is semisimple, then the fact that $\hom_{KM}(K\Omega,K)\neq 0$ (by Proposition~\ref{fixedset}) implies that the trivial representation is a subrepresentation of $K\Omega$.  But this means that $K\Omega^M\neq 0$ and so $M$ is a group.  Conversely, if $M$ is a group, then $K\Omega$ is semisimple by Maschke's theorem.
\end{proof}

Let us now interpret some of these results for associated actions of $M$.
A $K$-bilinear form $B\colon K\Omega\times K\Omega\to K$ is said to be \emph{$M$-invariant} if \[B(vm,wm)=B(v,w)\] for all $v,w\in K\Omega$ and $m\in M$.   Let $\bil_K(K\Omega)$ be the space of $K$-bilinear forms on $K\Omega$.  There is a natural left $KM$-module structure on $\bil_K(K\Omega)$ given by putting $(mB)(v,w) = B(vm,wm)$.  Then $\bil_K(K\Omega)^M$ is the space of $M$-invariant $K$-bilinear forms.  As $K$-bilinear forms are determined by their values on a basis, it is easy to see that $\bil_K(K\Omega)\cong K^{\Omega\times \Omega}$ as $KM$-modules.  Moreover, $\bil_K(K\Omega)^M\cong (K^{\Omega\times \Omega})^M$.  Thus we have proved the following.

\begin{Prop}
The dimension of the space of $M$-invariant $K$-bilinear forms on $K\Omega$ is the number of weak orbitals of $(\Omega,M)$.
\end{Prop}

Let $\Omega^{\{2\}}$ denote the subset of $P(\Omega)$ consisting of all $1$- and $2$-element subsets.  Then $\Omega^{\{2\}}$ is $M$-invariant and can be identified with the quotient of $\Omega^2$ by the equivalence relation putting $(\alpha,\omega)\equiv (\omega,\alpha)$ for all $\alpha,\omega$.  It is then easy to see that $K^{\Omega^{\{2\}}}$ is isomorphic as a $KM$-module to the space of symmetric $K$-bilinear forms on $K\Omega$ and hence $(K^{\Omega^{\{2\}}})^M$ is the space of $M$-invariant symmetric $K$-bilinear forms on $K\Omega$.  Thus we have proved:

\begin{Prop}
The  dimension of the space of $M$-invariant symmetric $K$-bilinear forms on $K\Omega$ is the number of weak orbits of $M$ on $\Omega^{\{2\}}$.
\end{Prop}

\subsection{The augmentation submodule}
If $(\Omega,M)$ is a finite transformation monoid, one always has the augmentation map $\varepsilon\colon K\Omega\to K$ given by \[\varepsilon(v)=\langle v,I_{\Omega}\rangle.\]  If $v=\sum_{\omega\in\Omega}c_{\omega}\omega$, then $\varepsilon(v)=\sum_{\omega\in\Omega}c_{\omega}$.  Clearly $I_{\Omega}$ is constant on weak orbits and so $\varepsilon$ is a $KM$-module homomorphism (where $K$ is given the trivial $KM$-module structure).  Thus $\ker \varepsilon$ is a $KM$-submodule, called the \emph{augmentation submodule}, and is denoted $\mathrm{Aug}(K\Omega)$.  A key fact, which plays a role in the \v{C}ern\'y conjecture literature, is that $m\in M$ is a constant map if and only if $m$ annihilates the augmentation submodule.  Indeed $\mathrm{Aug}(K\Omega)$ consists of those vectors $v=\sum_{\omega\in\Omega}c_{\omega}\omega$ such that $\sum_{\omega\in\Omega}c_{\omega}=0$.  If we fix $\omega_0\in \Omega$, then the set of differences $\omega-\omega_0$ where $\omega$ runs over $\Omega\setminus \{\omega_0\}$ is a basis for $\mathrm{Aug}(K\Omega)$.  Thus $m$ annihilates $\mathrm{Aug}(\Omega)$ if and only if $\omega m=\omega_0 m$ for all $\omega\in \Omega$, i.e., $m$ is a constant map. This has a generalization, due to the author and Almeida~\cite{mortality} (inspired by Rystsov~\cite{rystrank}), that we reproduce here for the reader's convenience.  First we need some notation.  If $X\subseteq \Omega$, let $[X] =\sum_{\omega\in X}\omega$.

\begin{Prop}\label{minrankspace}
Let $(\Omega,M)$ be a transformation monoid of degree $n$ and min-rank $r$.  Let $K\Omega_r$ be the subspace of $\mathrm{Aug}(K\Omega)$ spanned by the differences $[X]-[Y]$ with $X,Y\in \min_M\nolimits(\Omega)$.  Then $K\Omega_r$ is a $KM$-submodule with $\dim_K K\Omega_r\leq n-r$.  Moreover, if $M$ is transitive, then $m\in M$ annihilates $K\Omega_r$ if and only if $m\in I(M)$.
\end{Prop}
\begin{proof}
First observe that $\varepsilon([X]-[Y]) = r-r=0$ for $X,Y\in \min_M\nolimits(\Omega)$ and so $K\Omega_r\subseteq \mathrm{Aug}(K\Omega)$.
The $M$-invariance of $K\Omega_r$ follows from the fact that $\min_M(\Omega)$ is $M$-invariant and Proposition~\ref{kernelpartition}.  Fix $s\in I(M)$ and let $\ker s=\{P_1,\ldots,P_r\}$. Proposition~\ref{kernelpartition} shows that if $X\in \min_M(\Omega)$, then $|X\cap P_i|=1$ for $i=1,\ldots, r$.  But $|X\cap P_i| = \langle [X],I_{P_i}\rangle$ and so $K\Omega_r\subseteq \mathrm{Span}\{I_{P_i}\mid 1\leq i\leq r\}^{\perp}$.  Since $\{P_1,\ldots,P_r\}$ is a partition, the indicator functions $I_{P_1},\ldots,I_{P_r}$ trivially form a linearly independent subset of $K^\Omega$.  As our pairing is non-degenerate, we may conclude that $\dim_K K\Omega_r\leq n-r$.

Suppose now that $(\Omega,M)$ is transitive.  Then if $m\in I(M)$, trivially $Xm=\Omega m$ for any $X\in \min_M(\Omega)$.  Thus $m$ annihilates $K\Omega_r$.  Suppose that $m\notin I(M)$.  Then $m$ has rank at least $r+1$.  Choose $X\in \min_M(\Omega)$.  Then $|Xm|=r$ and hence $Xm$ is a proper subset of $\Omega m$.  Let $\alpha\in \Omega m\setminus Xm$ and suppose that $\alpha =\beta m$.  By transitivity of $I(M)$ on $\Omega$, we can find $n\in I(M)$ with $\beta\in \Omega n=Y$.  Then $([Y]-[X])m = [Y]m-[X]m\neq 0$ as the coefficient of $\alpha$ in $[Y]m$ is non-zero, whereas the coefficient of $\alpha$ in $[X]m$ is zero.  Thus $m$ does not annihilate $K\Omega_r$, completing the proof.
\end{proof}

Of course, when the min-rank is $1$ and $(\Omega,M)$ is transitive, then $I(M)$ consists of the constant maps and $K\Omega_{1}=\mathrm{Aug}(K\Omega)$.  On the other extreme, if the min-rank of $(\Omega,M)$ is $n$, that is, $(\Omega,M)$ is a permutation group, then $K\Omega_n=\{0\}$.  

Our next result generalizes a result from~\cite{synchgroups} for permutation groups.  Let us continue to assume that $K$ is a field of characteristic $0$.  Then a transformation monoid $(\Omega,M)$ is said to be a \emph{$KI$-monoid} if $\mathrm{Aug}(K\Omega)$ is a simple $KM$-module.  It is well known that for permutation groups being a $\mathbb CI$-group is equivalent to $2$-transitivity and being an $\mathbb RI$-group is equivalent to $2$-homogeneity~\cite{synchcoop}.  The case of $\mathbb QI$-groups has been studied in~\cite{synchgroups,dixonQI,synchcoop}.  The results of~\cite{synchgroups} imply that a $KI$-group $(\Omega,G)$ is primitive and if $f$ is any non-invertible map on $\Omega$, then $\langle G,f\rangle$ contains a constant map.  Here is the general case.

\begin{Thm}\label{KImonoid}
Let $(\Omega,M)$ be a $KI$-monoid. Then:
\begin{enumerate}
\item $(\Omega,M)$ is primitive;
\item If in addition $(\Omega,M)$ is transitive, then either it is a permutation group or contains a constant map.
\end{enumerate}
\end{Thm}
\begin{proof}
Let $\equiv$ be a non-trivial proper congruence on $\Omega$.  Functoriality of the transformation module construction and the observation that the trivial module $K$ is the transformation module associated to the trivial action of $M$ on a one-point set yield the commutative diagram
\[\xymatrix{K\Omega\ar[rr]^{\psi}\ar[rd]_{\varepsilon}&&K[\Omega/{\equiv}]\ar[ld]^{\varepsilon'}\\ &K&}\] with $\psi$ induced by the quotient map and with $\varepsilon,\varepsilon'$ the augmentations.  As $\equiv$ is proper and non-trivial it follows that $\ker \psi$ is a non-zero proper $KM$-submodule of $\ker \varepsilon=\mathrm{Aug}(K\Omega)$, contradicting that $(\Omega,M)$ is a $KI$-monoid.  Thus $(\Omega,M)$ is primitive.

To prove the second item, assume by way of contradiction that the min-rank $r$ of $(\Omega,M)$ satisfies $1<r<n$.  Since $r<n$, $\min_M\nolimits(\Omega)$ has at least two elements and so $K\Omega_r\neq 0$.  On the other hand, Proposition~\ref{minrankspace} shows that $K\Omega_r$ is a $KM$-submodule of $\mathrm{Aug}(K\Omega)$ of dimension at most $n-r<n-1=\dim \mathrm{Aug}(K\Omega)$.  This contradicts that $(\Omega,M)$ is a $KI$-monoid.
\end{proof}

\subsection{Partial transformation modules}
Next we consider the case of $M$-sets with zero.  Even if we start with a transformation monoid $(\Omega,M)$, consideration of the quotient $\Omega/\Lambda$ by an $M$-invariant subset $\Lambda$ will lead us to this case.  So suppose that $(\Omega,M)$ is a finite transformation monoid with $\Omega$ an $M$-set with zero.  For the moment we shall denote the zero of $\Omega$ by $\zeta$ to distinguish it from the zero element of $K\Omega$.  Define the \emph{partial transformation module} (or \emph{contracted transformation module}) \[K_0\Omega = K\Omega/K\zeta.\]  This is indeed a $KM$-module because $K\zeta$ is a $KM$-submodule. As a $K$-vector space, $K_0\Omega$ has basis the cosets $\alpha+K\zeta$ with $\zeta\neq\alpha\in \Omega$.  Thus from now on we identify $\zeta$ with the zero of $K_0\Omega$ and return to using $0$ for the distinguished sink of $\Omega$.  We identify $\alpha$ with the coset $\alpha+K\zeta$ for $0\neq \alpha\in \Omega$.  Said differently, we can view $K_0\Omega$ as a $K$-vector space with basis $\Omega\setminus\{0\}$.  The action of $m\in M$ on $\Omega$ is extended linearly, but where we now identify the zero of $\Omega$ with the zero of $K\Omega_0$.

An alternate viewpoint is the following (retaining the above notation).  The augmentation \[\varepsilon\colon K\Omega\to K\] splits via the map $K\to K\Omega$ given by $c\mapsto c\zeta$.  Thus \[K\Omega=\mathrm{Aug}(K\Omega)\oplus K\zeta\cong \mathrm{Aug}(K\Omega)\oplus K\] as a $KM$-module and so $K_0\Omega=K\Omega/K\zeta\cong \mathrm{Aug}(K\Omega)$.  The natural basis to take for $\mathrm{Aug}(K\Omega)$ consists of all differences $\omega-\zeta$ with $\omega\in \Omega\setminus \{\zeta\}$.  Then the action of $m\in M$ is given by $m(\omega-\zeta) = m\omega-\zeta$ and so this provides another model of $K_0\Omega$.

If $M$ contains the zero map $z$, then $z$ acts on $K_0\Omega$ as a zero and so $K_0\Omega$ is naturally a module for the \emph{contracted monoid algebra} $K_0M=KM/Kz$. Therefore,  we will continue to use $0$ to denote the zero map on $\Omega$ and we shall identify the zero of $M$ with the zero of $K_0M$ and view $K_0\Omega$ as a $K_0M$-module.  The representations of $K_0M$ are exactly the representations of $M$ that send $0$ to the zero matrix. In particular, the trivial $KM$-module $K$ is not a $K_0M$-module and hence $K_0\Omega$ does not contain the trivial representation as a constituent.  We record this as a proposition.

\begin{Prop}
Let $(\Omega,M)$ be transformation monoid where $\Omega$ is an $M$-set with zero and suppose that $M$ contains the zero map.  Then the trivial module is not a constituent of $K_0\Omega$ and in particular $K_0\Omega^M=0$.
\end{Prop}

Notice that if $\Omega$ is an $M$-set and $\Lambda$ is an $M$-invariant subset, then there is an isomorphism $K\Omega/K\Lambda\cong K_0[\Omega/\Lambda]$.  We shall use both notations as convenient.

Returning to the case of a transformation monoid $(\Omega,M)$ where $\Omega$ is an $M$-set with zero, we would like the analogue of the dual pairing \eqref{pairing}.  Let us again momentarily use the notation $\zeta$ for the zero of $\Omega$.   Let $K^{\Omega}_0$ be the subspace of all function $f\colon \Omega\to K$ such that $f(\zeta)=0$.  This is a $KM$-submodule because $f(\zeta)=0$ implies $mf(\zeta)=f(\zeta m)=f(\zeta)=0$ for all $m\in M$.  As $K^{\Omega}_0$ is the annihilator of $K\zeta$ with respect to the pairing \eqref{pairing}, it follows that the pairing descends to a non-degenerate dual pairing $K_0\Omega\times K^{\Omega}_0\rightarrow K$ given by
\[\langle \alpha,f\rangle =f(\alpha)\] for $\alpha\in \Omega\setminus \{0\}$ which is compatible with the $KM$-module structure. Alternatively, if we identify $K_0\Omega$ with $\mathrm{Aug}(K\Omega)$, then we can just restrict the original pairing \eqref{pairing}. We now return to writing $0$ for $\zeta$ and identify $K^{\Omega}_0$ with $K^{\Omega\setminus \{0\}}$.  The left action of $m\in M$ on $f\colon K^{\Omega\setminus \{0\}}\to K$ is then given by \[mf(\alpha) = \begin{cases} f(\alpha m) & \alpha m\neq 0\\ 0 & \text{else.}\end{cases}\]  The dual basis to $\Omega\setminus \{0\}$ consists of the functions $\delta_{\alpha}$ with $\alpha\in \Omega\setminus \{0\}$.
If $M$ contains the zero map $z$, then $z$ annihilates $K^{\Omega}_0$ (viewed as a subspace of $K^{\Omega}$) and hence $K^{\Omega}_0$ is a left $K_0M$-module.

Let us return to the case of a finite transformation monoid $(\Omega,M)$ (with or without zero).  Consider a strong orbit $\OO_s(\omega)$ of $M$ on $\Omega$.  Let $\Upsilon(\omega) = \omega M\setminus \OO_s(\omega)$.  Then $\omega M$ is an $M$-invariant subset of $\Omega$ and $\Upsilon(\omega)$ is an $M$-invariant subset of $\omega M$.  Thus we can form the quotient $0$-transitive $M$-set $\omega M/\Upsilon (\omega)$ and hence the partial transformation module $K_0[\omega M/\Upsilon (\omega)]\cong K\omega M/K\Upsilon (\omega)$ (where if $\Upsilon(\omega)=\emptyset$, we interpret $\omega M/\Upsilon (\omega)=\omega M$ and $K\Upsilon(\omega)=0$).  This module has a basis in bijection with $\OO_s(\omega)$.  Thus we can put a right $KM$-module structure on $K\OO_s(\omega)$ by putting, for $\alpha\in \OO_s(\omega)$,
\[\alpha m = \begin{cases} \alpha\cdot m & \alpha\cdot m\in \OO_s(\omega)\\ 0 & \alpha\cdot m\notin \OO_s(\omega)\end{cases}\]
where for the moment we use $\cdot$ to indicate the action in $\Omega$.  With this module structure, we have a $KM$-isomorphism $K\OO_s(\omega)\cong K_0[\omega M/\Upsilon (\omega)]$.  If one considers an unrefinable series of $M$-invariant subsets of $\Omega$ as per \eqref{series}, then one obtains a series
\[K\Omega = K\Omega_0\supset K\Omega_1\supset K\Omega_2\supset\cdots \supset K\Omega_k\supset \{0\}\]
with successive quotients the modules of the form $K\OO_s(\omega)$ with $\omega\in \Omega$.  In particular, every irreducible constituent of $K\Omega$ is a constituent of some $K\OO_s(\omega)$ with $\omega\in \Omega$.

\section{A brief review of monoid representation theory}
In this section we briefly review the theory of irreducible representations of finite monoids.  This theory was first developed by Munn, Ponizovsky and Clifford~\cite[Chapter 5]{CP}.  It was further refined and elaborated on by Rhodes and Zalcstein~\cite{RhodesZalc}, Lallement and Petrich~\cite{LallePet} and McAlister~\cite{McAlisterCharacter}.  In~\cite{myirreps} a modern functorial approach was adopted based on Green's theory~\cite[Chapter 6]{Greenpoly}; more in depth information can be found in~\cite{rrbg}.   See also~\cite{ZurJohnBen} for the analogue over semirings. The advantage of this approach is that it avoids reliance on technical semigroup theory and at the same time clarifies the situation by highlighting functoriality and adjunctions.

Fix a finite monoid $M$.  If $e\in E(M)$, define $I_e = \{m\in M\mid e\notin MmM\}$ and observe that $I_e$ is an ideal of $M$.  We follow the obvious conventions when $I_e=\emptyset$, that is, $e\in E(I(M))$. Define $A_e= KM/KI_e\cong K_0[M/I_e]$.  Stability immediately yields that $I_e\cap eMe=eMe\setminus G_e$.  Thus $eA_ee\cong KG_e$.  Hence by Green's theory~\cite{Greenpoly,myirreps} there are induction, restriction and coinduction functors between $KG_e$-modules and $A_e$-modules.  Viewing the category of $A_e$-modules as a full subcategory of the category of $KM$-modules, we have the following functors:
\begin{gather*}
\Ind_e\colon \modu KG_e\to \modu KM\\ \Res_e\colon \modu KM\to \modu KG_e\\ \Coind_e\colon \modu KG_e\to \modu KM
\end{gather*}
defined by
\begin{align*}
\Ind_e(V) &= V\otimes_{KG_e}e(KM/KI_e) = V\otimes_{KG_e}K_0[eM/eI_e]\\
\Res_e(V) &= Ve\\
\Coind_e(V) &= \hom_{KG_e}((KM/KI_e)e,V) = \hom_{G_e}(Me\setminus {I_ee},V).
\end{align*}

Moreover, we have the following results~\cite{myirreps,rrbg}.

\begin{Prop}\label{repfacts}
Let $e\in E(M)$. Let $K$ be any field (not necessarily characteristic zero).
\begin{enumerate}
\item If $V$ is a $KM$-module annihilated by $I_e$ and $W$ is a $KG_e$-module, then there are natural isomorphisms:
\begin{align*}
\hom_{KM}(\Ind_e(W),V)&\cong \hom_{KG_e}(W,\Res_e (V))\\
\hom_{KM}(V,\Coind_e(W))&\cong \hom_{KG_e}(\Res_e (V),W).
\end{align*}
\item The functors $\Res_e\Ind_e$ and $\Res_e\Coind_e$ are naturally isomorphic to the identity functor on $\modu KG_e$.
\item The functors $\Ind_e$, $\Res_e$ and $\Coind_e$ are exact and preserve direct sum decompositions.  Moreover, $\Ind_e$ and $\Coind_e$ preserve indecomposability.
\end{enumerate}
\end{Prop}
\begin{proof}
We just sketch the proof.  See~\cite{myirreps,rrbg,Greenpoly} for details.  The first part follows from the classical adjunction between tensor products and hom functors once one observes that $\Res_e(V)\cong\hom_{A_e}(eA_e,V)\cong V\otimes_{A_e} Ae$.  The second part is direct from Green-Morita theory~\cite[Chapter 6]{Greenpoly}; see also~\cite{myirreps}.  Let us turn to the last part.  The point here is that $eM\setminus eI_e$ is a free left $G_e$-set and $Me\setminus {I_ee}$ is a free right $G_e$-set~\cite{qtheor,CP}.  Thus $e(KM/KI_e)$ and $(KM/KI_e)e$ are free $KG_e$-modules and so $\Ind_e$ and $\Coind_e$ are exact.  As any additive functor preserves direct sum decompositions it remains to consider indecomposability.

To see that these functors preserve indecomposability, let $V$ be a $KG_e$-module and observe that (1) and (2) yield \[\hom_{KM}(\Ind_e(V),\Ind_e(V))\cong \hom_{KG_e}(V,\Res_e\Ind_e(V))\cong \hom_{KG_e}(V,V)\] and in fact this isomorphism is a ring isomorphism.  But a module is indecomposable if and only if the only idempotents in its endomorphism algebra are $0$ and $1$.  Thus $V$ is indecomposable if and only if $\Ind_e(V)$ is indecomposable.  The argument for $\Coind_e(V)$ is identical.
\end{proof}

From the theory of Green~\cite{Greenpoly,myirreps}, if $V$ is a simple $KG_e$-module, then $\Ind_e(V)$ has a unique maximal submodule $\rad(\Ind_e(V))$ that can be described as the largest submodule annihilated by $e$, or alternatively \[\rad(\Ind_e(V))=\{v\in \Ind_e(V)\mid vme=0, \forall m\in M\}.\]  The quotient $\til V=\Ind_e(V)/\rad(\Ind_e(V))$ is then a simple $KM$-module and $\til Ve\cong V$;  in fact, the image of the projection $\Ind_e(V)\to \til V$ under the restriction functor $\Res_e$ is the identity as $e$ annihilates $\rad(\Ind_e(V))$.  It turns out that all simple $KM$-modules are constructed in this way~\cite{myirreps}.

\begin{Thm}
Let $K$ be a field and $M$ a finite monoid.  Choose a transversal of idempotents $e_1,\ldots,e_m$ to the set of principal ideals generated by idempotents.  Let $\Irr(KG_{e_i})$ contain one simple $KG_{e_i}$-module from each isomorphism class.  Then the modules of the form $\til V=\Ind_{e_i}(V)/\rad(\Ind_{e_i}(V))$ where $V\in \Irr(KG_{e_i})$ and $1\leq i\leq m$ form a complete set of representatives of the isomorphism classes of simple $KM$-modules.
\end{Thm}

Recall that if $V$ is a $KM$-module, then $\rad(V)$ is the intersection of all the maximal submodules of $V$.  The quotient $V/\rad(V)$ is a semisimple module called the \emph{top} of $V$, denoted $\mathrm{top}(V)$.

The description of the radical of $\Ind_{e_i}(V)$ for $V$ a simple $KG_{e_i}$-module generalizes.

\begin{Prop}\label{computeradical}
Let $M$ be a finite monoid and $e\in E(M)$.  Suppose that $K$ is a field of characteristic zero and
$V$ is a $KG_e$-module.  Then
\begin{equation}\label{radicaleq}
\rad(\Ind_e(V)) = \{w\in \Ind_e(V)\mid wme=0, \forall m\in M\}
\end{equation}
is the largest submodule of $\Ind_e(V)$ annihilated by $e$.
\end{Prop}
\begin{proof}
Denote by $U$ the right hand side of \eqref{radicaleq}; it is clearly the largest $KM$-submodule of $\Ind_e(V)$ annihilated by $e$.  Let $V=\bigoplus_{i=1}^sm_iV_i$ be the decomposition of $V$ into simple $KG_e$-modules.  Then as \[\Ind_e(V)\cong \bigoplus_{i=1}^sm_i\Ind_e(V_i),\] and $\til V_i=\Ind_e(V_i)/\rad(\Ind_e(V_i))$, we have an exact sequence of $KM$-modules
\[0\longrightarrow \rad(\Ind_e(V))\longrightarrow \Ind_e(V)\rightarrow \bigoplus_{i=1}^sm_i\til V_i\longrightarrow 0.\]
Using the exactness of the restriction functor $\Res_e$ and the fact that it maps the projection $\Ind_e(V_i)\to \til V_i$ to the identity map $V_i\to V_i$, we see that $0=\Res_e (\rad(\Ind_e(V)))=\rad(\Ind_e(V))e$.  This shows that $\rad(\Ind_e(V))\subseteq U$.

For the converse, let $\p\colon \Ind_e(V)\to W$ be an epimorphism of $KM$-modules with $W$ a simple $KM$-module. Then $I_e$ annihilates $W$ and so by the adjunction, we have a non-zero morphism $V\to We$ and so $We\neq 0$.  Now $\p(U)$ is a submodule of $W$.  If it is non-zero, then $\p(U)=W$.  But then $We=\p(U)e=\p(Ue)=\p(0)=0$, a contradiction.  Thus $U\subseteq \ker \p$.  As $\p$ was arbitrary, we conclude that $U\subseteq \rad(\Ind_e(V))$.
\end{proof}

A fact we shall use later is that \[\Ind_e(V)eKM=V\otimes_{KG_e}e(KM/I_e)eKM=\Ind_e(V)\] because $e(KM/I_e)e=KG_e$ and $VKG_e=V$.

\section{The projective cover of a transformation module}

From now on we assume that the characteristic of our field $K$ is zero and we fix a finite monoid $M$.
An important special case of the above theory is when $e\in E(I(M))$.  In this case $I_e=\emptyset$ and so $\Ind_e(V) = V\otimes_{KG_e}eKM$ and $\Coind_e(V) = \hom_{G_e}(Me,V)$.  Moreover, the adjunctions of Proposition~\ref{repfacts} hold for all $KM$-modules $V$.  Observe that $\Ind_e(KG_e) = KG_e\otimes_{KG_e}eKM=eKM$ is a projective $KM$-module (as $KM = eKM\oplus (1-e)KM$).  Let
\[KG_e = \bigoplus_{i=1}^s d_iV_i\] be the decomposition of $KG_e$ into simple modules.  Then the decomposition
\[eKM = \Ind_e(KG_e) = \bigoplus_{i=1}^s d_i\Ind_e(V_i)\] establishes that the $\Ind_e(V_i)$ are projective modules.  Furthermore, $\Ind_e(V_i)$ is indecomposable by Proposition~\ref{repfacts}. Thus $\Ind_e(V_i)\to \til V_i$ is the projective cover of the simple module $\til V_i$.  We recall here that if $V$ is a module over a finite dimensional algebra $A$, then the projective cover $P$ of $V$ is a projective module $P$ together with an epimorphism $\pi\colon P\to V$ such that $\pi$ induces an isomorphism $\mathrm{top}(P)\to \mathrm{top}(V)$~\cite{assem}.  Equivalently, it is an epimorphism $\pi\colon P\to V$ with $\ker\pi\subseteq \rad(P)$.  The projective cover of a module is unique up to isomorphism~\cite{assem}.  The projective covers of the simple modules are the projective indecomposables.  We have thus proved:

\begin{Prop}\label{projcover1}
Let $K$ be a field of characteristic zero and $M$ a finite monoid.  Let $e\in E(I(M))$ and assume that $V_i$ is a simple $KG_e$-module.  Then the projection $\Ind_e(V_i)\to \til V_i$ is the projective cover of the simple $KM$-module $\til V_i$.
\end{Prop}

Note that if $\Lambda$ is a right $M$-set and $\Omega$ is a bi-$M$-$N$-set, then $K[\Lambda\otimes_M\Omega]\cong K\Lambda\otimes_{KM} K\Omega$ as a right $KN$-module, as is immediate from the universal property of tensor products.  In particular, if $e\in E(I(M))$ and $\Omega$ is a $G_e$-set, then one has $K[\ind_e(\Omega)]\cong \Ind_e(K\Omega)$.  Taking $\Omega$ to be the trivial $G_e$-set $\{\ast\}$, we then have $\Ind_e(K)\cong K[\ind_e(\{\ast\})] = K(\{\ast\}\otimes_{G_e} eM) = K(G_e\backslash eM)$.  Thus Proposition~\ref{projcover1} has the following consequence.

\begin{Cor}\label{trivialcover}
The projective cover of the trivial representation of $KM$ is the augmentation map $\varepsilon\colon K[G_e\backslash eM]\to K$ where $e\in E(I(M))$.  In particular, if $(\Omega,M)$ is a transitive transformation monoid with the maximal subgroup of $I(M)$ trivial, then $K\Omega$ is a projective indecomposable representation with simple top the trivial $KM$-module and radical $\mathrm{Aug}(K\Omega)$.
\end{Cor}
\begin{proof}
It just remains to verify the final statement.  But Proposition~\ref{constantmapcase} shows that in this case $\Omega\cong eM=G_e\backslash eM$.
\end{proof}

Let $A$ be any finite dimensional $K$-algebra and $P$ a projective indecomposable with corresponding simple module $S=P/\rad(P)$.  Then it is well known that, for any $A$-module $V$, the $K$-dimension of $\hom_A(P,V)$ is the multiplicity of $S$ as an irreducible constituent of $V$~\cite{assem}.  Hence we have the reciprocity result:

\begin{Prop}\label{constituentmultiplicity}
Suppose $e\in E(I(M))$ and $V_i$ is a simple $KG_e$-module.  Let $W$ be a $KM$-module.  Then the multiplicity of $\til V_i$ as a constituent of $W$ is the same as the multiplicity of $V_i$ as a constituent of $\Res_e (W)=We$.
\end{Prop}
\begin{proof}
Since $\Ind_e(V_i)$ is the projective cover of $\til V_i$, we have that the multiplicity of $\til V_i$ in $W$ is \[\dim_K \hom_{KM}(\Ind_e(V_i),W)=\dim_K\hom_{KG_e}(V_i,\Res_e(W))\] and this latter dimension is the multiplicity of $V_i$ in $We$.
\end{proof}

The advantage of this proposition is that one can then apply the orthogonality relations of group representation theory~\cite{curtis} in order to compute the multiplicity.  Applying this to the special case of the trivial representation of $KG_e$ yields:

\begin{Cor}
Let $(\Omega,M)$ be a transformation monoid.  The multiplicity of the trivial $KM$-module as an irreducible constituent of $K\Omega$ is the number of orbits of $G_e$ on $\Omega e$ where $e\in E(I(M))$.  This can be strictly larger than \[\dim_K\hom_{KM}(K\Omega,K)=|\pi_0(\Omega)|.\]
\end{Cor}
\begin{proof}
By standard group representation theory, the multiplicity of the trivial representation of $G_e$ in $K\Omega e$ is the number of orbits of $G_e$ on $\Omega e$~\cite{curtis,cameron,dixonbook}.  The final statement follows from Proposition~\ref{fixedset} and the example just after Proposition~\ref{weakorbitrestriction}.
\end{proof}

Next we want to establish the analogues of Propositions~\ref{projcover1} and~\ref{constituentmultiplicity} for the case of monoids with zero.
\begin{Prop}\label{projcover2}
Let $M$ be a finite monoid with zero containing a unique $0$-minimal ideal $I$ and let $K$ be a field of characteristic zero.  Let $0\neq e\in E(I)$ and suppose that $V$ is a simple $KG_e$-module.  Then $\Ind_e(V)$ is a projective indecomposable $KM$-module and the projection $\Ind_e(V)\to \til V$ is the projective cover.  Moreover, if $W$ is a $K_0M$-module, then the multiplicity of $\til V$ as a constituent in $W$ is the same as the multiplicity of $V$ as a constituent in $We$.
\end{Prop}
\begin{proof}
First observe that if $z$ is the zero of $M$, then $z$ and $1-z$ are central idempotents of $KM$ and so we have an isomorphism of $K$-algebras \[KM= (1-z)KM\oplus zK\cong K_0M\oplus K.\]  Thus $K_0M$ is a projective $KM$-module.  But $K_0M=eK_0M\oplus (1-e)K_0M$ and so $eK_0M$ is a projective $KM$-module. Suppose that $KG_e=\bigoplus_{i=1}^s d_iV_i$ is the decomposition into simple $KG_e$-modules. Then \[eK_0M = KG_e\otimes_{KG_e} eK_0M = \Ind_e(KG_e) = \bigoplus_{i=1}^s d_i\Ind_e(V_i)\] and thus each $\Ind_e(V_i)$ is a projective module.  Proposition~\ref{repfacts} then yields that $\Ind_e(V_i)$ is a projective indecomposable and hence the canonical projection $\Ind_e(V_i)\to \til V_i$ is the projective cover.

For the final statement, Proposition~\ref{repfacts} provides the isomorphism \[\hom_{KM}(\Ind_e(V),W)\cong \hom_{KG_e}(V,We).\]  The dimension of the left hand side is the multiplicity of $\til V$ as a constituent of $W$, whereas the dimension of the right hand side is the multiplicity of $V$ as a constituent in $We$.
\end{proof}

An immediate corollary is the following.
\begin{Cor}
Let $(\Omega,M)$ be a $0$-transitive finite transformation monoid such that $G_e$ is trivial for $0\neq e\in E(I)$ where $I$ is the $0$-minimal ideal of $M$.  Then $K_0\Omega$ is a projective indecomposable $KM$-module.
\end{Cor}
\begin{proof}
We know that $\Omega\cong eM$ from Proposition~\ref{aperiodicbottom} and so $K_0\Omega\cong eK_0M=\Ind_e(K)$ and hence is a projective indecomposable by Proposition~\ref{projcover2}.
\end{proof}

In~\cite{mobius2} it is proved that if $(\Omega,M)$ is a $0$-transitive transformation inverse monoid, then the module $K_0\Omega$ is semisimple and decomposes as follows. Let $e$ be a non-zero idempotent of the $0$-minimal ideal of $M$ and let $\bigoplus_{i=1}^sm_iV_i$ be decomposition of $K\Omega e$ into simple $KG_e$-modules.  Then \[K_0\Omega\cong \bigoplus_{i=1}^sm_i\til V_i.\]  For more general transformation monoids, we lose semisimplicity.  But we show here that the analogous result holds at the level of the projective cover.  Of course, in characteristic zero, inverse monoid algebras are semisimple~\cite{CP} and so the simple modules are the projective indecomposables.

\subsection{The transitive case}
We describe here the projective cover of $K\Omega$ when $(\Omega,M)$ is transitive (and in slightly more generality).

\begin{Thm}\label{mainprojcover}
Let $(\Omega,M)$ be a finite transformation monoid and $K$ a field of characteristic zero.  Let $e\in E(I(M))$ and suppose that $\Omega eM=\Omega$; this happens, for instance, if $(\Omega,M)$ is transitive.  Then the natural map \[\p\colon \Ind_e(K\Omega e)\to K\Omega\] induced by the identity map on $K\Omega e$ is the projective cover.
\end{Thm}
\begin{proof}
First we observe that $\p$ is an epimorphism because \[\p(\Ind_e(K\Omega e)) = \p(\Ind_e(K\Omega e)eKM) = K\Omega eM=K\Omega.\]  It remains to show that $\ker \p\subseteq \rad(\Ind_e(K\Omega e))$.  By Proposition~\ref{computeradical} this occurs if and only if $e$ annihilates $\ker \p$.  But we have an exact sequence
\[0\longrightarrow \ker \p\longrightarrow \Ind_e(K\Omega e)\xrightarrow{\,\,\p\,\,} K\Omega\longrightarrow 0\] and hence application of $\Res_e$, which is exact, and the fact that $\Res_e(\p)=1_{K\Omega e}$ yield an exact sequence
\[0\longrightarrow (\ker \p)e\longrightarrow K\Omega e\xrightarrow{1_{K\Omega e}} K\Omega e\longrightarrow 0.\] Thus $(\ker \p)e=0$, as required.
\end{proof}

As a corollary, we have the following description of $\mathrm{top}(K\Omega)$.

\begin{Cor}\label{cortomainproj}
Under the hypotheses of Theorem~\ref{mainprojcover} one has \[\mathrm{top}(K\Omega)\cong \bigoplus_{i=1}^sm_i\til{V_i}\]
where $K\Omega e=\bigoplus_{i=1}^sm_iV_i$ is the decomposition into simple $KG_e$-modules. In particular, if $(\Omega,M)$ is transitive (and hence $(\Omega e,G_e)$ is transitive), then  $\sum_{i=1}^sm_i^2$ is the rank of the permutation group $(\Omega e,G_e)$.
\end{Cor}
\begin{proof}
The first part is clear from Theorem~\ref{mainprojcover}; the second part follows from a well-known result in permutation group theory~\cite{dixonbook,cameron}.
\end{proof}

\subsection{The $0$-transitive case}
Our next result is the analogous theorem for the $0$-transitive case.  Observe that if $(\Omega,M)$ is a $0$-transitive finite transformation monoid and $e$ is a non-zero idempotent of the $0$-minimal ideal $I$, then $K_0\Omega e$ is the permutation module associated to the permutation group $(\Omega e\setminus\{0\},G_e)$.

\begin{Thm}\label{mainproj2}
Let $(\Omega,M)$ be a finite $0$-transitive transformation monoid and $K$ be a field of characteristic $0$.  Let $e\neq 0$ be an idempotent of the $0$-minimal ideal $I$ of $M$.    Then the natural homomorphism \[\p\colon \Ind_e(K_0\Omega e)\to K_0\Omega\] induced by the identity map on $K_0\Omega e$ is the projective cover.  In particular, if $KM$ is semisimple, then $\Ind_e(K_0\Omega e)\cong K_0\Omega$.
\end{Thm}
\begin{proof}
The homomorphism $\p$ is surjective by the computation \[\p(\Ind_e(K_0\Omega e)) = \p(\Ind_e(K_0\Omega e)eKM) = K_0\Omega eM = K_0\Omega\] where the last equality uses $0$-transitivity.  To show that $\p$ is the projective cover, we must show that $\ker \p$ is contained in $\rad(\Ind_e(K_0\Omega e))$, or equivalently by Proposition~\ref{computeradical}, that $e$ annihilates $\ker \p$.  This is proved exactly as in Theorem~\ref{mainprojcover}.  Applying the exact functor $\Res_e$ to the exact sequence
\[0\longrightarrow \ker \p\longrightarrow \Ind_e(K_0\Omega e)\xrightarrow{\,\,\p\,\,} K_0\Omega\longrightarrow 0\]  and using that $\Res_e(\p)=1_{K_0\Omega e}$ we obtain the exact sequence
\[0\longrightarrow (\ker \p)e\longrightarrow K_0\Omega e\xrightarrow{1_{K_0\Omega e}} K_0\Omega e\longrightarrow 0.\]  It follows that $(\ker \p)e=0$, completing the proof.
\end{proof}

In particular, Theorem~\ref{mainproj2} has as a special case the result in~\cite{mobius2} decomposing the partial transformation module associated to a $0$-transitive transformation inverse monoid.

Of course, we have the following analogue of Corollary~\ref{cortomainproj}.

\begin{Cor}\label{cortomainproj2}
Under the same assumptions as Theorem~\ref{mainproj2} one has \[\mathrm{top}(K_0\Omega)\cong \bigoplus_{i=1}^sm_i\til{V_i}\]
where $K_0\Omega e=\bigoplus_{i=1}^sm_iV_i$ is the decomposition into simple $KG_e$-modules.  Moreover, $\sum_{i=1}^sm_i^2$ is the rank of the permutation group $(\Omega e\setminus \{0\},G_e)$.
\end{Cor}

\section{Probabilities, Markov chains and Neumann's lemma}
A partition $\{P_1,\ldots, P_r\}$ on a finite set $\Omega$ is said to be \emph{uniform} if all the blocks have the same size, i.e., $|P_1|=\cdots=|P_r|$.  Let's consider a probabilistic generalization.  Recall that a \emph{probability distribution} on $\Omega$ is a function $\mu\colon \Omega\to [0,1]$ such that $\sum_{\omega\in \Omega}\mu(\omega)=1$.  The \emph{support} $\mathrm{supp}(\mu)$ is the set of elements $\omega\in \Omega$ with $\mu(\omega)\neq 0$.  One can then view $\mu$ as a probability measure on $\Omega$ by putting \[\mu(A) = \sum_{\omega\in A}\mu(\omega)\] for a subset $A\subseteq \Omega$.  The uniform distribution $U$ on $\Omega$ is defined by $U(\omega)=1/|\Omega|$ for all $\omega\in \Omega$.  Of course $U(A)=|A|/|\Omega|$.  Thus a partition is uniform if and only if each of its blocks are equiprobable with respect to the uniform distribution. More generally, if $\mu$ is a probability distribution on $\Omega$,  we shall say that the partition $\{P_1,\ldots, P_r\}$ of $\Omega$ is \emph{$\mu$-uniform} if $\mu(P_1)=\cdots=\mu(P_r)$.

P.~Neumann in his work on synchronizing groups~\cite{pneumann} showed that if $(\Omega,M)$ is a finite transformation monoid with transitive group of units $G$, then the kernel of each element of $I(M)$ is a uniform partition.  In this section we consider a generalization of his result.  Our results can also be viewed as a generalization of a result of Friedman from~\cite{Friedman}.

We shall need to introduce a few more notions from probability theory. If $f\colon \Omega\to \mathbb R$ is a \emph{random variable} on $\Omega$, that is a real-valued function, then the \emph{expected value} of $f$ with respect to the probability distribution $\mu$ is
\[E_{\mu}(f) = \sum_{\omega\in \Omega}f(\omega)\mu(\omega).\]

A \emph{Markov chain} with state set $\Omega$ is given by a stochastic matrix \[P\colon \Omega\times \Omega\to [0,1]\] called the \emph{transition matrix} of the chain.  The adjective ``stochastic'' means that each row is a probability distribution on $\Omega$, i.e., for any fixed $\alpha\in \Omega$, one has \[\sum_{\omega\in\Omega} P(\alpha,\omega)=1.\]   Viewing probability distributions on $\Omega$ as row vectors, it follows that if $\mu$ is a probability distribution, then so is $\mu P$ where \[\mu P(\alpha) = \sum_{\omega\in \Omega}\mu(\omega)P(\omega,\alpha).\]  In particular, if $\mu$ is an initial distribution on $\Omega$, then $\mu P^k$ is the distribution at the $k^{th}$-step of the Markov chain.  A distribution $\pi$ is said to be \emph{stationary} if $\pi P=\pi$.

To a Markov chain with state set $\Omega$ and transition matrix $P$ one associates a digraph (possibly with loop edges) by declaring $(\alpha,\beta)$ to be an edge if $P(\alpha,\beta)>0$.  The Markov chain is said to be \emph{irreducible} if the associated digraph is strongly connected.  The following is a classical theorem in Markov chain theory.

\begin{Thm}\label{Markovchaintheorem}
Let $P$ be the transition matrix of an irreducible Markov chain with state set $\Omega$.  Then $P$ has a unique stationary distribution $\pi$, which moreover has support $\Omega$.   Furthermore, \[\lim_{k\to \infty}\frac{1}{k}\sum_{i=0}^{k-1}P^i = \Pi\] where $\Pi$ is the $\Omega\times \Omega$ matrix whose rows are all equal to $\pi$.
\end{Thm}

Let $(\Omega,M)$ be a finite transformation monoid and suppose that $\mu$ is a probability distribution on $M$.  Then we can define a Markov chain with state space $\Omega$ by putting
\begin{equation}\label{defineMarkovop}
P(\alpha,\beta) = \sum_{\alpha m=\beta}\mu(m);
\end{equation}
so $P(\alpha,\beta)$ is the probability that an element $m\in M$ chosen randomly according to $\mu$ takes $\alpha$ to $\beta$.  To see that $P$ is stochastic, notice that
\[\sum_{\beta\in \Omega}P(\alpha,\beta) = \sum_{\beta\in\Omega}\sum_{\alpha m=\beta}\mu(m) = \sum_{m\in M}\sum_{\beta=\alpha m}\mu(m)=\sum_{m\in M}\mu(m)=1.\]
If $N=\langle\supp(\mu)\rangle$ is transitive on $\Omega$, then $P$ is the transition matrix of an irreducible Markov chain.  Indeed, the digraph associated to $P$ is the underlying digraph of the automaton with state set $\Omega$ and input alphabet $\supp(\mu)$.

Observe that if $\nu$ is a probability distribution on $\Omega$, we can identify it with the element \[\sum_{\omega\in \Omega}\nu(\omega)\omega\in \mathbb R\Omega.\]  Similarly, we can identify $\mu$ with the element \[\sum_{m\in M}\mu(m)m\in \mathbb RM.\]  Then one easily verifies that \[\nu \mu = \sum_{\omega\in \Omega,m\in M}\nu(\omega)\mu(m)\omega m,\] whereas the coefficient of $\beta$ in $\nu P$ is \[\sum_{\omega\in \Omega}\nu(\omega)P(\omega,\beta) = \sum_{\omega\in \Omega,\omega m=\beta}\nu(\omega)\mu(m).\]  Thus under our identifications, we see that $\nu \mu=\nu P$ and hence $\nu P^k=\nu \mu^k$.

Our next result is an ergodic theorem in this context.

\begin{Thm}[Ergodic theorem]\label{ergodic}
Let $(\Omega,M)$ be a finite transformation monoid and let $\nu$ be a probability distribution on $\Omega$.
Suppose that $\mu$ is a probability distribution on $M$ such that $N=\langle\supp(\mu)\rangle$ is transitive on $\Omega$ and let $P$ be the transition matrix of the irreducible Markov chain defined in \eqref{defineMarkovop}. Denote by $\pi$ the stationary distribution of $P$.  If $f\colon \Omega\to \mathbb R$ is a random variable such that \[E_{\nu}(mf)=E_{\nu}(f)\] for all $m\in N$, then the equality \[E_{\pi}(f)=E_{\nu}(f)\] holds.
\end{Thm}
\begin{proof}
We use here the dual pairing of $\mathbb R\Omega$ and $\mathbb R^{\Omega}$.  Notice that if $\theta$ is any probability distribution on $\Omega$, then viewing $\theta\in \mathbb R\Omega$, we have \[E_{\theta}(f) = \sum_{\omega\in \Omega}f(\omega)\theta(\omega) = \langle \theta,f\rangle.\]  Also observe that if $\lambda$ is any probability distribution with support contained in $N$, then $E_{\nu}(\lambda f)=E_{\nu}(f)$ where we view $\lambda\in \mathbb RM$.  Indeed, linearity of expectation implies that
\[E_{\nu}(\lambda f) = \sum_{m\in N}\lambda(m)E_{\nu}(mf) = \sum_{m\in N}\lambda(m)E_{\nu}(f)=E_{\nu}(f).\]

A simple calculation reveals that $\nu\Pi=\pi$ and so applying Theorem~\ref{Markovchaintheorem} and the above observations (with $\lambda=\mu^i$) yields
\begin{align*}
E_{\pi}(f)  &= \langle \pi,f\rangle
            = \langle \nu\Pi,f\rangle
            = \left\langle \nu\lim_{k\to \infty}\frac{1}{k}\sum_{i=0}^{k-1}P^i,f\right\rangle\\
            &=\lim_{k\to \infty}\frac{1}{k}\sum_{i=0}^{k-1}\langle \nu \mu^i,f\rangle
            =\lim_{k\to \infty}\frac{1}{k}\sum_{i=0}^{k-1}\langle \nu,\mu^if\rangle
            \\ &=\lim_{k\to \infty}\frac{1}{k}\sum_{i=0}^{k-1}E_{\nu}(\mu^if)=\lim_{k\to \infty}\frac{1}{k}\sum_{i=0}^{k-1}E_{\nu}(f)\\
            &=E_{\nu}(f)
\end{align*}
as required.
\end{proof}

As a consequence, we obtain the following result.
\begin{Lemma}\label{generalizedneumann}
Let $(\Omega,M)$ be a finite transformation monoid and let $\mu$ be a probability distribution on $M$ such that $N=\langle \supp(\mu)\rangle$ is transitive.  Let $P$ be the stochastic matrix \eqref{defineMarkovop} and let $\pi$ be the stationary distribution of the irreducible Markov chain with transition matrix $P$.
Suppose that $B$ and $S$ are subsets of $\Omega$ such that $|S\cap Bm\inv|=1$ for all $m\in N$.  Then $|S|\cdot \pi(B)=1$.
\end{Lemma}
\begin{proof}
Observe that taking $m=1$, we have $|S\cap B|=1$.
Let $\nu$ be the probability distribution on $\Omega$ given by $I_S/|S|$.  Then, for $m\in N$, we have \[E_{\nu}(mI_B) = E_{\nu}(I_{Bm\inv}) = \nu(Bm\inv) = |S\cap Bm\inv|/|S|=1/|S|=E_{\nu}(I_B).\]  Thus the ergodic theorem yields \[1/|S|=E_{\nu}(I_B)=E_{\pi}(I_B)=\pi(B)\] and so $1=|S|\cdot \pi(B)$ as required.
\end{proof}

A particular example is the case that $(\Omega,G)$ is a transitive permutation group and $\mu$ is the uniform distribution on $G$.  One easily verifies that $\pi$ is the uniform distribution on $\Omega$ (since each element of $G$ fixes the uniform distribution on $\Omega$ as an element of $\mathbb R\Omega$).  Thus the lemma says in this setting that if $S,B$ are subsets of $\Omega$ with $|S\cap Bg|=1$ for all $g\in G$, then $|S|\cdot |B|=|\Omega|$.  This is a result of P.~Neumann.

\begin{Thm}\label{genneumann}
Let $(\Omega,M)$ be a finite transformation monoid and let $\mu$ be a probability distribution on $M$ such that $\langle \supp(\mu)\rangle$ is transitive on $\Omega$.  Let $P$ be the transition matrix of the irreducible Markov chain defined in \eqref{defineMarkovop} and let $\pi$ be the stationary distribution of $P$. Let $s\in I(M)$ and suppose that $\ker s = \{B_1,\ldots, B_r\}$.  Then $|\Omega s|\cdot \pi(B_i)=1$ for $i=1,\ldots,r$.  In particular, $\ker s$ is $\pi$-uniform.
\end{Thm}
\begin{proof}
Observe that if $m\in M$, then $\Omega ms=\Omega s$ as all elements of $I(M)$ have the same rank.  Hence if $\omega_i = B_is$ (and so $B_i=\omega_is\inv$), for $i=1,\ldots,r$, then \[\ker ms = \{\omega_1(ms)\inv,\ldots, \omega_r(ms)\inv\} = \{B_1m\inv,\ldots, B_rm\inv\}.\]  Proposition~\ref{kernelpartition} now implies that $|\Omega s\cap B_im\inv|=1$ for all $1\leq i\leq r$. As $m$ was arbitrary, Lemma~\ref{generalizedneumann} yields $|\Omega s|\cdot \pi(B_i)=1$ for $i=1,\ldots, r$.
\end{proof}

As a consequence, we obtain Neumann's lemma~\cite{pneumann}.
\begin{Cor}[Neumann's lemma]\label{neumann}
Let $(\Omega,M)$ be a finite transformation monoid with a transitive group of units.  Then $\ker m$ is a uniform partition for all $m\in I(M)$.
\end{Cor}
\begin{proof}
Let $G$ be the group of units of $M$ and let $\mu$ be the uniform distribution on $G$.  Then, as observed earlier, $\pi$ is the uniform distribution on $\Omega$.  The result is now immediate from Theorem~\ref{genneumann}.
\end{proof}

We can now present Neumann's proof~\cite{Neumann} of a result of Pin~\cite{Pincerny}; it also can be deduced from Theorem~\ref{KImonoid} since a transitive permutation group of prime degree is a $\mathbb QI$-group, cf.~\cite{synchgroups}.

\begin{Prop}\label{pinresult}
Suppose that $(\Omega,M)$ is a transformation monoid with transitive group of units and $|\Omega|$ is prime.  Then either $M$ is a group or $M$ contains a rank $1$ transformation (i.e., a constant map).
\end{Prop}
\begin{proof}
The kernel of each element of $I(M)$ is a uniform partition.  Since $|\Omega|$ is prime, it follows that either each element of $I(M)$ is a permutation or each element of $I(M)$ is a constant map.  In the former case, $I(M)=M$ is a group; in the latter case $M$ contains a rank $1$ map.
\end{proof}

Neumann's lemma can be generalized to transformation monoids containing an Eulerian subset.  Let $(\Omega,M)$ be a finite transformation monoid.  Let us say that a subset $A\subseteq M$ is \emph{Eulerian} if $\langle A\rangle$ is transitive and, for each $\omega\in \Omega$, the equality
\begin{equation}\label{Eulerian}
|A|=\sum_{a\in A} |\omega a\inv|
\end{equation}
holds.  The reason for this terminology is that if one considers the automaton with input alphabet $A$, state set $\Omega$ and transition function $\Omega\times A\to \Omega$ given by $(\omega,a)\mapsto \omega a$, then the underlying digraph of the automaton (with multiple edges allowed) contains a directed Eulerian path precisely under the assumption that $A$ is Eulerian.  Eulerian automata were considered by Kari in the context of the Road Coloring Problem and the \v{C}ern\'y conjecture~\cite{Kari}.  Notice that if $A$ consists of permutations and $\langle A\rangle$ is transitive on $\Omega$, then it is trivially Eulerian because each $|\omega a\inv|=1$.  Thus the following theorem is a generalization of Neumann's lemma.

\begin{Thm}\label{eulerneumann}
Suppose that $(\Omega,M)$ is a finite transformation monoid containing an Eulerian subset $A$. Then $\ker m$ is a uniform partition for all $m\in I(M)$.  In particular, if $|\Omega|$ is prime, then either $M$ is a group or $M$ contains a rank $1$ map.
\end{Thm}
\begin{proof}
Suppose that $|A|=k$.  Define a probability distribution $\mu$ on $M$ by putting $\mu=(1/k)I_A$. Let  $P$ be the stochastic matrix \eqref{defineMarkovop}.  The corresponding Markov chain is irreducible,  let $\pi$ be its stationary distribution.  We claim that $\pi$ is the uniform distribution on $\Omega$. Theorem~\ref{generalizedneumann} will then imply that $\ker m$ is uniform for each $m\in I(M)$.  It is well known and easy to see that the uniform distribution is stationary for a Markov chain if and only if the transition matrix $P$ is doubly stochastic, meaning that the columns of $P$ also sum to $1$.  In our case, the sum of the entries of the column of $P$ corresponding to $\omega\in \Omega$ is
\[\sum_{\alpha\in\Omega}P(\alpha,\omega)=\sum_{\alpha\in \Omega}\sum_{\alpha m=\omega}\mu(m) = \sum_{m\in M}\mu(m)\cdot |\omega m\inv| = \frac{1}{k}\sum_{a\in A}|\omega a\inv| = 1\] where we have used \eqref{Eulerian}.

The final statement is proved exactly as in Proposition~\ref{pinresult}.
\end{proof}

Theorem~\ref{eulerneumann} holds more generally for any finite transformation monoid $(\Omega,M)$ such that there is a probability distribution $\mu$ on $M$ with $\langle\supp(\mu)\rangle$ transitive and the matrix $P$ from \eqref{defineMarkovop} doubly stochastic.  It is not hard to construct transformation monoids for which this occurs that do not contain Eulerian subsets.  The corresponding class of automata was termed \emph{pseudo-Eulerian} by the author in~\cite{averaging}.

\subsection{A Burnside-type lemma}
The classical Burnside lemma (which in fact was known to Cauchy and Frobenius) says that the number of orbits of a permutation group equals the average number of fixed points.  The best we can say for transformation monoids is the following, where $\mathrm{Fix}(m)$ is the fixed-point set of $m\in M$ and $\mathrm{Stab}(\omega)$ is the stabilizer of $\omega\in \Omega$.

\begin{Lemma}
Let $(\Omega,M)$ be a finite transformation monoid.  Suppose that $\mu$ is a probability distribution of $M$ and $\pi$ is a probability distribution on $\Omega$.  Let $F$ be the random variable defined on $M$ by $F(m)=\pi(\mathrm{Fix}(m))$ and let $S$ be the random variable defined on $\Omega$ by $S(\omega)=\mu(\mathrm{Stab}(\omega))$.  Then $E_{\mu}(F)=E_{\pi}(S)$.
\end{Lemma}
\begin{proof}
This is a trivial computation:
\begin{align*}
E_{\mu}(F) &= \sum_{m\in M} \pi(\mathrm{Fix}(m))\mu(m)= \sum_{\omega m=\omega}\pi(\omega)\mu(m)= \sum_{\omega\in\Omega} \mu(\mathrm{Stab}(\omega))\pi(\omega)\\
&= E_{\pi}(S)
\end{align*}
as required.
\end{proof}

The classical Burnside lemma is obtained by taking $M$ to be a group $G$, $\mu$ to be the uniform distribution on $G$ and $\pi$ to be the uniform distribution on $\Omega$:  one simply observes that $\mu(\mathrm{Stab}(\omega))=|\mathrm{Stab}(\omega)|/|G|=1/|\omega\cdot G|$.

Suppose that $|\Omega|=n$, $M=T_{\Omega}$ and one takes $\mu$ and $\pi$ to be uniform.  Clearly
$|\mathrm{Stab}(\omega)| = n^{n-1}$. Thus we have the well known result
\[\frac{1}{|T_{\Omega}|}\sum_{f\in T_{\Omega}}|\mathrm{Fix}(f)| = \frac{1}{n^n}\sum_{\omega\in \Omega}|\mathrm{Stab}(\omega)|=1\] just as in the case of the symmetric group $S_\Omega$.

\subsection*{Acknowledgments}
In the summer of 2008 I visited with Jorge Almeida in Porto, which led to our joint work~\cite{mortality}.  This paper grew out of initial discussions we had at that time.  

\def\malce{\mathbin{\hbox{$\bigcirc$\rlap{\kern-7.75pt\raise0,50pt\hbox{${\tt
  m}$}}}}}\def\cprime{$'$} \def\cprime{$'$} \def\cprime{$'$} \def\cprime{$'$}
  \def\cprime{$'$} \def\cprime{$'$} \def\cprime{$'$} \def\cprime{$'$}
  \def\cprime{$'$}


\begin{thebibliography}{10}

\bibitem{AMSV}
J.~Almeida, S.~Margolis, B.~Steinberg, and M.~Volkov.
\newblock Representation theory of finite semigroups, semigroup radicals and
  formal language theory.
\newblock {\em Trans. Amer. Math. Soc.}, 361(3):1429--1461, 2009.

\bibitem{mortality}
J.~Almeida and B.~Steinberg.
\newblock Matrix mortality and the \v {C}ern{\'y}-{P}in conjecture.
\newblock In {\em Developments in language theory}, volume 5583 of {\em Lecture
  Notes in Comput. Sci.}, pages 67--80. Springer, Berlin, 2009.

\bibitem{volkovc2}
D.~S. Ananichev and M.~V. Volkov.
\newblock Some results on \v {C}ern{\'y} type problems for transformation
  semigroups.
\newblock In {\em Semigroups and languages}, pages 23--42. World Sci. Publ.,
  River Edge, NJ, 2004.

\bibitem{volkovc3}
D.~S. Ananichev and M.~V. Volkov.
\newblock Synchronizing generalized monotonic automata.
\newblock {\em Theoret. Comput. Sci.}, 330(1):3--13, 2005.

\bibitem{volkovc1}
D.~S. Ananichev, M.~V. Volkov, and Y.~I. Zaks.
\newblock Synchronizing automata with a letter of deficiency 2.
\newblock {\em Theoret. Comput. Sci.}, 376(1-2):30--41, 2007.

\bibitem{synchcoop}
J.~Ara\'ujo, P.~J. Cameron, P.~M. Neumann, C.~Praeger, J.~Saxl, C.~Schneider,
  P.~Spiga, and B.~Steinberg.
\newblock Permutation groups and synchronizing automata.
\newblock Unpublished manuscript, October 2008.

\bibitem{synchgroups}
F.~Arnold and B.~Steinberg.
\newblock Synchronizing groups and automata.
\newblock {\em Theoret. Comput. Sci.}, 359(1-3):101--110, 2006.

\bibitem{assem}
I.~Assem, D.~Simson, and A.~Skowro{\'n}ski.
\newblock {\em Elements of the representation theory of associative algebras.
  {V}ol. 1}, volume~65 of {\em London Mathematical Society Student Texts}.
\newblock Cambridge University Press, Cambridge, 2006.
\newblock Techniques of representation theory.

\bibitem{beal}
M.-P. B{\'e}al.
\newblock A note on {C}erny's conjecture and rational series.
\newblock Unpublished, 2003.

\bibitem{PerrinBeal}
M.-P. B{\'e}al and D.~Perrin.
\newblock A quadratic upper bound on the size of a synchronizing word in
  one-cluster automata.
\newblock In {\em Developments in language theory}, volume 5583 of {\em Lecture
  Notes in Comput. Sci.}, pages 81--90. Springer, Berlin, 2009.

\bibitem{berstelperrinreutenauer}
J.~Berstel, D.~Perrin, and C.~Reutenauer.
\newblock {\em Codes and automata}, volume 129 of {\em Encyclopedia of
  Mathematics and its Applications}.
\newblock Cambridge University Press, Cambridge, 2010.

\bibitem{BHR}
P.~Bidigare, P.~Hanlon, and D.~Rockmore.
\newblock A combinatorial description of the spectrum for the {T}setlin library
  and its generalization to hyperplane arrangements.
\newblock {\em Duke Math. J.}, 99(1):135--174, 1999.

\bibitem{bjorner2}
A.~Bj{\"o}rner.
\newblock Random walks, arrangements, cell complexes, greedoids, and
  self-organizing libraries.
\newblock In {\em Building bridges}, volume~19 of {\em Bolyai Soc. Math.
  Stud.}, pages 165--203. Springer, Berlin, 2008.

\bibitem{bjorner1}
A.~Bj{\"o}rner.
\newblock Note: {R}andom-to-front shuffles on trees.
\newblock {\em Electron. Commun. Probab.}, 14:36--41, 2009.

\bibitem{Brown1}
K.~S. Brown.
\newblock Semigroups, rings, and {M}arkov chains.
\newblock {\em J. Theoret. Probab.}, 13(3):871--938, 2000.

\bibitem{Brown2}
K.~S. Brown.
\newblock Semigroup and ring theoretical methods in probability.
\newblock In {\em Representations of finite dimensional algebras and related
  topics in Lie theory and geometry}, volume~40 of {\em Fields Inst. Commun.},
  pages 3--26. Amer. Math. Soc., Providence, RI, 2004.

\bibitem{DiaconisBrown1}
K.~S. Brown and P.~Diaconis.
\newblock Random walks and hyperplane arrangements.
\newblock {\em Ann. Probab.}, 26(4):1813--1854, 1998.

\bibitem{cameron}
P.~J. Cameron.
\newblock {\em Permutation groups}, volume~45 of {\em London Mathematical
  Society Student Texts}.
\newblock Cambridge University Press, Cambridge, 1999.

\bibitem{strongtrans}
A.~Carpi and F.~d'Alessandro.
\newblock The synchronization problem for strongly transitive automata.
\newblock In {\em Developments in language theory}, volume 5257 of {\em Lecture
  Notes in Comput. Sci.}, pages 240--251. Springer, Berlin, 2008.

\bibitem{strongtrans2}
A.~Carpi and F.~d'Alessandro.
\newblock The synchronization problem for locally strongly transitive automata.
\newblock In {\em Mathematical Foundations of Computer Science}, volume 5734 of
  {\em Lecture Notes in Comput. Sci.}, pages 211--222. Springer, Berlin, 2009.

\bibitem{CP}
A.~H. Clifford and G.~B. Preston.
\newblock {\em The algebraic theory of semigroups. {V}ol. {I}}.
\newblock Mathematical Surveys, No. 7. American Mathematical Society,
  Providence, R.I., 1961.

\bibitem{CP2}
A.~H. Clifford and G.~B. Preston.
\newblock {\em The algebraic theory of semigroups. {V}ol. {II}}.
\newblock Mathematical Surveys, No. 7. American Mathematical Society,
  Providence, R.I., 1967.

\bibitem{curtis}
C.~W. Curtis and I.~Reiner.
\newblock {\em Representation theory of finite groups and associative
  algebras}.
\newblock Wiley Classics Library. John Wiley \& Sons Inc., New York, 1988.
\newblock Reprint of the 1962 original, A Wiley-Interscience Publication.

\bibitem{dixonQI}
J.~D. Dixon.
\newblock Permutation representations and rational irreducibility.
\newblock {\em Bull. Austral. Math. Soc.}, 71(3):493--503, 2005.

\bibitem{dixonbook}
J.~D. Dixon and B.~Mortimer.
\newblock {\em Permutation groups}, volume 163 of {\em Graduate Texts in
  Mathematics}.
\newblock Springer-Verlag, New York, 1996.

\bibitem{dubuc}
L.~Dubuc.
\newblock Sur les automates circulaires et la conjecture de \v {C}ern{\'y}.
\newblock {\em RAIRO Inform. Th{\'e}or. Appl.}, 32(1-3):21--34, 1998.

\bibitem{Eilenberg}
S.~Eilenberg.
\newblock {\em Automata, languages, and machines. {V}ol. {B}}.
\newblock Academic Press, New York, 1976.
\newblock With two chapters (``Depth decomposition theorem'' and ``Complexity
  of semigroups and morphisms'') by Bret Tilson, Pure and Applied Mathematics,
  Vol. 59.

\bibitem{Friedman}
J.~Friedman.
\newblock On the road coloring problem.
\newblock {\em Proc. Amer. Math. Soc.}, 110(4):1133--1135, 1990.

\bibitem{GM}
O.~Ganyushkin and V.~Mazorchuk.
\newblock {\em Classical finite transformation semigroups, an introduction}.
\newblock Number~9 in Algebra and Applications. Springer, 2009.

\bibitem{myirreps}
O.~Ganyushkin, V.~Mazorchuk, and B.~Steinberg.
\newblock On the irreducible representations of a finite semigroup.
\newblock {\em Proc. Amer. Math. Soc.}, 137(11):3585--3592, 2009.

\bibitem{Green}
J.~A. Green.
\newblock On the structure of semigroups.
\newblock {\em Ann. of Math. (2)}, 54:163--172, 1951.

\bibitem{Greenpoly}
J.~A. Green.
\newblock {\em Polynomial representations of {${\rm GL}\sb{n}$}}, volume 830 of
  {\em Lecture Notes in Mathematics}.
\newblock Springer-Verlag, Berlin, 1980.

\bibitem{Higginsbook}
P.~M. Higgins.
\newblock {\em Techniques of semigroup theory}.
\newblock Oxford Science Publications. The Clarendon Press Oxford University
  Press, New York, 1992.
\newblock With a foreword by G. B. Preston.

\bibitem{Howie}
J.~M. Howie.
\newblock The subsemigroup generated by the idempotents of a full
  transformation semigroup.
\newblock {\em J. London Math. Soc.}, 41:707--716, 1966.

\bibitem{howiebook}
J.~M. Howie.
\newblock {\em Fundamentals of semigroup theory}, volume~12 of {\em London
  Mathematical Society Monographs. New Series}.
\newblock The Clarendon Press Oxford University Press, New York, 1995.
\newblock Oxford Science Publications.

\bibitem{ZurJohnBen}
Z.~Izhakian, J.~Rhodes, and B.~Steinberg.
\newblock Representation theory of finite semigroups over semirings, April
  2010.
\newblock http://arxiv.org/abs/1004.1660.

\bibitem{Karicounter}
J.~Kari.
\newblock A counter example to a conjecture concerning synchronizing words in
  finite automata.
\newblock {\em Bull. Eur. Assoc. Theor. Comput. Sci. EATCS}, (73):146, 2001.

\bibitem{Kari}
J.~Kari.
\newblock Synchronizing finite automata on {E}ulerian digraphs.
\newblock {\em Theoret. Comput. Sci.}, 295(1-3):223--232, 2003.
\newblock Mathematical foundations of computer science (Mari{\'a}nsk{\'e}
  L{\'a}zn\v e, 2001).

\bibitem{actsbook}
M.~Kilp, U.~Knauer, and A.~V. Mikhalev.
\newblock {\em Monoids, acts and categories}, volume~29 of {\em de Gruyter
  Expositions in Mathematics}.
\newblock Walter de Gruyter \& Co., Berlin, 2000.
\newblock With applications to wreath products and graphs, A handbook for
  students and researchers.

\bibitem{PDT}
K.~Krohn and J.~Rhodes.
\newblock Algebraic theory of machines. {I}. {P}rime decomposition theorem for
  finite semigroups and machines.
\newblock {\em Trans. Amer. Math. Soc.}, 116:450--464, 1965.

\bibitem{KRannals}
K.~Krohn and J.~Rhodes.
\newblock Complexity of finite semigroups.
\newblock {\em Ann. of Math. (2)}, 88:128--160, 1968.

\bibitem{Arbib}
K.~Krohn, J.~Rhodes, and B.~Tilson.
\newblock {\em Algebraic theory of machines, languages, and semigroups}.
\newblock Edited by Michael A. Arbib. With a major contribution by Kenneth
  Krohn and John L. Rhodes. Academic Press, New York, 1968.
\newblock Chapters 1, 5--9.

\bibitem{LallePet}
G.~Lallement and M.~Petrich.
\newblock Irreducible matrix representations of finite semigroups.
\newblock {\em Trans. Amer. Math. Soc.}, 139:393--412, 1969.

\bibitem{Lawson}
M.~V. Lawson.
\newblock {\em Inverse semigroups}.
\newblock World Scientific Publishing Co. Inc., River Edge, NJ, 1998.
\newblock The theory of partial symmetries.

\bibitem{Lipscomb}
S.~Lipscomb.
\newblock {\em Symmetric inverse semigroups}, volume~46 of {\em Mathematical
  Surveys and Monographs}.
\newblock American Mathematical Society, Providence, RI, 1996.

\bibitem{Mac-CWM}
S.~{Mac Lane}.
\newblock {\em Categories for the working mathematician}, volume~5 of {\em
  Graduate Texts in Mathematics}.
\newblock Springer-Verlag, New York, second edition, 1998.

\bibitem{rrbg}
S.~W. Margolis and B.~Steinberg.
\newblock The quiver of an algebra associated to the {M}antaci-{R}eutenauer
  descent algebra and the homology of regular semigroups.
\newblock {\em Algebr. Represent. Theory}, to appear.

\bibitem{McAlisterCharacter}
D.~B. McAlister.
\newblock Characters of finite semigroups.
\newblock {\em J. Algebra}, 22:183--200, 1972.

\bibitem{selfsimilar}
V.~Nekrashevych.
\newblock {\em Self-similar groups}, volume 117 of {\em Mathematical Surveys
  and Monographs}.
\newblock American Mathematical Society, Providence, RI, 2005.

\bibitem{Neumann}
H.~Neumann.
\newblock {\em Varieties of groups}.
\newblock Springer-Verlag New York, Inc., New York, 1967.

\bibitem{pneumann}
P.~M. Neumann.
\newblock Primitive permutation groups and their section-regular partitions.
\newblock {\em Michigan Math. J.}, 58(1):309--322, 2009.

\bibitem{Pincerny}
J.-E. Pin.
\newblock Sur un cas particulier de la conjecture de {C}erny.
\newblock In {\em Automata, languages and programming (Fifth Internat. Colloq.,
  Udine, 1978)}, volume~62 of {\em Lecture Notes in Comput. Sci.}, pages
  345--352. Springer, Berlin, 1978.

\bibitem{pincernyconjecture}
J.-E. Pin.
\newblock Le probl{\`e}me de la synchronisation et la conjecture de \v
  {C}ern{\'y}.
\newblock In {\em Noncommutative structures in algebra and geometric
  combinatorics (Naples, 1978)}, volume 109 of {\em Quad. ``Ricerca Sci.''},
  pages 37--48. CNR, Rome, 1981.

\bibitem{twocomb}
J.-E. Pin.
\newblock On two combinatorial problems arising from automata theory.
\newblock In {\em Combinatorial mathematics (Marseille-Luminy, 1981)},
  volume~75 of {\em North-Holland Math. Stud.}, pages 535--548. North-Holland,
  Amsterdam, 1983.

\bibitem{Rees}
D.~Rees.
\newblock On semi-groups.
\newblock {\em Proc. Cambridge Philos. Soc.}, 36:387--400, 1940.

\bibitem{qtheor}
J.~Rhodes and B.~Steinberg.
\newblock {\em The {$q$}-theory of finite semigroups}.
\newblock Springer Monographs in Mathematics. Springer, New York, 2009.

\bibitem{RhodesZalc}
J.~Rhodes and Y.~Zalcstein.
\newblock Elementary representation and character theory of finite semigroups
  and its application.
\newblock In {\em Monoids and semigroups with applications (Berkeley, CA,
  1989)}, pages 334--367. World Sci. Publ., River Edge, NJ, 1991.

\bibitem{rystcom}
I.~Rystsov.
\newblock Reset words for commutative and solvable automata.
\newblock {\em Theoret. Comput. Sci.}, 172(1-2):273--279, 1997.

\bibitem{rystrank}
I.~C. Rystsov.
\newblock On the rank of a finite automaton.
\newblock {\em Kibernet. Sistem. Anal.}, (3):3--10, 187, 1992.

\bibitem{rystsov1}
I.~K. Rystsov.
\newblock Quasioptimal bound for the length of reset words for regular
  automata.
\newblock {\em Acta Cybernet.}, 12(2):145--152, 1995.

\bibitem{rystsov2}
I.~K. Rystsov.
\newblock On the length of reset words for automata with simple idempotents.
\newblock {\em Kibernet. Sistem. Anal.}, (3):32--39, 187, 2000.

\bibitem{Salomcerny}
A.~Salomaa.
\newblock Composition sequences for functions over a finite domain.
\newblock {\em Theoret. Comput. Sci.}, 292(1):263--281, 2003.
\newblock Selected papers in honor of Jean Berstel.

\bibitem{Schutzrep}
M.~P. Sch{\"u}tzenberger.
\newblock {$\overline{\mathscr D}$} repr{\'e}sentation des demi-groupes.
\newblock {\em C. R. Acad. Sci. Paris}, 244:1994--1996, 1957.

\bibitem{orbitoids}
R.~Scozzafava.
\newblock Graphs and finite transformation semigroups.
\newblock {\em Discrete Math.}, 5:87--99, 1973.

\bibitem{mobius1}
B.~Steinberg.
\newblock M{\"o}bius functions and semigroup representation theory.
\newblock {\em J. Combin. Theory Ser. A}, 113(5):866--881, 2006.

\bibitem{mobius2}
B.~Steinberg.
\newblock M{\"o}bius functions and semigroup representation theory. {II}.
  {C}haracter formulas and multiplicities.
\newblock {\em Adv. Math.}, 217(4):1521--1557, 2008.

\bibitem{averaging}
B.~Steinberg.
\newblock The averaging trick and the {C}erny conjecture, October 2009.
\newblock http://arxiv.org/abs/0910.0410.

\bibitem{mycerny}
B.~Steinberg.
\newblock \v{C}ern{\'y}'s conjecture and group representation theory.
\newblock {\em J. Algebr. Comb.}, 31(1):83--109, 2010.

\bibitem{Talwar3}
S.~Talwar.
\newblock Morita equivalence for semigroups.
\newblock {\em J. Austral. Math. Soc. Ser. A}, 59(1):81--111, 1995.

\bibitem{Tilson}
B.~Tilson.
\newblock Categories as algebra: an essential ingredient in the theory of
  monoids.
\newblock {\em J. Pure Appl. Algebra}, 48(1-2):83--198, 1987.

\bibitem{traht2}
A.~N. Trahtman.
\newblock An efficient algorithm finds noticeable trends and examples
  concerning the \v {C}erny conjecture.
\newblock In {\em Mathematical foundations of computer science 2006}, volume
  4162 of {\em Lecture Notes in Comput. Sci.}, pages 789--800. Springer,
  Berlin, 2006.

\bibitem{Trahtman}
A.~N. Trahtman.
\newblock The \v {C}ern{\'y} conjecture for aperiodic automata.
\newblock {\em Discrete Math. Theor. Comput. Sci.}, 9(2):3--10 (electronic),
  2007.

\bibitem{cerny}
J.~\v{C}ern{\'y}.
\newblock A remark on homogeneous experiments with finite automata.
\newblock {\em Mat.-Fyz. \v Casopis Sloven. Akad. Vied}, 14:208--216, 1964.

\bibitem{VolkovLata}
M.~V. Volkov.
\newblock Synchronizing automata and the \v{C}ern{\'y} conjecture.
\newblock In C.~Mart{\'i}n-Vide, F.~Otto, and H.~Fernau, editors, {\em Language
  and Automata Theory and Applications Second International Conference, LATA
  2008, Tarragona, Spain, March 13-19, 2008.}, volume 5196 of {\em Lecture
  Notes in Computer Science}, pages 11--27, Berlin / Heidelberg, 2008.
  Springer.

\end{thebibliography}
\end{document}